\documentclass[12pt]{amsart}
\newcounter{dummy}
\usepackage{amsmath,amsthm,amscd,amsfonts,wasysym,amssymb,epic,eepic,tikz-cd,mathtools,multicol,enumitem,mathrsfs,comment}
\usepackage[pagebackref,colorlinks=true,linkcolor=blue,citecolor=blue]{hyperref}
\usepackage[noabbrev]{cleveref}
\makeatletter
\newcommand\myitem[1][]{\item[#1]\refstepcounter{dummy}\def\@currentlabel{#1}}
\makeatother



\allowdisplaybreaks
\usepackage[margin=1in]{geometry}

\numberwithin{equation}{section}
\theoremstyle{plain}
\newtheorem{theorem}{Theorem}[section]
\newtheorem{lemma}[theorem]{Lemma}
\newtheorem{proposition}[theorem]{Proposition}
\newtheorem{corollary}[theorem]{Corollary}

\newtheorem{introtheorem}{Theorem}

\theoremstyle{definition}
\newtheorem{definition}[theorem]{Definition}
\newtheorem{example}[theorem]{Example}
\theoremstyle{remark}
\newtheorem{remark}[theorem]{Remark}
\renewcommand{\bar}{\overline}

\newcommand{\Q}{\mathbb Q}
\newcommand{\A}{\mathbb A}
\newcommand{\C}{\mathbb C}
\newcommand{\R}{\mathbb R}
\newcommand{\GL}{\mathrm{GL}}
\newcommand{\SL}{\mathrm{SL}}
\newcommand{\Sp}{\mathrm{Sp}}
\newcommand{\GSp}{\mathrm{GSp}}
\newcommand{\SO}{\mathrm{SO}}
\newcommand{\PGL}{\mathrm{PGL}}
\newcommand{\Sym}{\operatorname{Sym}}
\newcommand{\Z}{\mathbb{Z}}

\newcommand{\Std}{\mathrm{Std}}
\DeclareMathOperator{\Inv}{Inv}

\newcommand{\cM}{\mathcal M}
\newcommand{\tcM}{\widetilde{\mathcal M}}

\newcommand{\sgn}{\mathrm{sgn}}
\newcommand{\Ad}{\mathrm{Ad}}
\newcommand{\Def}{\mathrm{Def}}

\title[A Root-Theoretic Combinatorial Lemma]{Eisenstein Cohomology and Critical Values of Certain \(L\)-Functions: A Root-Theoretic Combinatorial Lemma}
\author{Farid HosseiniJafari}
\address{Department of Mathematics, Duke University, Durham, North Carolina, USA}
\email{farid.hj@duke.edu}
\date{}

\hypersetup{pdftitle={Criticality and the inverse Kostant problem for maximal parabolic subgroups},pdfauthor={Farid HosseiniJafari}}
\usepackage{booktabs,array}
\usepackage[T1]{fontenc}
\usepackage{lmodern}
\usepackage[english]{babel}

\begin{document}
\begin{abstract}
This paper gives a root-theoretic extension of the Harder--Raghuram combinatorial lemma,  relating criticality of adjacent $L-$values to degrees in nilradical cohomology, to maximal parabolics $P=MU$ of connected split reductive groups with irreducible, essentially tempered cohomological inducing data.

At real places, the paper classifies the maximal parabolics for which the two Kostant degrees in the boundary comparison can agree at a regular shift. For each such parabolic, an admissible shift makes both adjacent arguments critical in every Langlands--Shahidi grade exactly when both Kostant degrees equal $\dim U/2$. The classification recovers the two-block parabolics of type $A$ with even product of block sizes (Harder--Raghuram, proved by Weselmann) and the isotropic-line parabolic of type $D_{2m+1}$ for $m\ge2$ (Bhagwat--Raghuram). The only exceptional cases are $P_1$ and $P_6$ of type $E_6$. In all other cases, roots of $U$ fixed by the longest Levi Weyl element confine simultaneous criticality to a singular wall, generalizing the odd-block obstruction for $\GL_N$. Their number is the difference of the fundamental ranks of $M$ and $G$. When the Kostant degrees differ, as in Harder's $\GL_2$ setting, the paper determines which boundary degrees can still coincide.

At a complex place, an admissible shift is simultaneously critical exactly when the two conjugate Levi weights after the shift occur in Kostant degrees summing to $\dim U$. With a common relative center, total degree balance at a regular admissible shift characterizes simultaneous criticality at every archimedean place. For split $E_6$ over $\mathbb Q$, a half-spin family yields Eisenstein cohomology classes in degrees $20$ and $21$, including a degree-$20$ class with nonzero restrictions to both associate boundary face types.
\end{abstract}
\maketitle

\section{Introduction}\label{sec:introduction}
This paper determines when Kostant degrees in boundary cohomology characterize adjacent critical arguments of Langlands--Shahidi $L$-functions \cite{Shahidi2010}. The question begins with the \emph{inverse Kostant problem} in Harder's Eisenstein cohomology program \cite{HarderBook}. Given a Levi coefficient, does it occur in nilradical cohomology, and if so, for which dominant ambient weight and in what degree? The combinatorial lemma of Harder and Raghuram gives this question arithmetic significance. In their Rankin--Selberg setting, the shifts that realize the required Kostant degrees are precisely those that give two adjacent critical arguments, no more and no less \cite[Lemma~7.14]{HarderRaghuram2020}.

The graded formulation of the combinatorial lemma grew out of the author's work on split $G_2$ over totally imaginary fields. Its two maximal parabolics contribute two and three $L$-factors, respectively, to the Eisenstein constant term. A shift allowed by the adjoint-cube factor can be excluded by the determinant factor, which is evaluated at twice the argument. The combinatorial lemmas for $G_2$ \cite[Lemmas~6.1--6.2]{HosseiniJafariG2} thus suggest a general question. For a maximal parabolic of a connected split reductive group, when is equality of the two associate Kostant degrees equivalent to adjacent criticality in every grade? This paper proves the equivalence at complex places and classifies the obstruction at real places.

At a real place, the paper classifies the maximal parabolics for which some integral Levi weight, self-dual up to a character, admits a regular shift in middle Kostant degree. Harder and Raghuram conjectured the combinatorial lemma over $\mathbb Q$ for $n$ even and $n'$ odd \cite[Conj.~5.1]{HarderRaghuram2011}. Their monograph states the lemma over totally real fields for $nn'$ even and gives Weselmann's proof \cite[Lemma~7.14 and \S7.2.4]{HarderRaghuram2020}. Bhagwat and Raghuram treated the orthogonal line parabolic \cite[\S5.3]{BhagwatRaghuram2025} and identified the obstruction for two odd blocks of $\GL_N$ \cite[\S7.4.1]{HarderRaghuram2020}. The only exceptional cases are $P_1$ and $P_6$ of type $E_6$. For every other maximal parabolic, fixed inducing data admit at most one simultaneously critical shift. At that shift, if it exists, no dominant ambient weight realizes the shifted Levi coefficient in nilradical cohomology.

In the obstructed cases, the paper continues Harder's study of unequal Kostant lengths \cite[\S9.3.1, pp.~421--422]{HarderChapterNine}. A difference in the degrees of the Levi classes can compensate for unequal Kostant lengths, and the paper determines exactly which boundary degrees the two contributions can share. Combined with the archimedean pole calculation, this gives two arithmetic questions a precise cohomological setting. In the $\GL_4$ example, degree~$4$ is the common boundary degree in which to seek an interpretation of the first derivative of the finite $L$-function at its simple zero at $s=-1$. In the $\GSp_4$ Klingen example, the calculation singles out the critical value whose period normalization enters the author's conjecture on the denominator of an Eisenstein class \cite{HosseiniJafariKlingen}. Both examples are discussed in Section~\ref{subsec:arithmetic-unequal-degrees}.

\subsection{Setting}

Let $P=MU$ be a maximal standard parabolic of a connected split reductive group $G$ over a number field $F$. At an embedding $\eta:F\hookrightarrow\C$, let $\mathfrak u$ be the complex Lie algebra of $U$, and write $\cM_\lambda$ and $\cM_\mu^M$ for the irreducible algebraic representations of $G$ and $M$ with highest weights $\lambda$ and $\mu$. By Kostant's theorem, $\cM_\mu^M$ occurs in $H^\ell(\mathfrak u,\cM_\lambda)$ precisely when
\[
  \mu=w\cdot\lambda:=w(\lambda+\rho)-\rho,\qquad \ell=\ell(w),
\]
for a Kostant representative $w$, where $\rho$ is the half-sum of the positive roots. On the boundary face attached to $P$, a Levi class of degree $k$ contributes to degree $k+\sum_\eta\ell(w^\eta)$. Thus two Levi classes of the same degree on the associate sides can be compared in a common boundary degree when their total Kostant lengths agree. Harder varies the relative central character of the Levi coefficient to find such shifts
\cite[\S9.1.4, (9.88)--(9.90)]{HarderChapterNine}.

The link with $L$-values comes from the constant term. For a holomorphic Eisenstein class, boundary restriction involves the standard intertwining operator \cite{Schwermer1983,HarderBook}. For cuspidal inducing data $\sigma$, its normalization contains the quotient
\[
  \prod_j \frac{L^S(js,\sigma,r_j)}{L^S(js+1,\sigma,r_j)},
\]
where $S$ contains the archimedean places and the finite places at which $\sigma$ is ramified, and $r_j$ is the action of the dual Levi on grade $j$ of the dual nilradical \cite{Langlands1971,Shahidi2010}. Harder and Raghuram obtain rationality up to periods by comparing the complex operator with a map between rational boundary summands \cite{HarderRaghuram2020}. The combinatorial lemma supplies the coefficients and degrees needed to begin that comparison.

In the local results, the inducing representation is fixed, irreducible, essentially tempered, and cohomological for integral Levi coefficients. An admissible shift determines an integral Levi coefficient at each embedding and a point of evaluation $s$ (Definition~\ref{def:admissible-shift}). The shift is \emph{simultaneously critical} if, for every Langlands--Shahidi factor of grade $j$, both $js$ and $js+1$ are critical (Definition~\ref{def:simultaneous}). Here criticality at an argument $t$ means that the archimedean factor at $t$ and its dual at $1-t$ are both finite, following Deligne's criterion \cite[Def.~1.3]{Deligne1979}.

\subsection{Main theorems}

Throughout, the inducing data satisfy the hypotheses of the preceding subsection. Call a root of $U$ \emph{fixed} if the longest element of the Levi Weyl group fixes it. A shift is \emph{regular} if the shifted
Levi weight plus $\rho$ lies off every root hyperplane (Corollary~\ref{cor:shift-walls}).

\begin{introtheorem}[Combinatorial lemma at a real place]
\label{thm:intro-real}
At a real place there are two cases.
\begin{enumerate}
\item If $U$ has no fixed root, an admissible shift is simultaneously critical if and only if the shifted Levi coefficient occurs in $H^{\dim U/2}(\mathfrak u,\cM_\lambda)$ for some dominant weight $\lambda$ of $G$. The Kostant representative then has exactly half the roots of each grade as inversions. The maximal parabolics with no fixed root are precisely those in the following table. The type is that of the irreducible component containing the omitted simple root, with Bourbaki numbering.
\[
\begin{array}{c|l}
\text{Type} & \text{Maximal parabolic} \\ \hline
A_{N-1} & P_k,\quad 1\le k<N,\quad k(N-k)\text{ even} \\
D_{2m+1},\ m\ge2 & P_1,\quad\text{the stabilizer of an isotropic line} \\
E_6 & P_1,\ P_6
\end{array}
\]
\item If $U$ has a fixed root, the given inducing data admit at most one simultaneously critical shift. At any such shift, no dominant weight of $G$ realizes the shifted Levi coefficient in nilradical cohomology.
\end{enumerate}
\end{introtheorem}

Theorem~\ref{thm:real-criterion} gives the criticality criterion, including the Weil-sign test at the candidate shift in part~(2). Theorem~\ref{thm:realclass-classification} gives the classification in part~(1). The number of fixed roots is the fundamental rank of $M$ minus that of $G$ (Theorem~\ref{thm:fixed-roots-fundamental-rank}).

\begin{introtheorem}[Combinatorial lemma at a complex place]
\label{thm:intro-complex}
At a complex place, an admissible shift is simultaneously critical if and only if the shifted Levi coefficients at the two conjugate embeddings occur in nilradical cohomology for dominant ambient weights at their respective embeddings, in degrees whose sum is $\dim U$. The two Kostant representatives are then unique and complementary. This equality holds grade by grade: in grade $j$ the two inversion counts sum to $\dim r_j$, the number of roots in that grade.
\end{introtheorem}

For two general linear blocks, Theorem~\ref{thm:complex} recovers Raghuram's combinatorial lemma over totally imaginary fields \cite[Lemma~3.16]{Raghuram2025}. Equation~\eqref{eq:graded-balance} gives the balance in each grade. For fixed inducing data, the simultaneously critical shifts are exactly the admissible points in an explicit open interval, which may be empty (Proposition~\ref{prop:interval}).

The proof brings together Saper's count of Kostant inversions in pairs of roots \cite[Lemma~17.6]{SaperMicro} and Shahidi's description of archimedean local coefficients in terms of Artin factors \cite[\S3]{Shahidi1985}.  At a real place, a fixed root contributes a real Gamma factor whose parity is determined by the Weil group action on its root line. At a complex place, the paired root coordinates have a common midpoint after division by grade. Thus degree balance forces opposite signs in every pair, precisely the condition for adjacent criticality of the corresponding Gamma factors.

\begin{remark}\label{rem:intro-overlap}
Theorem~\ref{thm:complex} proves the full combinatorial lemma at complex places. Under its local hypotheses, simultaneous criticality in every Langlands--Shahidi grade is equivalent to the occurrence of the two conjugate Levi coefficients in Kostant degrees whose sum is $\dim U$. 

After completing this work, the author learned that Raghuram and Zou had independently proved the forward implication at their evaluation point, obtaining complementary Kostant representatives from their hyper-criticality condition \cite[Def.~5.8 and Lemma~7.15]{RaghuramZou2026}. Remark~\ref{rem:complex-attribution} compares the normalizations.
\end{remark}

In the real general linear setting of Harder and Raghuram, balance means that the Kostant representative has length $\dim U/2$ \cite[Lemma~7.14]{HarderRaghuram2020}. At a complex place, the corresponding condition is that the two conjugate Kostant degrees sum to $\dim U$. The following corollary shows when one equality of total degrees detects these local conditions at every archimedean place.

\begin{corollary}[Global balance and simultaneous criticality]
Suppose the relative centers agree at all archimedean places. An admissible shift that is regular at every embedding is simultaneously critical at every archimedean place if and only if the total Kostant degrees on the two associate sides agree.

More precisely, in each grade the absolute difference between the Kostant degrees summed over all embeddings on the two sides equals the total pole order of the archimedean factors at the adjacent arguments and of their duals at the reflected arguments.
\end{corollary}

The common relative center forces all nonzero local degree differences to have the same sign. Their sum therefore vanishes exactly when each local difference vanishes \textup{(Theorem~\ref{thm:global-converse})}. The grade-by-grade pole count is Corollary~\ref{cor:graded-defect-poles}.

\begin{introtheorem}[Half-spin Eisenstein classes]
\label{thm:intro-E6}
Let $G$ be split of type $E_6$ over $\mathbb Q$, either simply connected or adjoint. For every integer $t\ge1$, take the ambient highest weight $2t\rho$. At a suitable common finite level, holomorphic Eisenstein series induced from $P_1$ at the centered shift of the half-spin family give classes in degrees $20$ and $21$ with nonzero restriction to the inducing face. In degree $20$, a single class arising from one cuspidal inducing datum restricts nontrivially to both boundary face types $P_1$ and $P_6$.
\end{introtheorem}

The sixteen roots of $U$ form eight pairs, with one inversion in each pair. The resulting Kostant degree eight takes Levi classes of degrees $12$ and $13$ to boundary degrees $20$ and $21$. Theorem~\ref{thm:E6-ordinary} shows that the image of restriction to the inducing face contains, over $\overline{\mathbb Q}$, the entire inner Levi Hecke subspace in the specified length-eight summand. Corollary~\ref{cor:E6-two-faces} gives the two-face assertion.

\subsection{Arithmetic questions beyond equal Kostant degrees}\label{subsec:arithmetic-unequal-degrees}
Harder's program also asks what Eisenstein cohomology can detect when the two Kostant degrees differ \cite[\S9.3.1]{HarderChapterNine}. The degree difference then gives two pieces of information. On the cohomological side, the two boundary contributions meet in a common degree precisely when a difference of Levi cohomology degrees compensates for the difference of Kostant degrees (Lemma~\ref{lem:degree-overlap}). On the analytic side, at a regular shift with a common relative center, the absolute degree difference in each grade counts the poles of the direct and reflected archimedean factors at the adjacent arguments (Corollary~\ref{cor:graded-defect-poles}).

For the $(3,1)$-parabolic of $\GL_4/\mathbb Q$, the symmetric-square lift of Ramanujan's cusp form gives a common boundary degree $4=3+1=2+2$. Of the adjacent arguments $-2$ and $-1$, the first is critical. At the second, a simple archimedean pole and nonvanishing at the reflected argument imply that the finite $L$-function has a simple zero (Example~\ref{ex:gl4-unbalanced}). This singles out both a derivative and a boundary degree in which to seek its cohomological interpretation.

The Klingen parabolic leads to a different arithmetic question. For $\Sp_{2n}$ over a totally real field, with $n\ge2$, suitable globally generic cuspidal Levi data realize the singular shift as the unique simultaneously critical shift at every real place (Proposition~\ref{prop:klingen-wall}). For $\GSp_4/\mathbb Q$, the denominator question arises in the neighboring chamber of Kostant length two. For the algebraically twisted inducing representation at the cohomological parameter, only the upper argument is critical (Example~\ref{ex:klingen-mismatch}). The author's Klingen denominator conjecture relates this critical value, after period normalization, to the denominator of an Eisenstein class \cite{HosseiniJafariKlingen}, following Harder's conjecture for the Siegel parabolic of $\Sp_4$ \cite[\S1.5.1, (93)]{HarderDenominators}.

Section~\ref{sec:root-setup} develops the coefficient line and the inverse Kostant criterion. Sections~\ref{sec:real-places} and~\ref{sec:complex-criterion} prove the local lemmas. Section~\ref{sec:cohomological-degrees} treats the global defect and boundary degrees. Section~\ref{sec:realclass-E6} constructs the half-spin classes. Appendix~\ref{app:fixed-roots} verifies the fixed roots in the real classification.

\subsection*{Acknowledgments} The author thanks Leslie Saper for introducing him to his work on $\mathcal L$-modules and micro-support, Sug Woo Shin for helpful conversations on purity, and Jayce Getz for many helpful conversations and for suggestions that improved the manuscript. The author gratefully acknowledges partial support from the National Science Foundation through the RTG grant DMS-2231514. Any opinions, findings, and conclusions or recommendations expressed in this material are those of the author and do not necessarily reflect the views of the National Science Foundation.

\section{Graded root coordinates for the inverse Kostant problem}\label{sec:root-setup}
The same root coordinates compute the Kostant degree and the archimedean criticality test. Their slopes are the Langlands--Shahidi grades, and opposition centers them in pairs (Proposition~\ref{prop:root-centering}).

\subsection{The inverse Kostant problem and its degree}
Let $F$ be a number field and $G$ a connected split reductive group over $F$. Fix a split maximal torus $T$, a Borel subgroup $B\supset T$, and a maximal proper standard parabolic $P=MU_P\supset B$. Put $X=X^*(T)$ and $W=W_G$, let $W_M$ be the Levi Weyl group, and write $w_G$ and $w_M$ for the longest elements. For $H\in\{G,M\}$, let $\Phi_H^+$ and $\Delta_H$ be the positive and simple roots determined by $B\cap H$, with half-sum $\rho_H$; put $\rho=\rho_G$ and $\Delta=\Delta_G$. Via restriction to $T$, identify the character group $X^*(H)$ with the integral weights annihilating every coroot of $H$, and set
\begin{align*}
 \mathfrak a_{H,\mathbb Q}^*=X^*(H)\otimes_{\mathbb Z}\mathbb Q,
 \qquad
 \mathfrak a_H^*=\mathfrak a_{H,\mathbb Q}^*\otimes_{\mathbb Q}\mathbb R,
 \\
 X_H^+=\{\nu\in X:\langle\nu,\alpha^\vee\rangle\ge0
                         \text{ for all }\alpha\in\Delta_H\}.   
\end{align*}
 
Fix $\overline{\mathbb Q}\subset\mathbb C$ and write $\Sigma_F=\operatorname{Hom}(F,\mathbb C)$. At a chosen embedding, $\mathfrak u_P$ is the complexified Lie algebra of $U_P$. Let $\cM_\lambda$ be the irreducible algebraic $G$-module of highest weight $\lambda\in X_G^+$, and let $\tcM_\lambda$ denote its associated local system. For $\mu\in X_M^+$, we write $\cM_\mu^M$ and $\tcM_\mu^M$ for the corresponding Levi module and its associated local system. Put
\[
 \Psi_P=\Phi_G^+\setminus\Phi_M^+,
 \qquad \rho_P=\rho-\rho_M,
 \qquad d_P=|\Psi_P|=\dim U_P.
 \]
 
The set $W^P=\{w\in W:w^{-1}\Phi_M^+\subset\Phi_G^+\}$ consists of the minimal representatives for $W_M\backslash W$ \cite[\S6.2.5]{HarderBook}. Their length is $\ell(w)=|\Inv_P(w)|$, where $\Inv_P(w)=\{\beta\in\Psi_P:w^{-1}\beta\in-\Phi_G^+\}$ \cite{Humphreys1990}. With $w\cdot\lambda=w(\lambda+\rho)-\rho$, Kostant's theorem \cite{Kostant} gives
\begin{equation}\label{eq:Kostant}
 H^q(\mathfrak u_P,\cM_\lambda)
   =\bigoplus_{\substack{w\in W^P\\\ell(w)=q}}
       \cM^M_{w\cdot\lambda}
 \qquad(\lambda\in X_G^+),
\end{equation}
which underlies the boundary-cohomology constructions in \cite[\S6.2.5]{HarderBook} and \cite{Schwermer1983}. The inverse problem asks whether a given $\mu\in X_M^+$ has the form $w\cdot\lambda$ with $\lambda\in X_G^+$ and $w\in W^P$.

A weight $x\in X\otimes_{\mathbb Z}\mathbb R$ is \emph{$H$-regular} if $\langle x,\alpha^\vee\rangle\ne0$ for every $\alpha\in\Phi_H$, and \emph{strictly $H$-dominant} if its pairings with the simple coroots of $H$ are positive.

In their study of boundary cohomology, Harder and Schwermer use Kostant's theorem to identify the Levi coefficients and their degrees. The following inverse criterion uses the chamber argument in Clozel's account of their method \cite[\S3.3, (3.5)--(3.6)]{Clozel2022}. It makes integrality and uniqueness explicit and expresses the cohomological degree as a count of negative root pairings.

\begin{proposition}[Inverse Kostant criterion]\label{prop:inverse-kostant}
Let $\mu\in X_M^+$ and put $x=\mu+\rho$. There exist $\lambda\in X_G^+$ and $w\in W^P$ with $\mu=w\cdot\lambda$ if and only if $x$ is $G$-regular. In that case $w$ is the unique element of $W$ for which $w^{-1}x$ is strictly $G$-dominant, and
\[
 \lambda=w^{-1}x-\rho,
 \qquad
 \Inv_P(w)=\{\beta\in\Psi_P:\langle x,\beta^\vee\rangle<0\}.
\]
In particular,
\begin{equation}\label{eq:length-sign}
 \ell(w)=\#\{\beta\in\Psi_P:\langle x,\beta^\vee\rangle<0\}.
\end{equation}
\end{proposition}

\begin{proof}
If $\mu=w\cdot\lambda$ for dominant $\lambda$, then $x=w(\lambda+\rho)$ is regular. Conversely, a regular $x$ has a unique strictly dominant translate $y=w^{-1}x$. For $\alpha\in\Phi_M^+$, strict Levi dominance gives $\langle y,(w^{-1}\alpha)^\vee\rangle=\langle x,\alpha^\vee\rangle>0$. Thus $w^{-1}\alpha$ is positive and $w\in W^P$. The weight $\lambda=y-\rho$ is integral because $w^{-1}\rho-\rho$ belongs to the root lattice. Its simple-coroot pairings are nonnegative because those of $y$ are positive integers. Finally $\langle x,\beta^\vee\rangle=\langle y,(w^{-1}\beta)^\vee\rangle$ is negative exactly on $\Inv_P(w)$.
\end{proof}

The argument works for any standard parabolic and allows $\lambda=0$. Maximality makes $\mathfrak a_M^*/\mathfrak a_G^*$ a line. Regularity holds on a half-line of twists in the setting of \cite[\S3.3, Theorem~3.6]{Clozel2022}. Section~\ref{sec:chamber-intervals} determines the chamber intervals on the full line.

\subsection{Grading, coefficient shifts, and induction parameters} \label{sec:grading-shifts}
\subsubsection{The relative weight and the dual grading}
Let $\alpha_P$ be the unique element of $\Delta\setminus\Delta_M$. The fundamental weight $\gamma_P\in\operatorname{span}_{\mathbb Q}(\Phi_G)$ is characterized by $\langle\gamma_P,\alpha_P^\vee\rangle=1$ and $\langle\gamma_P,\alpha^\vee\rangle=0$ for $\alpha\in\Delta_M$. Its class spans $\mathfrak a_{M,\mathbb Q}^*/\mathfrak a_{G,\mathbb Q}^*$, and the root-span condition fixes the representative
\cite[\S9.1.3]{HarderChapterNine}. Since $W_M$ permutes $\Psi_P$, the weight $\rho_P$ is $W_M$-invariant and 
\begin{equation}\label{eq:relative-weight}
 \rho_P=f_P\gamma_P,
 \qquad
 f_P:=\langle\rho_P,\alpha_P^\vee\rangle
     =1-\langle\rho_M,\alpha_P^\vee\rangle
     \in\tfrac12\mathbb Z,
 \qquad f_P\ge1,
\end{equation}
the inequality because $\rho_M$ is a nonnegative combination of Levi simple roots. Thus $\gamma_P$ is the weight $\widetilde\alpha$ of the Langlands--Shahidi normalization \cite{Shahidi2010}.

\begin{definition}[Root grades]\label{def:root-grades}
The \emph{grade} $j(\beta)=\langle\gamma_P,\beta^\vee\rangle$ of $\beta\in\Psi_P$ is the coefficient of $\alpha_P^\vee$ in the simple-coroot expansion of $\beta^\vee$. It is a positive integer because $\beta$ and $\beta^\vee$ have the same support.

\end{definition}

In nonsimply laced types the grade need not equal the coefficient of $\alpha_P$ in $\beta$. For the complex dual group $\widehat G$ and its standard Levi $\widehat M$, the positive dual nilradical decomposes as 
\[
 \mathfrak u_P^\vee=\bigoplus_{j\ge1}V_j,
 \qquad
 V_j=\bigoplus_{\beta\in\Psi_{P,j}}\widehat{\mathfrak g}_{\beta^\vee},
 \qquad
 \Psi_{P,j}:=\{\beta\in\Psi_P:j(\beta)=j\}.
\]
The adjoint action of $\widehat M$ preserves each $V_j$; its restriction $r_j$ is the graded representation of the Langlands--Shahidi construction \cite{Shahidi2010}. Write $\widetilde r_j=r_j^\vee$ for its contragredient. Set 
\[
 d_j:=\dim V_j=|\Psi_{P,j}|,
 \qquad
 \ell_j(w):=|\Inv_P(w)\cap\Psi_{P,j}|,
 \qquad
 \ell(w)=\sum_j\ell_j(w).
\]
Products over $j$ range over the nonzero $V_j$.

\subsubsection{Coefficient shifts and algebraic admissibility}
For $\mathcal P=\operatorname{Res}_{F/\mathbb Q}P$ inside $\mathcal G=\operatorname{Res}_{F/\mathbb Q}G$, an integral Levi coefficient weight is a family $\mu=(\mu^\eta)_{\eta\in\Sigma_F}$ with $\mu^\eta\in X_M^+$. Let $\gamma_{P,\eta}$ be supported on the $\eta$-factor.

\begin{definition}[Coefficient shifts]\label{def:coefficient-shift}
Put
\[
 \gamma_{\mathcal P}:=\sum_{\eta\in\Sigma_F}\gamma_{P,\eta},
 \qquad
 \mu+b\gamma_{\mathcal P}
        =(\mu^\eta+b\gamma_P)_{\eta\in\Sigma_F},
 \qquad b\in\mathbb R.
\]
The real number $b$ is the \emph{coefficient-shift parameter}.
\end{definition}

This family uses the Galois-invariant direction $\gamma_{\mathcal P}$ and the same parameter $b$ at every embedding. Root calculations are performed in one factor with the embedding suppressed; their conditions must hold at every embedding. Allowing the ambient central character to vary, put 
\begin{equation}\label{eq:lattice}
 \Lambda_{\rm alg}:=
 \{b\in\mathbb R:b\gamma_P+z\in X^*(M)
           \text{ for some }z\in\mathfrak a_{G,\mathbb Q}^*\}.
\end{equation}

\begin{lemma}\label{lem:algebraic-shift-lattice}
There is an integer $d_{\rm alg}>0$ such that $\Lambda_{\rm alg}=d_{\rm alg}\mathbb Z$. For $\mu\in X_M^+$, the weight $\mu+b\gamma_P+z$ is integral for some $z\in\mathfrak a_{G,\mathbb Q}^*$ if and only if $b\in\Lambda_{\rm alg}$. Every such shift preserves the Levi coroot pairings of $\mu$.
\end{lemma}

\begin{proof}
Let $\mu\in X_M^+$. The set $\Lambda_{\rm alg}$ is a subgroup of $\mathbb R$. It lies in $\mathbb Z$, since $b=\langle b\gamma_P+z,\alpha_P^\vee\rangle$, and it is nonzero, since a positive multiple of $\gamma_P$ lies in $X^*(M)$. If $\mu+b\gamma_P+z$ is integral, then $b\gamma_P+z$ is integral and annihilates the Levi coroots, hence is a character of $M$; the converse is immediate.
\end{proof}

\begin{definition}[Algebraically admissible shifts] \label{def:admissible-shift}
A shift $b$ is \emph{algebraically admissible} if $b\in\Lambda_{\rm alg}$. Choose $z_b\in\mathfrak a_{G,\mathbb Q}^*$ so that $\theta_b=b\gamma_P+z_b\in X^*(M)$. Use the same character at every embedding and put $\mu_b^\eta=\mu^\eta+\theta_b\in X_M^+$, or $\mu_b$ when the embedding is suppressed. Then $\cM_{\mu_b^\eta}^M\simeq\cM_{\mu^\eta}^M\otimes\mathbb C_{\theta_b}$.
\end{definition}

For $G=\mathrm{GL}_N$ and $M=\mathrm{GL}_n\times\mathrm{GL}_{n'}$ with $N=n+n'$,
\[
 \gamma_P=\left(\frac{n'}N\mathbf1_n,-\frac nN\mathbf1_{n'}\right),
 \qquad
 b\gamma_P+\frac{bn}{N}\mathbf1_N=(b\mathbf1_n,0),
\]
so $\Lambda_{\rm alg}=\mathbb Z$ and the corrected shift is $\det^b$ on the first factor.

For $G=\mathrm{PGL}_2$ and $P=B$, the simple root $\alpha$ generates $X^*(T)$, $\gamma_P=\alpha/2$, and $\mathfrak a_{G,\mathbb Q}^*=0$. Thus $\Lambda_{\rm alg}=2\mathbb Z$. For a prescribed ambient central character, restrict the corrections in \eqref{eq:lattice} to those compatible with it.

\subsubsection{Induction parameters and root hyperplanes}
At a place $v$ of $F$, choose $n>0$ with $n\gamma_P\in X^*(M)$. Define the positive character $|\gamma_P|_v$ and, for $s\in\mathbb C$,
\[
 |\gamma_P|_v(m):=|(n\gamma_P)(m)|_v^{1/n},
 \qquad \chi_s(m):=|\gamma_P|_v(m)^s.
\]
The definition is independent of $n$; here $|\cdot|_v$ is the normalized absolute value ($|z|_{\mathbb C}=z\bar z$). Since the adjoint action on $\operatorname{Lie}U_P$ has determinant of weight $2\rho_P$, 
\[
 \delta_P=|\gamma_P|_v^{2f_P},
 \qquad \delta_P^{1/2}=\chi_{f_P}.
\]
For an admissible representation $\sigma$ of $M(F_v)$, realized smoothly with moderate growth at archimedean places, $I_P(s,\sigma)$ is normalized induction of $\sigma\otimes\chi_s$ with sections satisfying 
\[
f(mug)=\delta_P(m)^{1/2}\chi_s(m)\sigma(m)f(g).
\]
Writing ${}^{\rm a}\operatorname{Ind}$ for unnormalized induction, 
\begin{equation}\label{eq:induction-conversion}
 I_P(s,\sigma)
   ={}^{\rm a}\operatorname{Ind}_P^G
          (\sigma\otimes\chi_{s+f_P}),
 \qquad k_P:=-f_P.
\end{equation}
Thus $I_P(k_P,\sigma)={}^{\rm a}\operatorname{Ind}_P^G(\sigma)$; $k_P$ is the \emph{evaluation point}
\cite[Definition~2.1.1]{RaghuramPadic2023}.

\begin{proposition} \label{prop:coefficient-spectral-parameters}
At an archimedean place, suppose that $\sigma$ is irreducible and cohomological with coefficient $\cM_\mu^M$, in the convention $\sigma\otimes\cM_\mu^M$. Work in one embedding factor and put
\[
 \nu_M=-w_M(\mu+\rho_M),\qquad
 x(b)=\mu+\rho+b\gamma_P.
\]
The infinitesimal character of $I_P(s,\sigma)$ is represented, modulo $W_G$, by $\nu_P(s):=\nu_M+s\gamma_P$. For every real $b$, 
\begin{equation}\label{eq:coefficient-spectral}
 \nu_P(k_P-b)=-w_Mx(b),
 \qquad s=k_P-b,
 \qquad s+f_P=-b.
\end{equation}
At a complex place these formulas hold in both embedding factors.
\end{proposition}

\begin{proof}
By Wigner's lemma \cite[Ch.~I, Theorem~5.3(ii)]{BorelWallach}, $\sigma$ has the infinitesimal character of the dual coefficient, of highest weight $-w_M\mu$ and parameter $\nu_M$. Twisting by $\chi_s$ adds $s\gamma_P$, and normalized induction preserves the parameter, now modulo $W_G$ \cite[Lemma~3.2(2)]{Vogan2000}. Since $w_M\gamma_P=\gamma_P$ and $\rho=\rho_M+f_P\gamma_P$, we get $\nu_M-(f_P+b)\gamma_P=-w_M(\mu+\rho+b\gamma_P)$. 
\end{proof}

For $b=0$ and $\mu=w\cdot\lambda$, projecting \eqref{eq:coefficient-spectral} to $\mathfrak a_M^*$ along the span of the Levi roots gives the Levi-central component of $-w(\lambda+\rho)$. At $s=k_P$ its class modulo $\mathfrak a_G^*$ gives the cohomological evaluation point in the boundary construction \cite{Schwermer1983}
(see also \cite[(4.1)]{Clozel2022}).

\begin{remark}[Central corrections and real components] \label{rem:central-correction-parameter}
For the corrected coefficient of Definition~\ref{def:admissible-shift}, $-w_M(x(b)+z_b)=\nu_P(k_P-b)-z_b$: an exact comparison also changes the inducing parameter by $-z_b$, while all root pairings are unchanged. At a real place,
\[
 \theta_b(m)^{-1}
  =\operatorname{sgn}(\theta_b(m))
      |\gamma_P|_v(m)^{-b}|z_b|_v(m)^{-1},
\]
with $|z_b|_v$ defined by the same denominator-clearing convention. The component character $\varepsilon_b=\operatorname{sgn}\circ\theta_b$ is not detected by the infinitesimal parameter. For every $\beta\in\Psi_P$, 
\[
 \varepsilon_b(\beta^\vee(-1))
   =(-1)^{\langle\theta_b,\beta^\vee\rangle}
   =(-1)^{b j(\beta)}.
\]
Equivalently, the Langlands parameter of $\varepsilon_b$ takes values in the center of $\widehat M$ and acts on $V_j$ through the character $\operatorname{sgn}^{bj}$ of $W_{\mathbb R}$. This records the sign contributed by the algebraic twist; the full real Weil action also includes the action already present in the inducing parameter, and the real-place criteria below retain the actual Weil signs.
\end{remark}

For $\beta\in\Psi_P$, define
\begin{equation}\label{eq:q-H}
 q_\beta(b):=\langle x(b),\beta^\vee\rangle
       =\langle\mu+\rho,\beta^\vee\rangle+bj(\beta)
       =-\langle\nu_P(k_P-b),(w_M\beta)^\vee\rangle.
\end{equation}
For $\alpha\in\Phi_G$, the hyperplane $H_\alpha=\{x:\langle x,\alpha^\vee\rangle=0\}$ in
$X\otimes_{\mathbb Z}\mathbb R$ is a \emph{wall}. A \emph{regular shift} is a parameter $b$ for which $x(b)$ lies on no wall. At a regular admissible shift, the ambient highest weight is strictly dominant if and only if the dominant Weyl translate of $x(b)$ pairs to at least two with every simple coroot of $G$. At a regular admissible shift, the ambient highest weight is strictly dominant if and only if the dominant Weyl translate of $x(b)$ pairs to at least two with every simple coroot of $G$.

\begin{corollary}[Walls along the coefficient line]\label{cor:shift-walls}
For every real $b$, the weight $x(b)$ is strictly $M$-dominant. It is $G$-regular exactly when $b$ avoids the wall-crossing parameters
\begin{equation}\label{eq:wall-parameters}
 b_\beta=-\frac{\langle\mu+\rho,\beta^\vee\rangle}{j(\beta)}
 \qquad(\beta\in\Psi_P).
\end{equation}
On each component of $\mathbb R\setminus\{b_\beta:\beta\in\Psi_P\}$, the signs of all root pairings and the chamber containing $x(b)$ are constant.
\end{corollary}

\begin{proof}
For $\delta\in\Delta_M$, 
\[\langle x(b),\delta^\vee\rangle =\langle\mu,\delta^\vee\rangle+1>0.\] 
Thus $x(b)$ pairs positively with every positive Levi coroot. For each $\beta\in\Psi_P$, the pairing \eqref{eq:q-H} has positive slope $j(\beta)$ and unique zero $b_\beta$.
\end{proof}

\subsection{Opposition-purity, root centering and the sign count}\label{sec:purity}
We express the condition of conjugate contragredience modulo the Levi center and derive a common center for the graded root pairs, which supplies the sign condition for the criticality criteria. The term \emph{opposition-purity} refers to the opposition involution $-w_M$, which sends a Levi highest weight to the highest weight of the dual representation. Recall that
\[
 w_M\rho_M=-\rho_M,\qquad w_M\gamma_P=\gamma_P,
 \qquad \rho+w_M\rho=2f_P\gamma_P,
\]
and that $w_M$ permutes $\Psi_P$ with $j(w_M\beta)=j(\beta)$.

\begin{definition}[Opposition-purity]\label{def:opposition-purity}
Let $\mu=(\mu^\eta)_{\eta:F\hookrightarrow\mathbb C}$ be a family of $M$-dominant integral weights, and put $\bar\eta=\mathfrak c\circ\eta$, where $\mathfrak c$ is complex conjugation. We call $\mu$ \emph{locally opposition-pure} if, for every archimedean place $v$ represented by $\eta$, 
\begin{equation}\label{eq:local-purity}
 \mu^{\bar\eta}+w_M\mu^\eta=c_v\in\mathfrak a_{M,\mathbb Q}^*,
 \qquad a_v:=\langle c_v,\alpha_P^\vee\rangle.
\end{equation}
We call $\mu$ \emph{strongly opposition-pure} if every reindexed family
\[
 ({}^{\tau}\mu)^\eta:=\mu^{\tau^{-1}\circ\eta},
 \qquad \tau\in\operatorname{Aut}(\mathbb C),
\]
is locally opposition-pure. The central weight may depend on both $\tau$ and $v$.
\end{definition}

Local opposition-purity is the placewise derived-group form of Saper's conjugate self-contragredience condition \cite[Def.~7.2(ii), Prop.~8.4]{SaperMicro}. Indeed, \eqref{eq:local-purity} says precisely that the coefficients at $\eta$ and $\bar\eta$ become dual to one another on restriction to $M_{\mathrm{der}}$. Their highest weights therefore differ, after duality, by the Levi character $c_v$, which may vary with $v$. Saper imposes his condition on ${}^0L_P$. For a restriction of scalars, removing the rational split center can leave central directions linking different archimedean places \cite[\S0.11.5, (0.2.1)]{SaperMicro}.

Since $w_M$ fixes Levi characters, the relation is symmetric in the two embeddings. Maximality of $P$ gives
$c_v-a_v\gamma_P\in\mathfrak a_{G,\mathbb Q}^*$. Adding Levi characters $\zeta^\eta$ preserves local and strong opposition-purity and adds $\zeta^\eta+\zeta^{\bar\eta}$ to $c_v$. The common shift $\theta_b$ of Definition~\ref{def:admissible-shift} therefore adds $2b$ to $a_v$ and commutes with reindexing.

\begin{remark}[Local purity and the global center]
\label{rem:purity-centers}
For $M=\mathrm{GL}_n$, \eqref{eq:local-purity} reads $\mu_i^\eta+\mu_{n+1-i}^{\bar\eta}=\mathrm{w}_v\in\mathbb Z$ for $1\le i\le n$. For the coefficient of a cohomological cuspidal representation of $\GL_n(\mathbb A_F)$, Clozel's purity lemma gives a single purity weight $\mathrm w$ at all archimedean places \cite[Lemme~4.9]{Clozel1990}. Strong purity requires every arithmetic conjugate of the coefficient to be pure with the same weight $\mathrm w$ \cite[\S2.3.1]{Raghuram2025}. In the totally real setting of Harder and Raghuram \cite{HarderRaghuram2020}, this follows from purity because arithmetic conjugation permutes the real embeddings.

For example, the weights $(1,0)$ and $(2,0)$ at the two real embeddings of a real quadratic field are strongly opposition-pure, with distinct purity weights $1$ and $2$. Example~\ref{ex:nonCM} illustrates the other distinction: a family with one purity weight can lose the symmetry after reindexing. For $M=\mathrm{GL}_n\times\mathrm{GL}_{n'}$ the two block weights give $a_v=\mathrm{w}_v-\mathrm{w}'_v$.
\end{remark}

\begin{lemma} \label{lem:essential-unitarity-purity}
Let $\sigma_v$ be irreducible, essentially unitary, and cohomological with local coefficient weights $(\mu^\xi)_{\xi\mid v}$ in the convention of Proposition~\ref{prop:coefficient-spectral-parameters}. Then these weights satisfy \eqref{eq:local-purity}.
\end{lemma}

\begin{proof}
This is the standard coefficient-purity consequence of unitarity \cite[Ch.~II, Proposition~6.12(1) and \S6.13]{BorelWallach}; in our notation it reads as follows. Modulo Levi-central weights we may assume $\sigma_v$ unitary, so that its complex conjugate is its contragredient. Conjugation exchanges the two infinitesimal-character factors at a complex place and fixes the single factor at a real place, while contragredience negates them. The derived coordinates of $\nu_M^\xi=-w_M(\mu^\xi+\rho_M)$ are real, so uniqueness of the dominant representative gives $\nu_M^{\bar\eta}\equiv-w_M\nu_M^\eta$ modulo $\mathfrak a_M^*$. With $w_M\rho_M=-\rho_M$, this says that the integral weight $\mu^{\bar\eta}+w_M\mu^\eta$ is Levi-central.
\end{proof}

Local opposition-purity therefore holds automatically for the essentially tempered cohomological data in the local criteria. Proposition~\ref{prop:root-centering} identifies the common midpoint within each place in terms of $a_v$. The global balance criterion requires $a_v$ to be the same at every place. For automorphic data, Lemma~\ref{lem:global-center} gives the stronger conclusion that $c_v$ is common under its connected central-character hypothesis.

\begin{proposition}[Center of the graded root pairs]\label{prop:root-centering}
Assume that $\mu$ is locally opposition-pure, and restore the embedding index in \eqref{eq:q-H}: $q_\beta^\xi(b):=\langle\mu^\xi+\rho+b\gamma_P,\beta^\vee\rangle$. For every $b\in\mathbb R$ and $\beta\in\Psi_P$, 
\begin{equation}\label{eq:root-centering}
 q_\beta^\eta(b)+q_{w_M\beta}^{\bar\eta}(b)
   =2j(\beta)(b-b_v^\circ),\qquad b_v^\circ=-f_P-\frac{a_v}{2}.
\end{equation}
After division by $j(\beta)$, every pair has midpoint $b-b_v^\circ$. Write $T_v(b)=2(b-b_v^\circ)=a_v-2s$ for the normalized sum. 
\end{proposition}

\begin{proof}
By symmetry of \eqref{eq:local-purity} and Weyl invariance,
\[
 \begin{aligned}
 q_\beta^\eta(b)+q_{w_M\beta}^{\bar\eta}(b)
 &=\langle\mu^\eta+w_M\mu^{\bar\eta}
        +\rho+w_M\rho+2b\gamma_P,\beta^\vee\rangle\\
 &=\langle c_v+2(f_P+b)\gamma_P,\beta^\vee\rangle
 =j(\beta)(a_v+2f_P+2b),
 \end{aligned}
\]
the last step because $c_v-a_v\gamma_P$ is ambient central. Finally, $s=-f_P-b$.
\end{proof}

The identity holds for real shifts, and $z_b$ leaves the pairings unchanged. Since each $j(\beta)>0$, the paired sums are all nonnegative or all nonpositive. The following count applies when the coordinates are nonzero. 

\begin{lemma}\label{lem:common-sign}
Let $(x_i,y_i)_{i\in I}$ be a finite family of pairs of nonzero real numbers, and put $D=|I|$. Suppose that either $x_i+y_i\geq0$ for every $i$, or $x_i+y_i\leq0$ for every $i$. Then the total number of negative entries in the family is $D$ if and only if $x_iy_i<0$ for every $i$.
\end{lemma}
\begin{proof}
Let $n_-$ and $n_+$ count the pairs with two negative and two positive entries, respectively. The number of negative entries is $D+n_--n_+$. The common-sign hypothesis forces $n_-=0$ or $n_+=0$, so this number equals $D$ exactly when both counts vanish.
\end{proof}

\begin{remark} \label{rem:saper-dictionary}
In the notation of \cite[Prop.~8.4, Lemma~17.6 and Basic Lemma~24.2]{SaperMicro}, the involution $\tau'_P$ at a complex place sends $(\eta,\beta)$ to $(\bar\eta,w_M\beta)$, and $\Phi_w$ is the union of the two inversion sets. For a conjugate self-contragredient constituent occurring in nilradical cohomology, Saper counts Kostant inversions in pairs of roots to bound the number of inversions in each relative root space. For the effect of restriction of scalars on the Levi center, see the discussion following Definition~\ref{def:opposition-purity}.

Here the same pairing groups the dual root lines by Langlands--Shahidi grade.
\end{remark}

\subsection{Chambers and algebraically admissible shifts} \label{sec:chamber-intervals}
The signs of the root coordinates determine the interval of shifts realizing each Kostant representative and its length in every grade.

Fix an embedding and suppress it. Let $C=\{x:\langle x,\alpha^\vee\rangle>0 \text{ for every }\alpha\in\Delta\}$ be the strictly $G$-dominant chamber in $X\otimes_{\mathbb Z}\mathbb R$. Whenever $x(b)$ is $G$-regular, let $w(b)$ be the unique element with $w(b)^{-1}x(b)\in C$. The argument of Proposition~\ref{prop:inverse-kostant} gives, without algebraic admissibility,
\[
 w(b)\in W^P,\qquad
 \Inv_P(w(b))=\{\beta\in\Psi_P:q_\beta(b)<0\},
\]
and Corollary~\ref{cor:shift-walls} gives
\begin{equation}\label{eq:root-line}
 q_\beta(b)=j(\beta)(b-b_\beta),\qquad j(\beta)>0.
\end{equation}

\begin{proposition}[Chambers along the coefficient line]\label{prop:intervals}
For $w\in W^P$, the set of parameters with $x(b)$ $G$-regular and $w(b)=w$ is \begin{equation}\label{eq:Iw}
 I_w=(b_w^-,b_w^+),\qquad
 b_w^-=\max_{\beta\in\Psi_P\setminus\Inv_P(w)}b_\beta,
 \qquad b_w^+=\min_{\beta\in\Inv_P(w)}b_\beta,
\end{equation}
with $\max\varnothing=-\infty$ and $\min\varnothing=+\infty$. The interval is empty if $b_w^-\ge b_w^+$, and both endpoints are finite if $0<\ell(w)<d_P$. The nonempty intervals $I_w$ are precisely the components of $\mathbb R\setminus\{b_\beta:\beta\in\Psi_P\}$.
Whenever $x(b)$ is $G$-regular, 
\begin{equation}\label{eq:line-length}
 \ell_j(w(b))=\#\{\beta\in\Psi_{P,j}:b<b_\beta\},
 \qquad
 \ell(w(b))=\#\{\beta\in\Psi_P:b<b_\beta\}.
\end{equation}
There are $1+\#\{b_\beta:\beta\in\Psi_P\}\le d_P+1$ nonempty chamber intervals, with equality if and only if all crossing parameters are distinct.
\end{proposition}

\begin{proof}
By \eqref{eq:root-line}, $w(b)=w$ means $b<b_\beta$ for $\beta\in\Inv_P(w)$ and $b>b_\beta$ for the other $\beta\in\Psi_P$; with the positive Levi-root pairings, these signs specify $wC$. This gives \eqref{eq:Iw}, and both index sets are nonempty when $0<\ell(w)<d_P$. Signs are constant between consecutive crossing parameters and at least one changes at each, which gives the components and their number. Counting negative pairings in each grade proves \eqref{eq:line-length}.
\end{proof}

\begin{corollary}[Algebraic realization in a prescribed chamber] \label{cor:algebraic-chamber-realizations}
Let $b\in\Lambda_{\rm alg}$, with $z_b$ and $\mu_b=\mu+b\gamma_P+z_b$ as in Definition~\ref{def:admissible-shift}. For $w\in W^P$, there exists $\lambda_b\in X_G^+$ with $\mu_b=w\cdot\lambda_b$ if and only if $b\in I_w$. In that case,
\[
 \lambda_b=w^{-1}x(b)-\rho+z_b.
\]
Thus the algebraically admissible shifts with representative $w$ form $I_w\cap\Lambda_{\rm alg}$, which may be empty even when $I_w$ is nonempty.
\end{corollary}

\begin{proof}
The correction $z_b$ is fixed by $W$ and annihilates every coroot of $G$, so $\mu_b+\rho=x(b)+z_b$ has the root pairings of $x(b)$. Apply Proposition~\ref{prop:inverse-kostant}.
\end{proof}

For restriction of scalars, Definitions~\ref{def:coefficient-shift} and~\ref{def:admissible-shift} use the same parameter and central correction at every embedding. The shifts realizing a prescribed tuple $(w_\eta)_\eta$ are therefore $\Lambda_{\rm alg}\cap\bigcap_\eta I_{w_\eta}^{\eta}$, with $I_{w_\eta}^{\eta}$ computed using $\mu^\eta$.

\begin{remark}[Comparison with the relative-weight coordinate] \label{rem:harder-coordinate}
In the relative-weight coordinates of \cite[\S9.1.3]{HarderChapterNine}, uniquely write $\mu\equiv\mu^{(1,M)}+a_0\gamma_P$ modulo $\mathfrak a_{G,\mathbb Q}^*$, with $\mu^{(1,M)}\in\operatorname{span}_{\mathbb Q}(\Phi_M)$ and $a_0\in\mathbb Q$, and put $\widetilde\mu^{(1)}:=\mu^{(1,M)}+\rho_M$. Then 
\begin{equation}\label{eq:harder-coordinate}
 x(b)\equiv\widetilde\mu^{(1)}+b_{\rm H}\gamma_P
 \pmod{\mathfrak a_{G,\mathbb Q}^*},
 \qquad b_{\rm H}=a_0+f_P+b.
\end{equation}
Thus $b$ measures displacement from the fixed coefficient $\mu$, whereas $b_{\rm H}$ is the relative coordinate of $x(b)$; the endpoints of $I_w$ become $a_0+f_P+b_w^-$ and $a_0+f_P+b_w^+$. The relative coordinate of $\mu_b$ is $a_0+b$.
\end{remark}

\subsection{Arithmetic conjugation and passage to the associate parabolic}
\subsubsection{Testing and transporting purity}\label{sec:purity-transport}
Strong opposition-purity means that
\begin{equation}\label{eq:strong-purity}
 \mu^{\tau^{-1}\circ\bar\eta}+w_M\mu^{\tau^{-1}\circ\eta}
 \in\mathfrak a_{M,\mathbb Q}^*
 \qquad(\tau\in\operatorname{Aut}(\mathbb C),\ \eta\in\Sigma_F),
\end{equation}
where in general $\tau^{-1}\circ\bar\eta$ need not equal $\overline{\tau^{-1}\circ\eta}$.

\begin{proposition}\label{prop:finite-test}
Fix an embedding $F\subset\mathbb C$ and a finite Galois extension $L/\mathbb Q$ in $\mathbb C$ containing the normal closure of $F$. Put $\Gamma=\operatorname{Gal}(L/\mathbb Q)$, $c_L=\mathfrak c|_L$, and $\mathcal C=\{g^{-1}c_Lg:g\in\Gamma\}$. Then $\mu$ is strongly opposition-pure if and only if 
\begin{equation}\label{eq:finite-test}
 \mu^{\delta\circ\xi}+w_M\mu^\xi
 \in\mathfrak a_{M,\mathbb Q}^*
 \qquad(\xi\in\operatorname{Hom}(F,L),\ \delta\in\mathcal C).
\end{equation}
Thus it is tested by finitely many equalities in $(X\otimes_{\mathbb Z}\mathbb Q)/\mathfrak a_{M,\mathbb Q}^*$. 
\end{proposition}

\begin{proof}
Every embedding of $F$ into $\mathbb C$ has image in $L$. For $\tau\in\operatorname{Aut}(\mathbb C)$, put $g=\tau|_L$ and $\xi=g^{-1}\circ\eta$. Then $\tau^{-1}\circ\eta=\xi$ and $\tau^{-1}\circ\bar\eta=(g^{-1}c_Lg)\circ\xi$, so \eqref{eq:finite-test} implies \eqref{eq:strong-purity}. Conversely, every $g\in\Gamma$ extends to an automorphism $\tau$ of $\mathbb C$, and $\eta=g\circ\xi$ in \eqref{eq:strong-purity} gives \eqref{eq:finite-test}. 
\end{proof}

One can adapt the descent argument of \cite[Prop.~2.6]{Raghuram2025}
to weights modulo the Levi center. Put $\overline X_M:=(X\otimes_{\mathbb Z}\mathbb Q)/\mathfrak a_{M,\mathbb Q}^*$ and write $[\nu]$ for the class of a weight $\nu$. Since $w_M$ fixes Levi characters, $-w_M$ induces an involution on $\overline X_M$.

\begin{proposition} \label{prop:derived-purity-descent}
Let $F_1$ be the maximal CM subfield of $F$, if one exists, and otherwise the maximal totally real subfield of $F$. Then $\mu$ is strongly opposition-pure if and only if there is a family $(\kappa^\zeta)_{\zeta\in\Sigma_{F_1}}$ in $\overline X_M$ such that 
\[
 [\mu^\eta]=\kappa^{\eta|_{F_1}},\qquad
 \kappa^{\bar\zeta}=-w_M\kappa^\zeta.
\]
The descended family of classes is unique.
\end{proposition}

\begin{proof}
With the notation of Proposition~\ref{prop:finite-test}, let $N$ be the subgroup of $\Gamma$ generated by the products $\delta\delta'$ with $\delta,\delta'\in\mathcal C$, and put $E=L^N$. Then $N$ is normal, and all elements of $\mathcal C$ have the same image $t$ in $\Gamma/N$, which is central with $t^2=1$. Every subfield $K$ of $E$ is therefore stable under $t$, which is complex conjugation in every embedding of $K$; so $K$ is totally real if $t|_K=1$ and CM otherwise. Conversely, every element of $\mathcal C$ acts trivially on a totally real subfield of $L$ and as the CM involution on a CM subfield, so $N$ fixes both and they lie in $E$. Consequently $F_1=F\cap E$, and every embedding of $F_1$ has image in $E$.

Put $H=\operatorname{Gal}(L/F)$. Then $\operatorname{Gal}(L/F_1)=NH$, so the fibres of the restriction map $\operatorname{Hom}(F,L)=\Gamma/H\to\operatorname{Hom}(F_1,L)=\Gamma/NH$ are precisely the $N$-orbits. The finite test reads $[\mu^{\delta\xi}]=-w_M[\mu^\xi]$. Since $(-w_M)^2=1$, it implies $N$-invariance, hence descent to $F_1$, and its instance $\delta=c_L$ gives the relation for $\kappa$. Conversely, all $\delta\in\mathcal C$ agree with $c_L$ on $E$, so that relation for a descended family implies the finite test. Surjectivity of restriction on embeddings gives uniqueness.
\end{proof}

For totally real and CM fields $F_1=F$, so local and strong opposition-purity coincide; compare \cite[\S2.3.3]{Raghuram2025}. 

\begin{remark}[Comparison with usual strong purity]\label{rem:strong-purity-descent}
Over a totally imaginary field the full strongly pure $\mathrm{GL}_n$ weight descends to $F_1$ \cite[Prop.~2.6]{Raghuram2025}. For the corresponding result over general number fields, see
\cite[Prop.~3.3]{DarshanRaghuram2024}.

Proposition~\ref{prop:derived-purity-descent} gives descent of the weight modulo Levi characters for any split Levi. The relative coordinate along $\gamma_P$ can still vary over an embedding of $F_1$, so the critical intervals require a separate comparison. 
\end{remark}

\begin{example}[Local opposition-purity need not be strong] \label{ex:nonCM}
Adapting \cite[\S2.3.4]{Raghuram2025}, let $F=\mathbb Q(\sqrt[3]{2},\omega)$ with $\omega=e^{2\pi i/3}$, and identify its embeddings with $S_3$, complex conjugation $c=(23)$ acting on the left. For $M=\mathrm{GL}_2$, set $\kappa=(1,-1)$, and let $\mu^g=0$ for $g\in\{e,c\}$ and $\mu^g=\kappa$ otherwise. These weights are dominant integral, $\mu^{cg}=\mu^g$, and $w_M\kappa=-\kappa$, so every local opposition sum is zero. For an automorphism $\tau$ restricting to $d=(12)$ and $\eta=d$, however, 
\[
 ({}^{\tau}\mu)^\eta=0,\qquad
 ({}^{\tau}\mu)^{\bar\eta}=\mu^{dcd}=\kappa
 \notin\mathbb Q(1,1).
\]
Thus the family is pure of weight zero in the usual $\mathrm{GL}_2$ sense, but not strongly opposition-pure: here $N=A_3$ and $F_1=\mathbb Q(\omega)$, and the derived weights are not constant on the fibres of restriction to $F_1$. Since usual strong purity is necessary for cuspidal $\mathrm{GL}_2$-cohomology \cite[\S1.2]{BhagwatRaghuram2017}, this family has no cuspidal realization.
\end{example}

\begin{corollary}[Purity from conjugate realizations] \label{cor:realized-strong-purity}
Suppose that, for every $\tau\in\operatorname{Aut}(\mathbb C)$ and every archimedean place $v$, the coefficient determined by the components of ${}^{\tau}\mu$ at $v$ admits an irreducible essentially unitary cohomological representation of $M(F_v)$. Then $\mu$ is strongly opposition-pure. 
\end{corollary}

\begin{proof}
Apply Lemma~\ref{lem:essential-unitarity-purity} to every ${}^{\tau}\! \mu$.
\end{proof}

Let $\sigma$ be a cuspidal automorphic representation of $M(\mathbb A_F)$,
cohomological with coefficient family $\mu$. For each
$\tau\in\operatorname{Aut}(\mathbb C)$ there is a discrete automorphic
representation with finite part ${}^{\tau}\sigma_f$ (scalar conjugation
of the smooth representation) and nonzero cohomology for ${}^{\tau}\mu$.
This follows by applying \cite[Prop.~2.4]{GanRaghuram2013} to
$\mathbf M:=\operatorname{Res}_{F/\mathbb Q}M$.

A discrete automorphic representation, in the convention allowing split-central twists, is essentially unitary. Indeed, let $A_{\mathbf M}$ be the maximal $\mathbb Q$-split central torus of $\mathbf M$ and $\mathbf M(\mathbb A)^1$ the common kernel of the adelic absolute values of rational characters. The decomposition 
\[
 \mathbf M(\mathbb A)
 =A_{\mathbf M}(\mathbb R)^0\times\mathbf M(\mathbb A)^1
\]
allows a positive automorphic character, trivial on $\mathbf M(\mathbb A)^1$, to cancel the modulus of the split-central character \cite[\S3, p.~16]{Arthur2005}. The archimedean components of the conjugate representation are therefore essentially unitary, and Corollary~\ref{cor:realized-strong-purity} applies to it without twisting.

For $\mathrm{GL}_n$ one can take the conjugate representation to be cuspidal \cite{Clozel1990} (see also \cite[Theorem~3.1]{GanRaghuram2013}). For a general Levi, cuspidality and essential temperedness require
further input.

\subsubsection{Opposition on the associate parabolic}\label{sec:associate-parabolic}
The opposition automorphism $\Theta=-w_G$ preserves $\Delta$. Let $Q$ be the standard parabolic whose Levi $M_Q$ has simple roots $\Delta_{M_Q}=\Theta(\Delta_M)$; it is conjugate to the opposite of $P$ and may equal $P$. Define $W^Q$, $\Psi_Q$, and $\Inv_Q$ as for $P$. Following \cite[\S9.1.3]{HarderChapterNine}, put
\[
 w_P:=w_Mw_G,\qquad v_P:=w_P^{-1}=w_Gw_M,
 \qquad \iota:=-v_P=\Theta w_M.
\]

\begin{proposition}\label{prop:associate-complement}
The element $w_P$ is the longest element of $W^P$. The map
\[
 i_{P,Q}:W^P\longrightarrow W^Q,\qquad w\longmapsto w':=v_Pw
\]
is a bijection, and $\iota$ maps $\Psi_P$ bijectively onto $\Psi_Q$. Moreover,
\begin{equation}\label{eq:complement-length}
 \Inv_Q(w')=\iota\bigl(\Psi_P\setminus\Inv_P(w)\bigr),
 \qquad \ell(w')=d_P-\ell(w).
\end{equation}
Let $\gamma_Q$ and $j_Q$ be the fundamental weight and grade for $Q$, and $\ell_j^Q(u)=\#\{\beta'\in\Inv_Q(u):j_Q(\beta')=j\}$. Then $\gamma_Q=\Theta\gamma_P$ and
$j_Q(\iota\beta)=j(\beta)$. Consequently,
\begin{equation}\label{eq:graded-complement}
 \ell_j^Q(w')=d_j-\ell_j(w).
\end{equation}
\end{proposition}

\begin{proof}
Since $w_M\Phi_M^+=-\Phi_M^+$ and $w_M\Psi_P=\Psi_P$, we have $v_P\Phi_M^+=\Phi_{M_Q}^+$ and $v_P\Psi_P=-\Psi_Q$. Hence $w_P\in W^P$ has inversion set $\Psi_P$. For $w\in W^P$,
$(w')^{-1}\Phi_{M_Q}^+=w^{-1}\Phi_M^+\subseteq\Phi_G^+$, so $w'\in W^Q$; the inverse map is $u\mapsto w_Pu$. For $\beta\in\Psi_P$ one has $(w')^{-1}(\iota\beta)=-w^{-1}\beta$, which proves \eqref{eq:complement-length}. Finally $\Theta$ sends $\gamma_P$ to $\gamma_Q$, and $w_M\gamma_P=\gamma_P$ gives $\langle\gamma_Q,(\iota\beta)^\vee\rangle
 =\langle\Theta\gamma_P,(\Theta w_M\beta)^\vee\rangle =\langle\gamma_P, \beta^\vee\rangle$; intersecting with each grade proves \eqref{eq:graded-complement}.
\end{proof}

\begin{lemma} \label{lem:associate-coefficient-line}
For $\mu\in X_M^+$ put
\[
 \widetilde\mu:=v_P\mu,\qquad
 \mu':=v_P(\mu+\rho)-\rho
      =\widetilde\mu-2f_P\gamma_Q.
\]
Both weights belong to $X_{M_Q}^+$. If $\mu=w\cdot\lambda$, then $\mu'=w'\cdot\lambda$ for $w'=v_Pw$. Writing $x_P^\mu(b)=\mu+\rho+b\gamma_P$ and using the same notation for $Q$, one has 
\begin{equation}\label{eq:associate-lines}
 x_Q^{\mu'}(-b)
 =x_Q^{\widetilde\mu}(-b-2f_P)
 =v_Px_P^\mu(b).
\end{equation}
Let $I_u^R(\kappa)$ denote the chamber interval of Proposition~\ref{prop:intervals} for parabolic $R$ and base weight $\kappa$. Then
\begin{equation}\label{eq:associate-intervals}
 I_{w'}^Q(\mu')=-I_w^P(\mu),\qquad
 I_{w'}^Q(\widetilde\mu)=-2f_P-I_w^P(\mu).
\end{equation}
In the spectral convention of Proposition~\ref{prop:coefficient-spectral-parameters},
\[
 s=k_P-b,\qquad b''=-b-2f_P
 \quad\Longrightarrow\quad k_Q-b''=-s.
\]
\end{lemma}

\begin{proof}
Put $\rho_Q=\rho-\rho_{M_Q}$. The root transport gives $\rho_Q=\Theta\rho_P=f_P\gamma_Q$, so $f_Q=f_P$, $v_P\gamma_P=-\gamma_Q$, and $v_P\rho=\rho-2\rho_Q$. The weight $v_P\mu$ is integral and $M_Q$-dominant, and $2\rho_Q$ is a character of $M_Q$; this proves the assertions about the weights and \eqref{eq:associate-lines}. Since $(w')^{-1}v_P=w^{-1}$, the two sides of \eqref{eq:associate-lines} have the same strictly dominant representative, which gives
\eqref{eq:associate-intervals}. Finally $k_Q=k_P=-f_P$.
\end{proof}

An ambient central correction is fixed by $v_P$: $v_P(x_P^\mu(b)+z_b)-\rho=\mu'-b\gamma_Q+z_b$, and $v_P(b\gamma_P+z_b)=-b\gamma_Q+z_b$ identifies the algebraic shift lattices for $P$ and $Q$. Since $2f_P\gamma_Q=2\rho_Q\in X^*(M_Q)$, both transformations of $b$ in \eqref{eq:associate-lines} preserve algebraic admissibility.

\begin{remark}\label{rem:associate-versus-local-opposition}
The pairing $\beta\mapsto w_M\beta$ of Section~\ref{sec:purity} stays inside $\Psi_P$, whereas $\iota=\Theta w_M$ transports roots to $\Psi_Q$. These associate-parabolic identities hold without opposition-purity.
\end{remark}

\subsection{Local factors and intertwining normalization}
\label{sec:factor-conventions}
Let $Q=v_P\bar Pv_P^{-1}$ and $v_P=w_Gw_M$ as in Section~\ref{sec:associate-parabolic}. Fix an $F$-rational pinning compatible with $(T,B)$. Use its Tits representative $\dot v_P$ globally and at every place, and put $(v_P\sigma)(m')=\sigma(\dot v_P^{-1}m'\dot v_P)$ on $M_Q(F_v)$.

At infinity take the maximal compact $K_{G,v}$ attached to the pinned Cartan involution and at unramified places the pinned hyperspecial subgroup. The element $\dot v_P$ belongs to each: its simple-reflection factors are images of $\bigl(\begin{smallmatrix}0&1\\-1&0\end{smallmatrix}\bigr)$ under the pinned $\mathrm{SL}_2$ homomorphisms \cite[proof of Proposition~2.8]{HowardMillson2004}. Thus conjugation carries $K_{M,v}=K_{G,v}\cap M(F_v)$ to $K_{M_Q,v}=K_{G,v}\cap M_Q(F_v)$, including when $Q=P$. It carries $K_{M,\infty}=\prod_{v\mid\infty}K_{M,v}$, its component characters, and $A_M=A_{\mathbf M}(\mathbb R)^0$ to their counterparts for $M_Q$. This identifies the relative complexes for either compact convention, with or without the split factor. In particular, $v_P\sigma$ is cohomological for $\cM^{M_Q}_{\widetilde\mu}$ with $\widetilde\mu=v_P\mu$, since $v_P\Phi_M^+=\Phi_{M_Q}^+$.

\begin{remark}[Which compact subgroup]\label{rem:compact-conventions}
Use $K_{M,\infty}A_M$, or $K_{M,\infty}^{\circ}A_M$ on the connected cover. The full-compact complex is the $\pi_0(K_{M,\infty})$-invariant subcomplex of the latter, with the component group acting on cochains. A central $z\in Z(M(F_\infty))\cap K_{M,\infty}$ acts trivially on their exterior-power domain, so the full complex vanishes unless $\omega_{\sigma_\infty}(z)\omega_\mu(z)=1$ \cite[Ch.~I, \S5.1]{BorelWallach}. For $\sgn^\epsilon|\cdot|^t\otimes\cM_u$ on $\GL_1(\R)$ this gives $\epsilon+u\in2\Z$ on the full compact. Section~\ref{sec:realclass-E6} uses the full compact, with $\zeta_M(K_{M,\infty})=\{\pm1\}$.

\end{remark}

Since $v_P\gamma_P=-\gamma_Q$, the standard intertwining operator is 
\begin{equation}\label{eq:standard-intertwiner}
 \begin{aligned}
 A_v(s,\sigma)&:I_P(s,\sigma)\longrightarrow I_Q(-s,v_P\sigma),\\
 (A_v(s,\sigma)f)(g)
   &:=\int_{U_Q(F_v)}f(\dot v_P^{-1}ug)\,du,
 \end{aligned}
\end{equation}
defined by absolute convergence and meromorphic continuation; for these conventions at finite places see \cite[\S2.2]{RaghuramPadic2023}. A reverse operator uses $\dot v_P^{-1}$, which may differ from the Tits representative of $w_P$ by an element of order two in $T(F)$. By Lemma~\ref{lem:associate-coefficient-line}, a source shift $b$ for $\mu$ corresponds to the target shift $b''=-b-2f_P$ for $\widetilde\mu$.

At an unramified place $v$, assume that $\sigma$ is irreducible and has a $K_{M,v}$-fixed vector. Give $U_Q(F_v)\cap K_{G,v}$ volume one. Let $f_v^\circ,f_v^{\prime\circ}$ be the spherical sections with the same fixed vector at the identity. If $q_v$ is the residue-field size and $t_{\sigma,v}$ the Satake parameter, use the convention $\beta^\vee(t_{\xi,v})=\xi(\beta^\vee(\varpi_v))$ for unramified torus characters and set
\[
 L_v(t,\sigma,r_j)=\det(1-q_v^{-t}r_j(t_{\sigma,v}))^{-1}.
\]
The positive dual nilradical gives the representations $r_j$. Grouping the Gindikin--Karpelevich factors by their grades \cite[Ch.~4]{Langlands1971} gives
\begin{equation}\label{eq:adjacent-tuple}
 \mathcal R_{P,v}(s,\sigma)
   :=\prod_j\frac{L_v(js,\sigma,r_j)}{L_v(js+1,\sigma,r_j)},
\end{equation}
and
\begin{equation}\label{eq:gindikin-karpelevich}
 A_v(s,\sigma)f_v^\circ
   =\mathcal R_{P,v}(s,\sigma)f_v^{\prime\circ}.
\end{equation}
Thus $N_v(s,\sigma):=\mathcal R_{P,v}(s,\sigma)^{-1}A_v(s,\sigma)$ sends $f_v^\circ$ to $f_v^{\prime\circ}$ as a meromorphic identity. For the upper triangular Borel of $\mathrm{GL}_2$ and $\sigma=\xi_1\boxtimes\xi_2$, this convention gives $\mathcal R_{B,v}(s,\sigma)
=L_v(s,\xi_1\xi_2^{-1})/L_v(s+1,\xi_1\xi_2^{-1})$.

At an archimedean place, use the local Langlands parameter $\phi_{\sigma_v}$. Since $M$ is split, extend $r_j$ to ${}^LM=\widehat M\times W_{F_v}$ trivially on the Weil factor, and define
\[
 L_v(t,\sigma,r_j):=L(t,r_j\circ\phi_{\sigma_v})
\]
by Artin local factors \cite[\S3]{Shahidi1985}, and $\mathcal R_{P,v}$ and $N_v$ by the same formulas. Real-place factors thus retain the action of the full Weil group, including its component outside $\mathbb C^\times$.

\begin{remark}[Comparison of normalizations]\label{rem:factor-conventions}
At infinity the action on a cohomological constituent requires a local
calculation. For $\GL_{2m}(\R)$ with equal blocks, $m$ even, and
tempered inducing data cohomological for trivial Levi coefficients,
the full intertwining operator has a pole at the evaluation point.
Its restriction to the unique irreducible submodule is holomorphic there
\cite[\S4.2]{Clozel2022}.

One must match the Weil representations, the spectral parameter and the direction of the operator \cite[(1.0.1), \S2.3]{RaghuramPadic2023} \cite{Shahidi1990,Shahidi2010}. Whittaker normalization also requires generic data and compatible functionals and measures.
\end{remark}

\subsubsection{Simultaneous criticality and the complex Gamma factor}
\begin{definition}[Simultaneous criticality]\label{def:simultaneous}
Let $v$ be an archimedean place, and let $s\in k_P-\mathbb Z$, the spectral lattice in our cohomological normalization. We call $s$ \emph{simultaneously critical at $v$} if, for every nonempty grade $j$, all four factors 
\begin{equation}\label{eq:four-factors}
 L_v(js,\sigma_v,r_j),\quad
 L_v(js+1,\sigma_v,r_j),\quad
 L_v(1-js,\sigma_v,\widetilde r_j),\quad
 L_v(-js,\sigma_v,\widetilde r_j)
\end{equation}
are finite.
 
We call $s$ \emph{simultaneously critical} if this condition holds at every archimedean place. For an algebraically admissible shift $b\in\Lambda_{\rm alg}$, we call $b$ a \emph{critical twist at $v$}, respectively a \emph{critical twist}, if $s=k_P-b$ is simultaneously critical at $v$, respectively at every archimedean place. 
\end{definition}

For the complex Gamma factor write $\Gamma_\C(t)=2(2\pi)^{-t}\Gamma(t)$. If a character $z^A\bar z^B$ satisfies $A+js=-q$ and $B+js=-q'$, its local factor is $L(t)=\Gamma_\C(t+\max(A,B))$ \cite[p.~13]{Tate1979}. Substitution gives the four arguments in \eqref{eq:four-factors}, in order,
\begin{equation}\label{eq:complex-arguments}
 -\min(q,q'),\quad 1-\min(q,q'),\quad 1+\max(q,q'),\quad\max(q,q').
\end{equation}
For integral $q,q'$, all four factors are finite if and only if $\min(q,q')<0<\max(q,q')$, or equivalently $qq'<0$. At a regular shift, this is equivalent to exactly one of the two coordinates being negative, so the pair contributes one inversion to the count in \eqref{eq:length-sign}. This calculation also applies to the $\Gamma_\C$ factor attached to an exchanged pair of real Weil weights \cite[(4.6)]{Knapp1994}. 

\begin{remark}
At a complex place, Proposition~\ref{prop:complex-dictionary} shows that, for $b\in\Lambda_{\rm alg}$, every Gamma argument in \eqref{eq:four-factors} is an integer. Hence the four factors are finite precisely when all these arguments are positive. At a real place the one-dimensional Weil characters introduce the parity test of Lemma~\ref{lem:parity}. 
\end{remark}

\begin{remark}\label{rem:simultaneous-names}
At the evaluation point for each maximal parabolic of split $G_2$, simultaneous criticality specializes to the adjacent-criticality conditions of \cite[\S5, Lemmas~6.1--6.2]{HosseiniJafariG2}. The Langlands--Shahidi
grading expresses these conditions uniformly as criticality of $js$ and $js+1$ in every grade. For compatibility with the critical-point lattices, see \cite[\S7.2]{HosseiniJafariG2}.
\end{remark}

\section{Real places: fixed lines and the regular classification} \label{sec:real-places}
For a prescribed cohomological representation and admissible shift, the real combinatorial lemma tests adjacent criticality on each Weil orbit. Its first consequence is a classification of the parabolics for which the regular case can occur. The Klingen example realizes the other possibility, a singular critical wall, for global cuspidal inducing data.

\subsection{The combinatorial lemma at real places}\label{subsec:real-criterion}
Let $\mathcal F_P=\{\beta\in\Psi_P:w_M\beta=\beta\}$ and $n_{\rm fix}(P)=|\mathcal F_P|$. Write $\mathcal O_P$ for the set of two-element $w_M$-orbits in $\Psi_P$.
Write $W_\R=\C^\times\sqcup J\C^\times$, $JzJ^{-1}=\bar z$, $J^2=-1$. Fix an irreducible, essentially tempered representation $\sigma_v$ of $M(\R)$, cohomological with coefficients in $\cM_\mu^M$ for an integral $M$-dominant weight $\mu$. Lemma~\ref{lem:essential-unitarity-purity} gives local opposition-purity:
\[
 \mu+w_M\mu=c_v\in\mathfrak a_{M,\mathbb Q}^*,\qquad
 a_v=\langle c_v,\alpha_P^\vee\rangle.
\]
Write $\phi_{\sigma_v}$ for its local Langlands parameter. Since $M$ is split, we identify its dual-group part with a homomorphism into $\widehat M$. Keep $b\in\Lambda_{\rm alg}$, $s=k_P-b$, and the coordinates $x(b)$ and $q_\beta(b)$ of \eqref{eq:q-H}. The integral coefficient is $\mu_b=\mu+\theta_b$, where $\theta_b=b\gamma_P+z_b$ is as in Definition~\ref{def:admissible-shift}; the correction $z_b$ does not change these root coordinates. The default analytic test uses the fixed representation $\sigma_v$.

At a real embedding, Proposition~\ref{prop:root-centering} reads
\begin{equation}\label{eq:real-centering}
 q_\beta(b)+q_{w_M\beta}(b)=2j(\beta)(b-b_v^\circ).
\end{equation}
For a fixed root the two coordinates coincide, so
\begin{equation}\label{eq:fixed-wall}
 q_\beta(b)=j(\beta)(b-b_v^\circ).
\end{equation}
Thus every fixed root vanishes at the same shift.

Definition~\ref{def:simultaneous} concerns the fixed representation $\sigma_v$, that is, the absolute-value family of Remark~\ref{rem:twist-families}; this remains our primary convention. For the algebraic family, we call $b\in\Lambda_{\rm alg}$ a critical twist at $v$ if the four factors of $\sigma_v\otimes\theta_b^{-1}$ at $jk_P$ and $jk_P+1$ are finite in every grade.

For a fixed root $\beta$, write $(-1)^{\epsilon_\beta(0)}$ for the action of $J$ on its dual root line. A \emph{twist family} is one of these two, with fixed-line signs $\epsilon_\beta(b)=\epsilon_\beta^{\rm abs}(b)$ or $\epsilon_\beta(b)=\epsilon_\beta^{\rm alg}(b)$. The root coordinates are the same for both.

\begin{theorem}[Real-place criticality and the fixed-wall alternative] \label{thm:real-criterion}
Let $\sigma_v$ be an irreducible, essentially tempered representation of $M(\R)$, cohomological for $\cM_\mu^M$ with $\mu\in X_M^+$. Let $b\in\Lambda_{\rm alg}$ and $s=k_P-b$, and fix a twist family. Then the shift $b$ is a critical twist at $v$ for the chosen family if and only if 
\begin{enumerate}
\item $q_\beta(b)q_{w_M\beta}(b)<0$ for each $\{\beta,w_M\beta\}\in\mathcal O_P$;
\item $q_\beta(b)=0$ and $\epsilon_\beta(b)=1$ for every $\beta\in\mathcal F_P$.
\end{enumerate}
The critical twists for the chosen family form the \emph{critical set at $v$}.
\begin{enumerate}
\renewcommand{\labelenumi}{(\alph{enumi})}
\item If $\mathcal F_P=\varnothing$, criticality is equivalent to a dominant integral ambient realization
\[
 \mu_b=u\cdot\lambda,\qquad u\in W^P,\qquad
 \ell(u)=d_P/2.
\]
Whenever these equivalent conditions hold, $x(b)$ is regular and $\ell_j(u)=d_j/2$ for every nonempty grade $j$, where $d_j=|\Psi_{P,j}|$.
\item If $\mathcal F_P\ne\varnothing$, the critical set is either empty or the singleton $\{b_v^\circ\}$. It is the singleton precisely when 
\begin{equation}\label{eq:singleton}
 b_v^\circ\in\Lambda_{\mathrm{alg}},\qquad
 \epsilon_\beta(b_v^\circ)=1\ (\beta\in\mathcal F_P),\qquad
 q_\beta(b_v^\circ)\ne0\ (\beta\in\Psi_P\setminus\mathcal F_P).
\end{equation}
At this point $x(b_v^\circ)$ is singular, and the corrected coefficient $\mu_{b_v^\circ}$ has no dominant ambient Kostant realization.
\end{enumerate}
\end{theorem}

The root conditions do not depend on the twist family, so in both families the critical set is either an interval on the lattice or a possible singleton; only the sign test at $b_v^\circ$ differs. The proof is given after the endpoint tests in Lemma~\ref{lem:parity} of Section~\ref{subsec:real-weil}.

\subsection{Classification of the regular real branch} \label{sec:realclass-classification}
For which maximal parabolics can the regular case of Theorem~\ref{thm:real-criterion} occur for some integral Levi weight? The following definition states the question entirely in terms of the coefficient.
\begin{definition}[Regular root configuration] \label{def:regular-root-configuration}
A regular root configuration is a pair $(\mu,b)$ with $\mu\in X_M^+$, $b\in\Lambda_{\rm alg}$, and
$\mu+w_M\mu\in\mathfrak a_{M,\mathbb Q}^*$, such that 
\[
 q_\beta(b)q_{w_M\beta}(b)<0
 \qquad(\beta\in\Psi_P).
\]
\end{definition}
These inequalities imply $G$-regularity of $x(b)$ and exclude fixed normal roots. Existence in this subsection is quantified over coefficients and chambers. Application to a representation also requires the cohomological local hypotheses of Theorem~\ref{thm:real-criterion} and the chosen compact-component convention. 

\begin{theorem}[Classification by the omitted node] \label{thm:realclass-classification}
For a maximal parabolic $P_k$ of a connected split reductive $G/\R$, mark the omitted node $k$ in its component $\Phi$. A regular root configuration exists exactly at the zeros in the table (Bourbaki numbering). 
\[
\begin{array}{lc@{\qquad}lc}
\toprule
 \Phi & n_{\rm fix}(P_k) & \Phi & (n_{\rm fix}(P_k))_k\\
\midrule
 A_{N-1} & \mathbf1_{\{k,\ N-k\ {\rm odd}\}}
      & E_6 & (0,1,1,1,1,0)\\
 B_n,\ C_n & \lceil k/2\rceil
      & E_7 & (1,4,3,3,4,2,3)\\
 D_n,\ k\le n-2 & \lfloor k/2\rfloor+2\mathbf1_{\{k,\ n-k\ {\rm odd}\}}
      & E_8 & (2,4,4,4,4,3,3,1)\\
 D_n,\ k=n-1,n & \lfloor n/2\rfloor & F_4 & (1,2,2,1)\\
 &&G_2 &(1,1)\\
\bottomrule
\end{array}
\]
Thus the regular nodes are $k$ with $k(N-k)$ even in $A_{N-1}$, the line node $1$ in $D_{2m+1}$ for $m\ge2$, and nodes $1,6$ in $E_6$. We use $D_n$ for $n\ge4$, with $D_3=A_3$ and $D_2=A_1\times A_1$. \end{theorem}

A listed parabolic can have an empty critical interval for a prescribed coefficient. For the line parabolic of $\SO_{2n+2}$ with $n$ even, the two root coordinates involving the last orthogonal weight coordinate coincide at every shift when $\mu_n=0$. Hence the critical interval is empty for this coefficient (Section~\ref{subsec:classical-tests}).

The proof follows Theorem~\ref{thm:fixed-roots-fundamental-rank}, which relates the fixed-root count to the fundamental rank.

In Saper's notation, the involutions on the split based root datum are $\tau_G=-w_G$, $\tau_P=-w_M$, and $\tau'_P=-\tau_P=w_M$ \cite[(8.2.3)]{SaperMicro}. Thus for $u\in W^P$ the condition $u\tau_Gu^{-1}=\tau_P$ becomes $u w_Gu^{-1}=w_M$. It characterizes the fundamental Weyl elements. Equivalently, $w_M$ interchanges $\Inv_P(u)$ and its complement in $\Psi_P$. Such an element has length $d_P/2$ \cite[Lemmas~8.7(ii)--(iii) and 8.8(i)]{SaperMicro}.

The following proposition realizes the corresponding chambers by integral coefficients on every isogeny form.

\begin{proposition}[Opposition and existence] \label{prop:realclass-opposition}
A regular root configuration with chamber representative $u\in W^P$ satisfies
\[
 w_M=u w_Gu^{-1}.
\]
Regular root configurations exist if and only if $\mathcal F_P=\varnothing$. When they exist, one can take $\mu=u\rho-\rho$ in the root lattice, $b=0$, and the trivial ambient coefficient $\lambda=0$. Consequently this coefficient construction works for every isogeny form. 
\end{proposition}

\begin{proof}
Let \(x\in uC\) be strictly Levi dominant and regular, with opposite signs on every normal pair. Then \(-w_Mx\) is strictly Levi dominant because \(-w_M\) preserves the positive Levi roots. On a normal root \(\beta\), its pairing is \(-\langle x,(w_M\beta)^\vee\rangle\), which has the same sign as \(\langle x,\beta^\vee\rangle\).
Thus \(x\) and \(-w_Mx\) lie in the same chamber. Since \(-C=w_GC\), this gives \(uC=w_Mu w_GC\). The Weyl group acts freely on chambers, so \(w_M=u w_Gu^{-1}\). A fixed normal root would be paired with itself and is therefore excluded. 

Conversely suppose $\mathcal F_P=\varnothing$. In the semisimple root space the $-1$ eigenspace $V_-$ of $w_M$ contains $\rho_M$. No normal-root functional vanishes on all of $V_-$, since such a root would be fixed by $w_M$. Choose $x\in V_-$ near $\rho_M$, strictly Levi dominant and off the normal-root hyperplanes. Its opposed coordinates have opposite signs. Applying the chamber argument above to this \(x\) gives \(w_Mu\rho=uw_G\rho=-u\rho\). Thus $u\rho$ lies in $V_-$ and in the same chamber as $x$. Its Levi simple-coroot pairings are positive integers, so $\mu=u\rho-\rho$ is $M$-dominant and belongs to the root lattice. Since $\mu+w_M\mu=-\rho-w_M\rho=-2\rho_P$, the pair $(\mu,0)$ satisfies the purity and sign conditions with $\mu=u\cdot0$.
\end{proof}

The trivial ambient coefficient in the construction is just a starting point. For any such $u$, put $\lambda_t=2t\rho$ and $\mu_t=(2t+1)u\rho-\rho=u\cdot\lambda_t$ for integers $t\ge1$. Both weights are in the root lattice because $2\rho$ and $u\rho-\rho$ are. The normal signs persist and $\mu_t+w_M\mu_t=-2\rho_P$; every simple Levi pairing of $\mu_t$ is at least $2t$ since the corresponding pairing of $u\rho$ is a positive integer. Thus every classified case admits strictly dominant ambient and Levi highest weights, for every isogeny form. 

\begin{theorem}[Fixed normal roots and fundamental rank] \label{thm:fixed-roots-fundamental-rank}
Let $G$ be a connected split reductive group over $\mathbb R$, let $P=MU_P$ be a maximal proper parabolic, and choose a positive root system for which $P$ is standard. Put 
\[
 \mathcal F_P=\{\beta\in\Phi_G^+\setminus\Phi_M^+:w_M\beta=\beta\},
 \qquad n_{\rm fix}(P)=|\mathcal F_P|.
\]
Recall that $d_P=|\Psi_P|=\dim_{\mathbb R}U_P(\R)$. For a real reductive group $H$ write 
\[
 \ell_0(H)=\operatorname{rk}_{\mathbb C}\mathfrak h_{\mathbb C}
       -\operatorname{rk}_{\mathbb C}\mathfrak k_{H,\mathbb C},
 \qquad
 q_0(H)=\frac{\dim_{\mathbb R}(H(\mathbb R)/K_{H,v})-\ell_0(H)}2,
\]
where $K_{H,v}$ is a maximal compact subgroup, using the choice of Section~\ref{sec:root-setup} for $H=G,M$. All split central directions are retained in these invariants; see \cite[Ch.~III, \S\S4.1, 4.3]{BorelWallach}. Then $\mathcal F_P$ is pairwise strongly orthogonal and \begin{equation}\label{eq:fixed-root-rank-degree}
 n_{\rm fix}(P)=\ell_0(M)-\ell_0(G),\qquad
 q_0(G)=q_0(M)+\frac{d_P+n_{\rm fix}(P)}{2}.
\end{equation}
Choose a Cartan involution $\theta$ stabilizing $M$, with $K_{M,v}=K_{G,v}\cap M(\mathbb R)$. There is an identification of the roots of a $\theta$-stable fundamental Cartan $H_M\subset M(\mathbb R)$ with the chosen abstract root system under which its Cartan involution acts by $-w_M$. Its real roots in $G$ are precisely $\mathcal F_P\sqcup(-\mathcal F_P)$. In particular, the following are equivalent:
\begin{enumerate}
\item $\mathcal F_P=\varnothing$;
\item $H_M$ is fundamental in $G(\mathbb R)$;
\item $\ell_0(M)=\ell_0(G)$;
\item $M$ contains a $G(\R)$-conjugate of $Z_G(A_f)$, where $A_f$
is the connected split part of a fundamental Cartan of $G(\R)$;
\item $P$ contains a fundamental real parabolic.
\end{enumerate}
Condition~(5) is equivalent to the existence of a fundamental Weyl element \cite[Lemma~8.7(i)--(ii)]{SaperMicro}. For the fundamental-parabolic criterion, see also \cite[\S1.8]{BorelCasselman1983}.
\end{theorem}

\begin{proof}
Strong orthogonality is Lemma~\ref{lem:fixed-roots-strongly-orthogonal} in Appendix~\ref{app:fixed-roots}. To identify the Cartan involution, start from a split Cartan of $M$. Its Cartan involution acts by $-1$. Kostant's cascade is a strongly orthogonal set of positive Levi roots whose commuting reflections multiply to $w_M$ \cite[Prop.~1.10]{KostantCascade2012}.

Perform the real Cayley transforms through those roots. For a real root $\alpha$, transport by its Cayley transform changes $\theta$ to $s_\alpha\theta$: the transform fixes $\ker\alpha$ and changes the Cartan-involution sign on the coroot line \cite[p.~392, (6.67a)--(6.68c)]{Knapp2002}. Starting from $\theta=-1$, the commuting transforms therefore give the involution $-w_M$. No Levi root is fixed by $w_M$, since $w_M$ interchanges positive and negative Levi roots. Thus the resulting Cartan has no real Levi roots and is fundamental in $M$ by \cite[Proposition~6.70]{Knapp2002}. As all fundamental Cartans are conjugate, this gives the asserted identification for any choice of $H_M$, after conjugacy. 

For the identification of the Cartan involution with $\tau_P$ on a fundamental Levi Cartan, see \cite[\S8.5]{SaperMicro}. The real roots of this Cartan in $G$ are those $\beta$ with $(-w_M)\beta=-\beta$, namely the $w_M$-fixed roots. There are no fixed Levi roots, and $w_M$ preserves the positive normal roots, so these roots are precisely $\pm\mathcal F_P$. 

Saper uses a strongly orthogonal subset of the positive real roots to obtain a dimension bound under the equal-rank hypothesis \cite[proof of Lemma~25.2, (25.2.2)]{SaperMicro}. For the maximal parabolics considered here, the full set $\mathcal F_P$ of positive real roots is strongly orthogonal by Lemma~\ref{lem:fixed-roots-strongly-orthogonal}. 

Perform real Cayley transforms through every root in $\mathcal F_P$. For strongly orthogonal real roots, the corresponding rank-one root subalgebras commute. Hence each unused root remains real after the other transforms, and the transforms commute. If $s_{\mathcal F}=\prod_{\beta\in\mathcal F_P}s_\beta$, the transported Cartan involution is $\theta'=-w_M s_{\mathcal F}$. On the real root space $V$, its \emph{minus-one} eigenspace is  
\[
 V^{\theta'=-1}=V^{w_M}\cap\mathcal F_P^\perp.
\]
There is no root in this space: such a root would be $w_M$-fixed and hence one of $\pm\mathcal F_P$, but could not be orthogonal to itself. The transformed Cartan is therefore fundamental in $G$. Every transform reduces the split dimension by one. The split dimension of a fundamental Cartan of $H$ is $\ell_0(H)$; this follows by taking its compact part to be a maximal torus of $K_{H,v}$. Thus $n_{\rm fix}(P)=\ell_0(M)-\ell_0(G)$, including the unchanged split center of $G$.

Finally, $G(\mathbb R)=U_P(\mathbb R)M(\mathbb R)K_{G,v}$ gives $\dim(G/K_{G,v})=d_P+\dim(M/K_{M,v})$. Substitution into the definition of $q_0$ gives the second identity. This proves (1)--(3), since $H_M$ is fundamental in $G$ exactly when it has no real roots in $G$. To see that (5) is equivalent to (1), use the fundamental Weyl element criterion \cite[Lemma~8.7]{SaperMicro}. Such an element interchanges inversions and noninversions, which excludes fixed roots. Conversely, when there are no fixed roots, Proposition~\ref{prop:realclass-opposition} constructs such an element. For (4), when $H_M$ is fundamental in $G$, its split part $A_f$ contains the connected split center of $M$, hence $Z_G(A_f)\subseteq M$. Conversely, the Levi in (4) contains a fundamental $G$-Cartan. This Cartan is also fundamental in $M$, since its compact dimension is already maximal in $G$. Conjugacy of fundamental Cartans in $M$ gives (2). 
\end{proof}

\begin{proof}[Proof of Theorem~\ref{thm:realclass-classification}]
By Proposition~\ref{prop:realclass-opposition} absence of fixed roots is equivalent to existence, including integral coefficients on every isogeny form. Theorem~\ref{thm:fixed-roots-fundamental-rank} gives $n_{\rm fix}=\ell_0(M)-\ell_0(G)$. The ambient split center and every component away from $k$ cancel in this difference. On the component containing the omitted node, deleting $k$ gives one split central direction and the simple Levi components $M_i$, hence the count is $1+\sum_i\ell_0(M_i)-\ell_0(G_{\rm active})$. The Cayley-transform argument in Theorem~\ref{thm:fixed-roots-fundamental-rank} identifies the fundamental rank of a split semisimple group $H$ with the dimension of the $+1$ eigenspace of $w_H$ on its root space. The longest-element actions described in \cite[\S3]{ConradBruhat} therefore give $\lceil N/2\rceil-1$ for $A_{N-1}$, $n\bmod2$ for $D_n$, $2$ for $E_6$, and zero for the other types. Substitution gives the table; Appendix~\ref{app:fixed-roots} identifies the individual roots.
\end{proof}

For $\SL_N$ a fundamental Levi has blocks of size two and, for odd $N$, one of size one. In odd $D_n$ the line Levi has compact Cartan on its $D_{n-1}$ derived group. In $E_6$ the $D_4$ on $\{2,3,4,5\}$ has compact rank four, that of $E_6(\R)$, giving a fundamental Levi inside both end Levis. For the classification of the Weyl-group involutions appearing in Proposition~\ref{prop:realclass-opposition}, see \cite[Theorem~A]{Richardson1982}. 

When $c_H$ central dimensions are removed, the reduced defect is $\delta_H=\ell_0(H)-c_H$; removing them from both dimension and defect leaves $q_0(H)$ unchanged. In particular $n_{\rm fix}(P)=\delta_M-\delta_G+1$ when the full split centers are removed. The two neighboring lengths at a critical fixed wall, given by Proposition~\ref{prop:face} and the sign conditions of Theorem~\ref{thm:real-criterion}, send $q_0(M)$ to $q_0(G)$ and $q_0(G)-n_{\rm fix}(P)$.

For the Borel--Casselman cohomological criterion with a specified ambient module, see \cite[Prop.~8.9]{SaperMicro}. Saper's inversion bounds \cite[Lemma~17.6]{SaperMicro} follow a paired-root argument of Saper--Stern \cite{SaperStern1990}. For related bounds, see \cite[Lemma~15.2 and Props.~18.2--18.3]{Saper2003} and, in rank one, \cite[Prop.~2.6]{Casselman1983}. 

\emph{Classification by group.}
For $G$ of positive semisimple rank, at least one maximal parabolic admits a regular root configuration precisely when the diagram has a component of type $A_r$ with $r\ge2$, odd $D_n$ with $n\ge5$, or $E_6$. Every maximal parabolic admits one precisely when every component has type $A_{2m}$, $m\ge1$. A split torus has no proper parabolic. Maximality is essential. For the Borel of split $\SL_3$ all three positive roots are fixed since $M$ is a torus, whereas $\ell_0(M)-\ell_0(G)=1$.

\begin{corollary}[Critical twists and ambient Kostant realizations] \label{cor:real-adjacent-pairs}
Let $v$ be a real place, let $\sigma_v$ and $\mu$ be as in Theorem~\ref{thm:real-criterion}, and fix a twist family. If $n_{\rm fix}(P)>0$, any critical twist at $v$ must equal $b_v^\circ$. If such a twist exists, its corrected integral coefficient $\mu_{b_v^\circ}$ has no dominant ambient Kostant realization. Consequently, simultaneous criticality can coexist with an ordinary ambient Kostant realization only when the active simple component and omitted node belong to the list in Theorem~\ref{thm:realclass-classification}. 
\end{corollary}

\begin{proof}
The first two assertions are Theorem~\ref{thm:real-criterion}(b). If $n_{\rm fix}(P)=0$, part~(a) of that theorem gives an ambient realization for every critical twist, and Theorem~\ref{thm:realclass-classification} lists these parabolics. 
\end{proof}

Existence of a regular configuration concerns coefficients. The cuspidal inducing data in the exceptional regular case are constructed in Lemma~\ref{lem:E6-centered-existence}.

\begin{remark}[Residues and jets at the fixed wall]\label{rem:wall-jets} 
Keep the inducing datum fixed and set $s_0=k_P-b_v^\circ$ at a critical fixed wall. The induced infinitesimal character is represented by $\nu_P(s_0)=-w_M(\mu+\rho+b_v^\circ\gamma_P)$. Its Weyl orbit is singular. The representation generated by a finite Laurent jet of this fixed inducing family has the same generalized infinitesimal character. Indeed, clear the pole and differentiate $DE(s)=\chi_{\nu_P(s)}(D)E(s)$ for $D\in Z(U(\mathfrak g_\C))$. On the resulting finite jet, $D-\chi_{\nu_P(s_0)}(D)$ lowers jet order and is nilpotent. This applies also to residues and subquotients.

Every finite-dimensional algebraic $\cM_\lambda^\vee$ has regular shifted infinitesimal character. Choose a central element $D$ that separates this character from the singular one. The two actions of $D$ on the corresponding Ext groups agree; their difference is a nonzero scalar plus a nilpotent operator, hence invertible. These Ext groups, and therefore cohomology with algebraic coefficients, vanish. This is the generalized-character version of the infinitesimal-character obstruction recalled in \cite[Ch.~I, Th.~5.3(ii)]{BorelWallach}. For residual cohomology at regular infinitesimal character, see \cite{GrbacSchwermer2010}. 
\end{remark}

If $d_P$ is odd, the involution of $\Psi_P$ has a fixed root, so simultaneous criticality requires the common wall. To compute the interval when $\mathcal F_P=\varnothing$, put
\begin{equation}\label{eq:delta}
 \Delta_{v,\beta}=q_\beta(b)-q_{w_M\beta}(b)
   =\langle\mu-w_M\mu+2\rho_M,\beta^\vee\rangle .
\end{equation}
The difference is independent of $b$. The entries of each pair have opposite signs exactly when the absolute value of their sum is smaller than that of their difference. By \eqref{eq:real-centering} the admissible critical shifts therefore satisfy
\[
 |b-b_v^\circ|<\min_{\beta\in\Psi_P}
                   \frac{|\Delta_{v,\beta}|}{2j(\beta)} .
\]

At a critical fixed wall, Proposition~\ref{prop:face} applies with vanishing roots $\mathcal F_P$. Put $n_{\rm fix}=|\mathcal F_P|$ and $n_{{\rm fix},j}=|\mathcal F_P\cap\Psi_{P,j}|$. Each nonfixed pair contributes one inversion on either side of the wall. All $n_{{\rm fix},j}$ fixed roots of grade $j$ change sign from negative to positive as $b$ increases. Thus the neighboring lengths are $(d_P+n_{\rm fix})/2$ on the left and $(d_P-n_{\rm fix})/2$ on the right, and in grade $j$ they are $(d_j\pm n_{{\rm fix},j})/2$. These are real chamber statements, not assertions about lattice points adjacent to the wall. 

\begin{example}[A second grade can force a wall]\label{ex:sp4-wall}
For the Siegel parabolic of $\Sp_4$, $\gamma_P=e_1+e_2$ and $f_P=3/2$. The exchanged roots $2e_i$ have coroots $e_i$ and grade one, whereas $e_1+e_2$ has grade two and carries the determinant on the dual Levi $\GL_2$. Take $\mu=(0,-1)$ and $b=-1$. Then $x=(1,-1)$, $s=-1/2$ and the coordinates are $(1,-1;0)$. The cohomological two-dimensional Weil parameter has exponent difference $\mu_1-\mu_2+1=2$. On its exchanged weight lines $J$ has determinant $-(-1)^2=-1$ (use $J^2=-1$), so the fixed determinant line is odd. Its parity is unchanged by this grade-two algebraic twist. Both orbits have arguments $1,2,2,1$. The tuple is critical, but the shift is singular and its neighboring chamber lengths are two and one. 
\end{example}

\subsection{The real Weil action and the endpoint tests}\label{subsec:real-weil}
We determine the real Weil action on the graded root spaces and translate adjacent criticality into root-sign and parity conditions, distinguishing exchanged root pairs from fixed root lines.

The real factor of $\sgn^\epsilon|\cdot|^h$ is 
\begin{equation}\label{eq:real-factor}
 L(t,\sgn^\epsilon|\cdot|^h)=\Gamma_\R(t+h+\epsilon),\qquad
 \Gamma_\R(t)=\pi^{-t/2}\Gamma(t/2).
\end{equation}
For $A,B\in\R$ with $A-B\in\Z$ and $A\ne B$, the induced Weil representation has factor \begin{equation}\label{eq:complex-factor}
 L\!\left(t,\operatorname{Ind}_{\C^\times}^{W_\R}
                 (z^A\bar z^B)\right)=\Gamma_\C(t+\max(A,B))
\end{equation}
\cite[(3.6), p.~404]{Knapp1994}.

For a fixed root $\beta$, write the character of its real Weil line as $\sgn^\epsilon|\cdot|^h$. Since $w_M\beta=\beta$, its exponent is $h=A_\beta=B_\beta$ in Proposition~\ref{prop:real-dictionary}. Thus \eqref{eq:exponents} gives $h+js=-q_\beta(b)$, and the four arguments in \eqref{eq:four-factors} are 
\begin{equation}\label{eq:real-arguments}
 -q_\beta(b)+\epsilon,\quad 1-q_\beta(b)+\epsilon,\quad 1+q_\beta(b)+\epsilon,\quad q_\beta(b)+\epsilon.
\end{equation}

\begin{lemma}[Real parameter alignment]\label{lem:alignment}
Under the hypotheses on $\sigma_v$ in Theorem~\ref{thm:real-criterion}, after conjugating into the dual Levi, the parameter satisfies
\[
 \phi_{\sigma_v}(z)=z^{\nu_M}\bar z^{w_M\nu_M},\qquad
 \nu_M=-w_M(\mu+\rho_M),\qquad z\in\C^\times.
\]
Moreover $\phi_{\sigma_v}(J)$ normalizes the chosen dual torus and induces $w_M$ on its roots. Thus it exchanges the dual root lines labeled by $(w_M\beta)^\vee$ and $\beta^\vee$.
\end{lemma}

\begin{proof}
The infinitesimal-character formula for real parameters \cite[\S2, Lemma~1]{NairPrasad2021} permits us to write $\phi_{\sigma_v}(z)=z^{\nu_M}\bar z^{\nu'}$, with the first exponent equal to the displayed strictly $M$-dominant representative. The second exponent is real, since $\nu_M-\nu'$ is an integral cocharacter. Regularity of $\nu_M$ implies that the centralizer of $\phi_{\sigma_v}(\C^\times)$ in the dual Levi is its chosen torus: no dual root is trivial on this image. The Weil relation implies that $\phi_{\sigma_v}(J)$ normalizes that torus. Let $v_0\in W_M$ be its Weyl action; the same relation gives $\nu'=v_0\nu_M$.

Essential temperedness makes the real radial exponent $\nu_M+\nu'$ central. Thus $\nu'$ is strictly antidominant on the Levi coroots. Since $\nu'=v_0\nu_M$ and $\nu_M$ is strictly dominant, uniqueness of the antidominant element of its Weyl orbit gives $\nu'=w_M\nu_M$. Regularity then implies $v_0=w_M$.
\end{proof}

\begin{proposition}[Real root--Gamma dictionary]\label{prop:real-dictionary}
Under the hypotheses on $\sigma_v$ in Theorem~\ref{thm:real-criterion}, on the dual root line labeled by $(w_M\beta)^\vee$, the restriction to $\C^\times$ has exponents
\[
 A_\beta=-\langle\mu+\rho_M,\beta^\vee\rangle,
 \qquad B_\beta=-\langle\mu+\rho_M,(w_M\beta)^\vee\rangle.
\]
At $s=k_P-b$ they satisfy
\begin{equation}\label{eq:exponents}
 A_\beta+js=-q_\beta(b),\qquad
 B_\beta+js=-q_{w_M\beta}(b).
\end{equation}
For a two-element orbit, substitute $q=q_\beta(b)$ and $q'=q_{w_M\beta}(b)$ in \eqref{eq:complex-arguments} to obtain the four $\Gamma_\C$ arguments. Each fixed root contributes the four real Gamma factors with arguments \eqref{eq:real-arguments}, where the sign exponent $\epsilon=\epsilon_\beta\in\{0,1\}$ is defined by the action $(-1)^{\epsilon_\beta}$ of $J$ on that line.
\end{proposition}

\begin{proof}
Pair the two exponents in Lemma~\ref{lem:alignment} with $(w_M\beta)^\vee$. Since $\langle\rho_P,\beta^\vee\rangle=jf_P$, adding $js=-jf_P-jb$ gives \eqref{eq:exponents}. On a two-element orbit the two characters are exchanged by $J$. When $A_\beta\ne B_\beta$, their span is the induced two-dimensional Weil representation of \eqref{eq:complex-factor}. When $A_\beta=B_\beta=h$, it instead splits into $|\cdot|^h$ and $\sgn|\cdot|^h$: the exchange matrix has eigenvalues $1,-1$, because $\phi_{\sigma_v}(J)^2$ acts as the character at $-1$, which is one. The identity 
\[
 \Gamma_\R(t+h)\Gamma_\R(t+h+1)=\Gamma_\C(t+h)
\]
then gives exactly the same product formula, also for the dual. On a fixed root line $A_\beta=B_\beta$, and its extension has one sign $\epsilon_\beta$. Apply \eqref{eq:real-factor} and the dual character formula. 
\end{proof}

\begin{lemma}[Criticality of either endpoint]\label{lem:parity}
Fix $b\in\Z$, put $s=k_P-b$, and consider a single grade $j$. For a complex Gamma factor use its paired coordinates $q,q'$ and put $q_{\min}=\min(q,q')$, $q_{\max}=\max(q,q')$. For the real Weil line of a fixed root $\beta$, the coordinate is $q_\beta(b)\in\Z$ and $\epsilon$ denotes the sign exponent. The entries give the conditions for each endpoint to be critical on that factor. 
\[
\begin{array}{lcc}
\toprule
\text{Factor} & js & js+1\\
\midrule
\Gamma_{\C}\ \text{(complex place or exchanged real pair)}
 & q_{\min}<0\le q_{\max} & q_{\min}\le0<q_{\max}\\
\Gamma_{\R},\ q_\beta(b)>0 & q_\beta(b)-\epsilon\ \text{odd} & q_\beta(b)-\epsilon\ \text{even}\\
\Gamma_{\R},\ q_\beta(b)<0 & q_\beta(b)-\epsilon\ \text{even} & q_\beta(b)-\epsilon\ \text{odd}\\
\Gamma_{\R},\ q_\beta(b)=0 & \epsilon=1 & \epsilon=1\\
\bottomrule
\end{array}
\]
An endpoint is critical for the grade at $v$ if it passes every factor at $v$. Global criticality requires these tests at all archimedean places, as in Definition~\ref{def:simultaneous}. On a fixed real line the two endpoints are critical together precisely when
\begin{equation}\label{eq:parity-condition}
 q_\beta(b)=0\quad\text{and}\quad\epsilon=1 .
\end{equation}
\end{lemma}
\begin{proof}
In \eqref{eq:complex-arguments} and \eqref{eq:real-arguments} the first and third factors test $js$, the second and fourth test $js+1$. At integral arguments $\Gamma_\C$ has poles at the nonpositive integers and $\Gamma_\R$ at the nonpositive even integers. For $q_\beta(b)>0$ the possible poles come from $-q_\beta(b)+\epsilon,1-q_\beta(b)+\epsilon$, respectively; for $q_\beta(b)<0$ they come from $1+q_\beta(b)+\epsilon,q_\beta(b)+\epsilon$. At $q_\beta(b)=0$, each endpoint tests $\epsilon$ and $1+\epsilon$.
\end{proof}

\begin{proof}[Proof of Theorem~\ref{thm:real-criterion}]
Proposition~\ref{prop:real-dictionary} gives one complex Gamma factor per exchanged pair and one real Gamma factor per fixed root in each local $L$-factor. The endpoint tests in Lemma~\ref{lem:parity} give the first equivalence. Suppose there are no fixed roots. Every orbit then has two elements. Opposite signs give exactly $d_P/2$ negative normal-root coordinates and no zero coordinate. The Levi coordinates of $x(b)$ are strictly positive. Proposition~\ref{prop:inverse-kostant}, applied to $\mu_b$, consequently gives a unique $u\in W^P$ and a dominant integral $\lambda$, with $\ell(u)=d_P/2$. Each grade has one inversion per orbit, giving $\ell_j(u)=d_j/2$.

Conversely, an ambient realization implies regularity, and \eqref{eq:length-sign} counts its length by negative normal-root coordinates. The sums in \eqref{eq:real-centering} have one common weak sign. Lemma~\ref{lem:common-sign}, applied to the $d_P/2$ two-element orbits, therefore forces opposite signs in every orbit when the length is $d_P/2$. This proves criticality. 

If there is a fixed root, its pole test forces $q_\beta(b)=0$ and $\epsilon_\beta(b)=1$. Since every grade is positive, \eqref{eq:fixed-wall} forces $b=b_v^\circ$. At this point a nonfixed pair sums to zero by \eqref{eq:real-centering}, so its entries have opposite signs exactly when $q_\beta(b_v^\circ)\ne0$. This proves \eqref{eq:singleton}. Finally, a zero root coordinate makes $x(b_v^\circ)$ singular. Every $\lambda+\rho$ for dominant integral $\lambda$ is regular, so its Weyl orbit cannot contain $x(b_v^\circ)+z_{b_v^\circ}$. 
\end{proof}

In particular, for $s\in k_P-\Z$ all Gamma arguments at a real place are integers. Definition~\ref{def:simultaneous} is then the local automorphic analogue of Deligne's criticality criterion \cite[Def.~1.3]{Deligne1979} at $js$ and $js+1$. 

\begin{remark}[Algebraic and absolute-value twists]\label{rem:twist-families}
Use the character $\theta_b$ of Definition~\ref{def:admissible-shift}. At $k_P$, both $\sigma_v\otimes\chi_{-b}$ and $\sigma_v\otimes\theta_b^{-1}$ give the same root coordinates as $\sigma_v$ at $s=k_P-b$; their signs need not coincide. The absolute-value twist is trivial on $J$. By Remark~\ref{rem:central-correction-parameter}, an algebraic integral twist carries grade-$j$ sign $(-1)^{jb}$, so it changes $\epsilon_\beta$ by $jb\pmod2$. In particular,
\[
 \epsilon_\beta^{\rm abs}(b)=\epsilon_\beta(0),\qquad
 \epsilon_\beta^{\rm alg}(b)\equiv\epsilon_\beta(0)+j(\beta)b\pmod2 .
\]
Simultaneous inverse algebraic twisting always cancels the change in the coefficient tensor, including its component character; absolute-value twisting need not. The real pole test always uses the sign of the specified family. 
\end{remark}

For a grade containing a single fixed line, its sign can be read from $\det(r_j(\phi_{\sigma_v}(J)))$. On each exchanged pair $O$ with exponent difference $m_O=A_\beta-B_\beta$, the determinant of $J$ is $(-1)^{m_O+1}$ since $J^2=-1$. Dividing the total determinant by the product of these signs gives the fixed-line sign. With several fixed lines, this calculation gives the product of their signs.

\subsection{The general-linear real parabolic}\label{subsec:gl-real}
Let $G=\GL_N$, $N=n+n'$, and $M=\GL_n\times\GL_{n'}$. Suppress the real embedding and write $\mu=(\mu_1,\ldots,\mu_n;\mu'_1,\ldots,\mu'_{n'})$. The normal roots and their coordinates are \begin{equation}\label{eq:gl-coordinates}
 \beta_{it}=e_i-e_{n+t},\qquad
 q_{it}(b)=\mu_i-\mu'_t+n+t-i+b.
\end{equation}
Every grade is one, $d_P=nn'$, and 
\[
 \gamma_P=(\underbrace{n'/N,\ldots,n'/N}_{n};
            \underbrace{-n/N,\ldots,-n/N}_{n'}),\qquad
 f_P=N/2.
\]
Although $\gamma_P$ need not be integral, adding a $G$-central correction makes the twist add $b$ to the first block and zero to the second. Conversely the omitted simple-coroot pairing forces $b\in\Z$, so $\Lambda_{\mathrm{alg}}=\Z$. 

Suppose the two block purity weights are $\mathrm{w}_1,\mathrm{w}_2$. Then
\[
 a_v=\mathrm{w}_1-\mathrm{w}_2,\qquad
 T_v(b)=\mathrm{w}_1-\mathrm{w}_2+N+2b,\qquad
 w_M\beta_{it}=\beta_{n+1-i,n'+1-t}.
\]

The block reversals fix the normal index $(i,t)$ exactly when both blocks have a middle entry. If $nn'$ is even, there are no fixed roots: the condition is an ambient realization of length $nn'/2$. For strongly inner data with disjoint cuspidal parameters over a totally real field, applying this criterion at every real place with $b=0$ recovers the combinatorial lemma of Harder--Raghuram, proved by Weselmann \cite[Lemma~7.14 and \S\S~7.2.4.]{HarderRaghuram2020}.

For two odd blocks the middle character of each real parameter is a line. Their tensor quotient has the sum of their parities modulo two. Criticality then requires 
\[
 b=-(\mathrm{w}_1-\mathrm{w}_2+N)/2\in\Z,\quad
 \epsilon_{\rm middle}+\epsilon'_{\rm middle}\equiv1\pmod2.
\]
Each of the remaining pairs must have nonzero gap. At this shift, \(x(b)=\mu+\rho+b\gamma_P\) is singular because its middle root coordinate vanishes. The middle signs here are those for the actual fixed representation at the shifted spectral argument. For an algebraic twist at $k_P$ use the transported signs. This recovers the odd--odd criticality and nonrealization statements of Bhagwat--Raghuram \cite[Lemma~7.14 and \S7.2.4]{HarderRaghuram2020}. Away from the wall Corollary~\ref{cor:regular-boundary-window} gives the common boundary degrees; see Example~\ref{ex:gl4-unbalanced} for $(3,1)$. 

For real pure weights it is useful to center each block separately:
\[
 A_i=\mu_i-\frac{\mathrm{w}_1}{2}+\frac{n+1-2i}{2},\qquad
 A'_t=\mu'_t-\frac{\mathrm{w}_2}{2}+\frac{n'+1-2t}{2}.
\]
Then $A_{n+1-i}=-A_i$, $A'_{n'+1-t}=-A'_t$, and 
\[
 q_{it}=A_i-A'_t+\frac{T_v}{2},\qquad
 \Delta_{v,it}=2(A_i-A'_t).
\]
When no fixed root occurs, the criterion is $|T_v|<\min_{i,t}2|A_i-A'_t|$. For a specified pair of block weights, a zero gap makes this interval empty even when no root is fixed.

In the totally real setting of Harder--Raghuram, write $A_{i,v}$ and $A'_{t,v}$ for the centered coordinates at $v$. The cuspidal parameters of the two blocks are $2A_{i,v}$ and $2A'_{t,v}$ \cite[\S7.4.1]{HarderRaghuram2020}. The cuspidal width, denoted there by $\ell(\mu,\mu')$, is distinct from the Kostant length $\ell(w)$. By \cite[Def.~7.4]{HarderRaghuram2020}, it equals
\[
  \ell(\mu,\mu')
    =\min_{v\mid\infty}\ \min_{\substack{1\le i\le n\\1\le t\le n'}}
        |\Delta_{v,it}|.
\]

The Hodge calculation gives the same minimum. The differences $p-q$ in the tensor-product Hodge types are sums of cuspidal parameters. Since $A'_{n'+1-t,v}=-A'_{t,v}$, reversing the second block identifies the smallest $|p-q|$ with the smallest $|\Delta_{v,it}|$. When the parameters are disjoint, this width equals $q_{\min}-p_{\max}$ \cite[\S\S7.1.3 and~7.2.2.1]{HarderRaghuram2020}. For the general-linear combinatorial lemma over a totally imaginary field, see Corollary~\ref{cor:Raghuram}.

\begin{example}[$\GL_3$ and two admissible twists]
Take $P=(2,1)$ and $\mu=(0,-2;0)$ over $\R$. The two coordinates are $2+b,-1+b$, hence the critical interval is $(-2,1)$ and its lattice points are $b=-1,0$. Both lengths are one. For $b=0$, the dominant ambient weight is $(0,-1,-1)$. For $b=-1$, choose the central correction making $\mu_b=(-1,-3;0)$; its dominant ambient weight is $(-1,-1,-2)$. Each follows by sorting $\mu_b+\rho$ and subtracting $\rho=(1,0,-1)$. Selecting the common Levi bottom degree one gives total boundary degree two, the degree identity of Section~\ref{sec:cohomological-degrees}. 
\end{example}

\subsection{The Klingen wall and two classical comparisons}\label{subsec:classical-tests}

\emph{Klingen.} On $\Sp_{2n}$ ($n\ge2$), let $M=\GL_1\times\Sp_{2n-2}$ and $\mu=(u;\mu_2,\ldots,\mu_n)$. All grades are one; $\gamma_P=e_1$, $f_P=n$, $d_P=2n-1$. Set $z=u+n+b$ and $t_i=\mu_i+n+1-i>0$. The fixed root $2e_1$ has coordinate $q_{2e_1}(b)=z$, and the exchanged roots $e_1\pm e_i$ have coordinates $z\pm t_i$. Take the cohomological discrete-series parameter of the symplectic factor and the character $\sgn^\delta|\cdot|^{-u}$, $\delta\equiv u\bmod2$. The standard orthogonal parameter has determinant one; its $n-1$ two-dimensional blocks have exponents $2t_i$ and determinant $-1$ at $J$. Its remaining line therefore has sign $(-1)^{n-1}$. After tensoring with the character, the fixed sign is $(-1)^{u+n-1}$. The wall $b_v^\circ=-u-n$ is critical for the absolute-value family exactly when $u+n$ is even and for the algebraic family for every $u$. The other gaps are nonzero and the neighboring lengths are $n,n-1$ (defects $1,-1$).

At a complex place the same root $2e_1$ has gap $d=u^\eta-u^{\bar\eta}$; the others have gaps $d\pm2t_i$. The center is $-n-(u^\eta+u^{\bar\eta})/2$ in the $b$ coordinate. 

\begin{proposition}[A critical Klingen wall for global data]\label{prop:klingen-wall}
Let $F$ be totally real, $n\ge2$, and let $P$ be the Klingen parabolic of $\Sp_{2n}$, with Levi $M=\GL_1\times\Sp_{2n-2}$. At each real place $v$, fix a dominant integral symplectic weight $\mu'_v$ whose cohomological discrete-series packet contains an integrable generic member. Then there is a globally generic cuspidal automorphic representation of $M(\A_F)$ which at every real place $v$ is essentially tempered and cohomological with weight $\mu^v=(-n;\mu'_v)$, and for which $b=0$ is the unique critical twist at every real place for both twist families. At this twist $x(0)$ is singular. Sufficiently deep weights $\mu'_v$ satisfy the hypothesis \cite{Muic2016}. \end{proposition}

\begin{proof}
The group $\mathbf H=\operatorname{Res}_{F/\Q}\Sp_{2n-2}$ is semisimple and quasi-split over $\Q$, and $\mathbf H(\R)=\prod_{v\mid\infty}\Sp_{2n-2}(\R)$ is connected with compact Cartan. Genericity and integrability pass to the external tensor product of the chosen real representations. Weak approximation in $F$ lets us prescribe the real signs of the simple-root coefficients of a global Whittaker character so that each chosen member is generic for its localization. Mui\'c's globalization theorem, which uses the method of Khare--Larsen--Savin, now gives a globally generic cuspidal representation of $\mathbf H(\mathbb A_\Q)$ with the prescribed archimedean component \cite{Muic2016,KhareLarsenSavin2008}. The selected generic discrete series at infinity have cohomology with the chosen coefficients \cite[Ch.~II]{BorelWallach}.

On $\GL_1$ take $\omega_{F(i)/F}^{\,n}|\cdot|_{\A_F}^{\,n}$, where $\omega_{F(i)/F}$ is the quadratic Hecke character. At every real place this is $\sgn^n|\cdot|^n$, so $u=-n$ and $\delta\equiv n\bmod 2$, and it cancels the coefficient character. The character is trivial on the Levi unipotent used for the Whittaker coefficient.

At every real place $\mu^v+w_M\mu^v=(-2n;0,\ldots,0)$, so $b_v^\circ=0$, which is admissible. The fixed sign is $(-1)^{u+n-1}=-1$, and the other root coordinates at the wall are $\pm t_i$ with $t_i>0$. The two twist families have the same signs at $b=0$. Theorem~\ref{thm:real-criterion}(b) now applies at every real place.
\end{proof}

\emph{Siegel.} For the Siegel parabolic of $\Sp_{2m}$ put $t_i=\mu_i+m+1-i$. The roots $2e_i$ have grade one and coordinate $t_i+b$; the roots $e_i+e_k$ ($i<k$) have grade two and coordinate $t_i+t_k+2b$. Here $r_1=\Std_m$, $r_2=\Lambda^2\Std_m$ when nonzero, $d_P=m(m+1)/2$, and $f_P=(m+1)/2$. The longest Levi element reverses indices, so its fixed roots are the mirror sums and, for odd $m$, the single middle diagonal. Suppose $\mu_i+\mu_{m+1-i}=\mathrm{w}$ and take an essentially tempered cohomological parameter with full-$K$ central character compatible with the coefficient. At the wall $b_v^\circ=-(\mathrm{w}+m+1)/2$ the numbers $t_i+b_v^\circ$ strictly decrease and are antisymmetric: precisely the fixed coordinates vanish. Each mirror root $e_i+e_{m+1-i}$ with $i<m+1-i$ corresponds to the determinant line of the exchanged block with angular exponent $t_i-t_{m+1-i}$. Since $J^2=-1$, the determinant at $J$ is $(-1)^{t_i-t_{m+1-i}+1}=(-1)^{\mathrm{w}+m}$ by the purity relation. For even $m$, admissibility and odd parity therefore both require $\mathrm{w}$ odd; exactly then $b_v^\circ$ works for either twist family.

For $m=2h+1$, write $c=\mu_{h+1}=\mathrm{w}/2$. Then $b_v^\circ=-c-h-1$ is integral and all mirror lines are odd. Let $(-1)^{\epsilon_0}$ be the action of $J$ on the middle Weil line of the standard parameter. Each of the other $h$ two-dimensional blocks has determinant $-1$ at $J$, so the determinant of the standard parameter at $J$ is $(-1)^{\epsilon_0+h}$. Compatibility of the central character with the coefficient at $-I$ gives $\epsilon_0+h\equiv\sum_i\mu_i\equiv c\pmod2$. Thus the absolute-value family is critical exactly when $c+h$ is odd; the algebraic twist adds $b_v^\circ$ to the middle sign, making it odd for every $c$. All critical walls here are singular.

\emph{Connected even orthogonal group.} For the line parabolic of $\SO_{2n+2}$, $M=\GL_1\times\SO_{2n}$ ($n\ge2$). Put $\mu=(u;\mu_1,\ldots,\mu_n)$, $z=u+n+b$, $t_i=\mu_i+n-i$. All $2n$ coordinates are $z\pm t_i$ in grade one, with $t_1>\cdots>t_{n-1}>|t_n|$. For even $n$ the longest element of $D_n$ negates all $t_i$; hence criticality is $|z|<\min_i|t_i|=|\mu_n|$, equivalently an ambient realization of length $n$. There are $\max(0,2|\mu_n|-1)$ integral twists; a zero gap can therefore exclude all points despite the absence of fixed roots.

For even $n$, Bhagwat and Raghuram prove the combinatorial lemma for the isotropic-line parabolic of the split group $\mathrm O(2n+2)$ \cite[\S5.3]{BhagwatRaghuram2025}. The inequality $|z|<|\mu_n|$ recovers the root-theoretic part of their criterion. The comparison for the disconnected group also involves component and packet data. 

For odd $n$ the last coordinate is fixed by the Levi action and purity forces $\mu_n=0$. The two fixed roots have coordinate $z$: both signs must be odd at $z=0$, and the neighboring lengths are $n+1,n-1$. 

\subsection{The two \texorpdfstring{$G_2$}{G2} parabolics at a real place}\label{sec:g2-real-endpoints}
Write $\alpha$ for the short root and $\beta$ for the long root, and label the parabolic by the retained root. The unique fixed root is $2\alpha+\beta$ for $P_\beta$ and $3\alpha+2\beta$ for $P_\alpha$. In both cases its grade is two and its Weil line is the determinant of the standard $\GL_2$ parameter. For Levi simple-coroot pairing $n\ge0$, the angular difference is $n+1$, so $J$ acts on the determinant by $(-1)^n$. Writing the omitted simple-coroot pairing as $a$, the fixed coordinates are $2a+3n+5+2b$ and $2a+n+3+2b$, respectively. Their zeros lie in $\Lambda_{\rm alg}=\Z$ precisely when $n$ is odd. At such a zero the determinant is odd for either twist family. The exchanged pairs have absolute gaps $3(n+1),n+1$ for $P_\beta$ and $n+1,n+1$ for $P_\alpha$. Thus, for either maximal parabolic and either twist family, Theorem~\ref{thm:real-criterion} gives an empty critical set when $n$ is even and a single critical shift $b_v^\circ$ when $n$ is odd. At that shift the ambient infinitesimal character is singular.

Over a totally real field, trivial inducing central character gives the factor $\zeta_F(2s)$ in grade two. At $b=0$ its adjacent arguments are $(-5,-4)$ for $P_\beta$ and $(-3,-2)$ for $P_\alpha$. The negative odd argument is critical \cite[Prop.~23(i)]{RaghuramHecke2022}. At the negative even argument \(-2m\), the archimedean factor \(\Gamma_{\mathbb R}(s)^{[F:\Q]}\) of \(\zeta_F(s)\) has a pole of order \([F:\Q]\). The reflected value \(\zeta_F(1+2m)\) is nonzero by the Euler product, so the functional equation gives a zero of exactly that order for the finite zeta function at \(-2m\) \cite[\S4.2.4]{RaghuramHecke2022}. 

\section{Complex places: the critical interval and its converse} \label{sec:complex-criterion}
At complex places, we characterize the algebraic coefficient shifts for which the graded normalizing tuple has adjacent critical arguments. We prove that this condition is equivalent to ambient Kostant realizations with complementary total degrees, and hence to balance in every grade. We then determine the critical shifts explicitly as the algebraic lattice points
of an interval on the coefficient line.

Let $v=\{\eta,\bar\eta\}$ be a complex place. Throughout this section $\sigma_v$ is an irreducible, essentially tempered representation of $M(F_v)$, cohomological for the integral $M$-dominant coefficient pair $(\mu^\eta,\mu^{\bar\eta})$ in the convention $\sigma_v\otimes\cM_{\mu,v}$, where $\cM_{\mu,v}=\cM^M_{\mu^\eta}\boxtimes\cM^M_{\mu^{\bar\eta}}$ on the two complexified embedding factors. Since $\sigma_v$ is essentially unitary, Lemma~\ref{lem:essential-unitarity-purity} gives \eqref{eq:local-purity}, with $c_v$ and $a_v$ as there.

For $b\in\Lambda_{\rm alg}$ put $s=k_P-b$ and $x^\xi(b)=\mu^\xi+\rho+b\gamma_P$, and let
\[
 \mu_b^\xi=\mu^\xi+\theta_b=x^\xi(b)-\rho+z_b\in X
 \qquad(\xi=\eta,\bar\eta)
\]
be the weights of Definition~\ref{def:admissible-shift}. We refer to these as the \emph{local hypotheses at $v$}.

The opposite-sign test in the next theorem uses the four Gamma arguments of \eqref{eq:complex-arguments}.

\begin{theorem}[Combinatorial lemma at complex places] \label{thm:complex}
Let $v=\{\eta,\bar\eta\}$ be a complex place and $\sigma_v$ an irreducible, essentially tempered representation of $M(F_v)$, cohomological for the integral $M$-dominant coefficient pair $(\mu^\eta,\mu^{\bar\eta})$. Let $b\in\Lambda_{\mathrm{alg}}$, with
$\mu_b^\xi$ as in Definition~\ref{def:admissible-shift}. The following are equivalent:
\begin{enumerate}
\item $s=k_P-b$ is simultaneously critical at $v$.
\item For every $\beta\in\Psi_P$, $q^\eta_\beta(b)q^{\bar\eta}_{w_M\beta}(b)<0$.
\item At both $\xi=\eta,\bar\eta$ the integral Levi weight admits a realization $\mu_b^\xi=u^\xi\cdot\lambda^\xi$ with $u^\xi\in W^P$ and $\lambda^\xi\in X_G^+$, and the two Kostant degrees satisfy $\ell(u^\eta)+\ell(u^{\bar\eta})=d_P$.
\end{enumerate}
When these equivalent conditions hold, the chamber representatives $u^\eta,u^{\bar\eta}$ are unique,
\begin{equation}\label{eq:complementary-inversions}
 \Inv_P(u^{\bar\eta})
     =w_M\bigl(\Psi_P\setminus\Inv_P(u^\eta)\bigr),
 \qquad\text{equivalently}\qquad
 u^{\bar\eta}=w_Mu^\eta w_G,
\end{equation}
and hence
\begin{equation}\label{eq:graded-balance}
 \ell_j(u^\eta)+\ell_j(u^{\bar\eta})=d_j
                       \qquad\text{for every grade }j.
\end{equation}
\end{theorem}

\begin{remark}[Relation to Raghuram and to Raghuram--Zou] \label{rem:complex-attribution}
For $\GL_N$, Corollary~\ref{cor:Raghuram} recovers Raghuram's combinatorial lemma under the corresponding hypotheses and normalizations \cite[Lemma~3.16]{Raghuram2025}. Lemma~\ref{lem:common-sign} gives the form of Saper's paired-root count used here, with $(\eta,\beta)$ paired with $(\bar\eta,w_M\beta)$ \cite[Lemmas~8.7(iii) and~17.6]{SaperMicro}.

At $s_0=k_P$ ($b=0$ here), Raghuram and Zou independently deduce complementary Kostant representatives from hyper-criticality \cite[Lemma~7.15 and (7.17)--(7.18)]{RaghuramZou2026}. For $\alpha=w_M\beta$, the weights they denote by $\xi$ and $\xi'$ satisfy
\[
 (\langle\xi,\alpha^\vee\rangle,\langle\xi',\alpha^\vee\rangle)
 =(-q^\eta_\beta(0),-q^{\bar\eta}_{w_M\beta}(0))
\]
by \cite[(7.7), (7.9)--(7.13)]{RaghuramZou2026}, and both $js_0$ and $js_0+1$ lie in their critical-point lattice $c_\alpha+\Z$ \cite[Def.~5.8 and (7.14)]{RaghuramZou2026}. Under their hypotheses, evaluating the $\theta_b^{-1}$-twisted datum at $s_0$ amounts to evaluating the fixed datum at $s=k_P-b$ \cite[(5.6)]{RaghuramZou2026}. Theorem~\ref{thm:complex} also proves the converse at every admissible $b$ by deriving the opposite-sign inequalities from total degree balance. Proposition~\ref{prop:interval} identifies the critical shifts as the admissible points of an explicit open interval on the coefficient line. 
\end{remark}

\subsection{From the parameter to a sign test}
Essential temperedness determines how the two infinitesimal characters are paired in the archimedean Artin factors.
\begin{lemma}\label{lem:complex-alignment}
Write $\nu_M^\xi=-w_M(\mu^\xi+\rho_M)$, as in Lemma~\ref{lem:essential-unitarity-purity}. Up to conjugation in $\widehat M$, the Langlands parameter of $\sigma_v$ on $W_{F_v}=\C^\times$ is 
\[
 z\longmapsto z^{\nu_M^\eta}\,\bar z^{\,w_M\nu_M^{\bar\eta}},
 \qquad
 \nu_M^\eta-w_M\nu_M^{\bar\eta}
   =\mu^{\bar\eta}-w_M\mu^\eta+2\rho_M\in X.
\]
\end{lemma}
\begin{proof}
The infinitesimal character fixes the two exponents up to independent $W_M$-orbits \cite[\S11, Cor.~3]{NairPrasad2021}. After conjugation the first is $\nu_M^\eta$ and the second is $v_0\nu_M^{\bar\eta}$ for some $v_0\in W_M$; both are real in the cohomological convention.

Essential temperedness makes the radial exponent $\nu_M^\eta+v_0\nu_M^{\bar\eta}$ Levi-central. Hence $-v_0\nu_M^{\bar\eta}$ has the Levi coroot pairings of $\nu_M^\eta$ and is strictly $M$-dominant. So is $-w_M\nu_M^{\bar\eta}$, and both lie in the $W_M$-orbit of $-\nu_M^{\bar\eta}$; they are therefore equal, and regularity gives $v_0=w_M$. The displayed difference is computed from $w_M\rho_M=-\rho_M$. 
\end{proof}

In particular $\nu_M^\eta+w_M\nu_M^{\bar\eta}=-c_v$, which is local purity again. 

\begin{proposition}[Root--Gamma dictionary]\label{prop:complex-dictionary}
Let $b\in\mathbb R$ and $s=k_P-b$. Fix $\beta\in\Psi_{P,j}$ and write $(q,q')=(q^\eta_\beta(b),q^{\bar\eta}_{w_M\beta}(b))$. On the line $\widehat{\mathfrak g}_{(w_M\beta)^\vee}$ of $V_j$, the representation $r_j\circ\phi_{\sigma_v}$ of $\C^\times$ acts by $z^A\bar z^B$ with 
\[
 A=jf_P-q^\eta_\beta(0),\qquad B=jf_P-q^{\bar\eta}_{w_M\beta}(0).
\]
The contributions of this line to \eqref{eq:four-factors} are the four $\Gamma_\C$ factors of \eqref{eq:complex-arguments}. If $b\in\Z$, then $q$ and $q'$ are integers, and all four are finite if and only if $qq'<0$.
\end{proposition}
\begin{proof}
Pair the exponents of Lemma~\ref{lem:complex-alignment} with $(w_M\beta)^\vee$ and use $\langle\rho_P,\beta^\vee\rangle=jf_P$. This gives the stated exponents and the identities $A+js=-q$ and $B+js=-q'$.

The local factor of $z^A\bar z^B$ is $\Gamma_\C(t+\max(A,B))$, where $\Gamma_\C(t)=2(2\pi)^{-t}\Gamma(t)$ (see \cite[p.~13]{Tate1979} and \cite[(4.6)]{Knapp1994}). Applying this formula to the character and its inverse gives the four Gamma arguments. For $b\in\Z$ they are integers, all positive precisely when $\min(q,q')<0<\max(q,q')$. A pole cannot cancel in a product of Gamma functions, which have no zeros.
\end{proof}

\subsection{Proof of the complex combinatorial lemma}
\begin{proof}[Proof of Theorem~\ref{thm:complex}]
Proposition~\ref{prop:complex-dictionary} gives $(1)\Leftrightarrow(2)$. Under (2), strict Levi dominance makes both shifted weights regular. Proposition~\ref{prop:inverse-kostant} gives the realizations, and the one inversion in each pair gives (3) and the first form of \eqref{eq:complementary-inversions}. 

Under (3), both shifted vectors are regular. To prove (2), count inversions in pairs as in Saper's argument \cite[Lemma~17.6]{SaperMicro}. By \eqref{eq:length-sign}, their total length is the number of negative coordinates in the $d_P$ pairs 
\[
 (q^\eta_\beta(b),q^{\bar\eta}_{w_M\beta}(b)),\qquad\beta\in\Psi_P.
\]
By \eqref{eq:root-centering}, the paired sums are $2j(\beta)(b-b_v^\circ)$. Their signs agree because every grade is positive. Lemma~\ref{lem:common-sign} therefore forces opposite signs in every pair from the total of $d_P$ negative entries, proving (2). At $b=b_v^\circ$ each pair sums to zero and its entries are nonzero by regularity.

For $u\in W^P$ the element $w_Muw_G$ lies in $W^P$, and $\beta\in\Psi_P$ is an inversion of it exactly when $w_M\beta$ is not an inversion of $u$; since an element of $W^P$ is determined by its inversion set, this gives the second form of \eqref{eq:complementary-inversions}. Finally $w_M$ preserves grades, so the complementary inversion sets give \eqref{eq:graded-balance}.
\end{proof}

Without the common-sign hypothesis, equal positive numbers of entirely negative and entirely positive pairs can cancel in the total count, giving a balanced but noncomplementary family. 

For example, on $\GL_3$ with $P=B$, which is not maximal, take the two shifted weights $(0,1,-1)$ and $(-1,1,0)$, both permutations of $\rho$. On $e_1-e_2,e_2-e_3,e_1-e_3$ the pairs are $(-1,-2),(2,1),(1,-1)$. The two Kostant lengths are $1$ and $2$, so their sum equals the number of normal roots. Nevertheless, the first pair is entirely negative and the second entirely positive: their contributions cancel in the total count. More generally, Lemma~\ref{lem:common-sign} applies whenever no paired coordinate vanishes and the paired sums are either all nonnegative or all nonpositive. Under these hypotheses, total degree balance forces opposite signs in every pair. 

In higher relative rank this is an additional chamber condition on the relative central parameter: its pairings with all positive normal coroots must lie in the same closed half-line. In relative rank one purity supplies that condition in \eqref{eq:root-centering}. Even for a maximal parabolic, without the integrality of Proposition~\ref{prop:complex-dictionary} the pair of coordinates $(q,q')=(1/2,1/2)$ passes all four Gamma tests despite having equal signs.

\subsection{The critical interval and admissible shifts}
We now determine the shifts satisfying the preceding criterion. The graded root inequalities define an explicit interval, whose intersection with $\Lambda_{\rm alg}$ gives the critical algebraic shifts. Intersecting these intervals over the complex places determines whether a common shift exists. For split $G_2$, the critical range in each grade and its variation under Tate twists are computed in \cite[\S\S5.1--5.2 and \S6.3]{HosseiniJafariG2}. The normal-root coordinates give uniform bounds for all maximal parabolics.

The difference
\begin{equation}\label{eq:gap}
 \Delta_{v,\beta}=q^\eta_\beta(b)-q^{\bar\eta}_{w_M\beta}(b)
 =\langle\mu^\eta-w_M\mu^{\bar\eta}+2\rho_M,\beta^\vee\rangle
\end{equation}
is independent of $b$ and equals $B-A$. By \eqref{eq:root-centering},
\begin{equation}\label{eq:product}
 4qq'=4j(\beta)^2(b-b_v^\circ)^2-\Delta_{v,\beta}^2.
\end{equation}
Consequently the entries have opposite signs exactly when $2j(\beta)|b-b_v^\circ|<|\Delta_{v,\beta}|$. For integral $b$ the sum and difference of the two coordinates have the same parity, so a strict gap between their absolute values is at least two. Set
\begin{equation}\label{eq:radius}
 R_v=\min_{\beta\in\Psi_P}
           \frac{|\Delta_{v,\beta}|}{j(\beta)}.
\end{equation}
The critical width $R_v$ is the minimum of the widths in the individual grades. The interval has radius $R_v/2$.

\begin{remark}[Roots fixed by Levi opposition]\label{rem:fixed-root}
If $\beta\in\Psi_P$ satisfies $w_M\beta=\beta$, then $\langle\rho_M,\beta^\vee\rangle=0$ and $\Delta_{v,\beta}=\langle\mu^\eta-\mu^{\bar\eta},\beta^\vee\rangle$. Thus, when such a root exists, equal conjugate weights $\mu^\eta=\mu^{\bar\eta}$ give $R_v=0$ and admit no critical twist at $v$.

This obstruction occurs for the Siegel and Klingen parabolics of $\Sp_{2n}$ and for both maximal parabolics of $G_2$. For Klingen the fixed root $2e_1$ has gap $u^\eta-u^{\bar\eta}$, so equality of the two $\GL_1$ weights already forces $R_v=0$ (Section~\ref{subsec:classical-tests}).
\end{remark}

\begin{proposition}[Critical interval]\label{prop:interval}
The real parameters for which all paired root coordinates have opposite signs form
\begin{equation}\label{eq:critical-interval}
 \mathcal I_v=(b_v^\circ-R_v/2,\,b_v^\circ+R_v/2).
\end{equation}
This is empty when $R_v=0$. The algebraically admissible critical shifts at $v$ are exactly $\mathcal I_v\cap\Lambda_{\mathrm{alg}}$. In the spectral variable, the same interval is 
\begin{equation}\label{eq:spectral-interval}
 \frac{a_v-R_v}{2}<s<\frac{a_v+R_v}{2}.
\end{equation}
On a nonempty $\mathcal I_v$, the two Kostant representatives $u^\eta,u^{\bar\eta}$ are constant, and $\mathcal I_v=I^\eta_{u^\eta}\cap I^{\bar\eta}_{u^{\bar\eta}}$ in the notation of Proposition~\ref{prop:intervals}; its endpoints are wall-crossing parameters \eqref{eq:wall-parameters} of one of the two embeddings. 
\end{proposition}

\begin{proof}
By \eqref{eq:product}, opposite signs in all pairs are equivalent to $2|b-b_v^\circ|<R_v$. This gives \eqref{eq:critical-interval}, and the substitution $s=k_P-b$ gives \eqref{eq:spectral-interval}. The theorem proves the assertion on the algebraic lattice. No root coordinate can vanish inside this interval, so its sign is constant there. The chamber criterion then gives constant representatives, so $\mathcal I_v$ lies in both chamber intervals. On their intersection every root coordinate has the sign it has on $\mathcal I_v$, which gives equality. At an endpoint some product $qq'$ vanishes, so one coordinate does.
\end{proof}

For $\Lambda_{\rm alg}=d_{\rm alg}\Z$ the number of lattice points in an open interval $(b_-,b_+)$ is 
\[
 \max\{0,\lceil b_+/d_{\rm alg}\rceil-\lfloor b_-/d_{\rm alg}\rfloor-1\};
\]
when nonempty, the extreme indices $n=b/d_{\rm alg}$ are $\lfloor b_-/d_{\rm alg}\rfloor+1$ and $\lceil b_+/d_{\rm alg}\rceil-1$.

Over a totally imaginary field take $b_-=\max_v(b_v^\circ-R_v/2)$ and $b_+=\min_v(b_v^\circ+R_v/2)$. A zero local width, disjoint open intervals, or an intersection missing the lattice therefore rules out a common twist. The count concerns coefficients; cuspidal existence is not part of it.

\begin{example}[A higher grade narrows the complex interval]\label{ex:sp4-complex}
For the Siegel parabolic of $\Sp_4$, use the normal roots $(2e_1,e_1+e_2,2e_2)$ of grades $(1,2,1)$, and $\mu^\eta=(-1,-3)$, $\mu^{\bar\eta}=(0,-2)$. Here $\rho=(2,1)$, $\gamma_P=(1,1)$, $f_P=3/2$, and $c_v=(-3,-3)$, so $a_v=-3$. At $b=0$ the paired coordinates are $(1,-1),(-1,1),(-2,2)$. The widths in grades one and two are $2$ and $1$, respectively, and the center $b_v^\circ=0$ gives the interval $(-1/2,1/2)$ and the unique lattice point $b=0$.

The two shifted vectors are signed permutations of $\rho$, so the two ambient coefficients are trivial and the lengths $2,1$ sum to $3$. The root fixed by Levi opposition at the middle of the list pairs two distinct embeddings here (at a real place its Weil character enters the parity test of Lemma~\ref{lem:parity}).
\end{example}

\begin{example}[A non-general-linear Levi with odd $d_P$] \label{ex:Sp6}
For the Klingen parabolic of $\Sp_6$, take $M=\GL_1\times\Sp_4$, $\gamma_P=e_1$, $f_P=3$, and $d_P=5$. All normal roots have grade one. In the root order
\[
 e_1-e_2,\ e_1-e_3,\ 2e_1,\ e_1+e_3,\ e_1+e_2,
\]
choose the coefficient pair
\[
 \mu^\eta=(-2;1,1),\qquad\mu^{\bar\eta}=(-4;1,1).
\]
Then $c_v=-6e_1$, $a_v=-6$, and $b_v^\circ=0$. The coordinates at $b=0$ are 
\[
 q^\eta=(-2,-1,1,3,4),\qquad
 (q^{\bar\eta}_{w_M\beta})_\beta=(2,1,-1,-3,-4).
\]
Thus $\Delta_v=(-4,-2,2,6,8)$, $R_v=2$, and $\mathcal I_v=(-1,1)$. Exactly $b=0$ lies in its algebraic lattice $\Z$. At $s=-3$, the selected arguments are $-3,-2$. The degrees are $2+3=5$, and both ambient weights are zero. This realizes the odd-dimensional complex-place case without a general-linear assumption on the derived Levi. The calculation applies to essentially tempered cohomological inducing data with these weights. 
\end{example}

The root coordinates and their chambers are defined for real $b$. For integral $b$, the four Gamma factors are finite precisely when the paired coordinates have opposite signs. An integral ambient Kostant realization also requires $b\in\Lambda_{\rm alg}$, as the following example shows. 

On $\PGL_2$ write $\omega=\alpha/2$ and take $\mu^\eta=0$, $\mu^{\bar\eta}=-4\omega$. At $b=1$ the paired coordinates are $2,-2$ and the four Gamma arguments $2,3,3,2$. But the shifted coefficients $\omega,-3\omega$ are outside $X=2\Z\omega$; there is no $G$-central correction. Thus even a cohomological torus parameter passing all integral Gamma tests need not give an integral ambient realization off $\Lambda_{\rm alg}$.

\begin{corollary}[The critical interval in Harder's coordinate] \label{cor:harder-interval}
Assume the local hypotheses of Theorem~\ref{thm:complex} at a complex place, or of Theorem~\ref{thm:real-criterion} at a real place. Call a regular admissible twist \emph{balanced at $v$} if its Kostant representatives satisfy $\ell(u^\eta)+\ell(u^{\bar\eta})=d_P$ at a complex place, respectively $\ell(u)=d_P/2$ at a real place. A regular admissible twist is balanced precisely when it is simultaneously critical.

At a complex place write $b_{\rm H}^{\xi}=a_0^\xi+f_P+b$ for the two coordinates in \eqref{eq:harder-coordinate}. The balanced twists are precisely those $b\in\Lambda_{\rm alg}$ for which 
\[
 |b_{\rm H}^{\eta}+b_{\rm H}^{\bar\eta}|<
 \min_{\beta\in\Psi_P}\frac{|\Delta_{v,\beta}|}{j(\beta)}.
\]
At a real place with $\mathcal F_P=\varnothing$, put $h_\beta=\langle\widetilde\mu^{(1)},\beta^\vee\rangle$. The shifted central coordinates of all balanced realizations form the set 
\[
 \left\{b_{\rm H}\in a_0+f_P+\Lambda_{\rm alg}:\quad
 |b_{\rm H}|<\min_{\beta\in\Psi_P}\frac{|h_\beta|}{j(\beta)}\right\}.
\]
If $\mathcal F_P\ne\varnothing$, simultaneous criticality requires $b_{\rm H}=0$. It occurs there exactly when the twist is admissible, every fixed line is odd, and all two-element orbits have nonzero gap. The shifted weight at that point has no ambient Kostant realization. 
\end{corollary}
\begin{proof}
The coefficient of $\gamma_P$ in the central projection of $\mu^\eta+w_M\mu^{\bar\eta}$ is $a_0^\eta+a_0^{\bar\eta}$. Hence the normalized sum $T_v(b)$ in \eqref{eq:root-centering} equals $b_{\rm H}^{\eta}+b_{\rm H}^{\bar\eta}$. Apply Proposition~\ref{prop:interval}.

At a real place opposition purity on the Levi root span gives $w_M\widetilde\mu^{(1)}=-\widetilde\mu^{(1)}$. The two coordinates in an orbit are $j(\beta)b_{\rm H}+h_\beta$ and $j(\beta)b_{\rm H}-h_\beta$. Their signs are opposite precisely when $j(\beta)|b_{\rm H}|<|h_\beta|$. On a fixed root $h_\beta=0$, so Theorem~\ref{thm:real-criterion} gives the last assertion. 
\end{proof}

\subsection{General-linear and exceptional normalizing tuples} \label{sec:complex-example}
The $\GL_N$ specialization has a single grade. For $G_2$, a higher grade can restrict the common critical interval. In both settings the root coordinates also give the archimedean Euler quotient at each admissible critical shift. 

\begin{corollary}[Raghuram's general-linear combinatorial lemma]\label{cor:Raghuram}
Let $F$ be totally imaginary, $G=\GL_N$, $N=n+n'$, and let $P$ have Levi $M=\GL_n\times\GL_{n'}$. Take cohomological inducing data $\sigma_1\boxtimes\sigma_2$ satisfying the local hypotheses of Theorem~\ref{thm:complex} at every infinite place, with integral coefficient family $\mu=\mu_1\oplus\mu_2$. At $b=0$, the points $-N/2$ and $1-N/2$ are both critical for $L(s,\sigma_1\times\sigma_2^\vee)$ if and only if $\mu$ has a dominant integral ambient realization with $\ell(u^\eta)+\ell(u^{\bar\eta})=nn'$ at every complex place \cite[Lemma~3.16]{Raghuram2025}. 
\end{corollary}

\begin{proof}
The $nn'$ normal roots $e_i-e_{n+t}$, with $1\le i\le n$ and $1\le t\le n'$, all have grade one. Thus
\[
 r_1=\Std_n\otimes\Std_{n'}^\vee,\qquad
 L_v(s,\sigma_{1,v}\boxtimes\sigma_{2,v},r_1)
 =L_v(s,\sigma_{1,v}\times\sigma_{2,v}^\vee).
\]
Since $\rho_P=(N/2)\gamma_P$, we have $k_P=-N/2$. The assertion therefore follows from Theorem~\ref{thm:complex}, applied at each complex place.

For the cohomological cuspidal data in \cite[Lemma~3.16]{Raghuram2025}, strong purity gives local opposition-purity on both Levi factors. The archimedean components are irreducible, essentially tempered, and cohomological \cite[\S2.5, Prop.~2.13]{Raghuram2025}, so the local hypotheses used above hold at every complex place.
\end{proof}

In this single grade, the cuspidal width defined in \cite[(3.11)]{Raghuram2025} is $\min_v R_v=\min_{v,\beta}|\Delta_{v,\beta}|$. At $b=0$ the criterion reads $|\mathrm{w}_1-\mathrm{w}_2+N|\le\min_v R_v-2$ for the two global purity weights. The gap of two follows from integrality and parity. 

For split $G_2$, let $\alpha$ be short and $\beta$ long, and label each maximal parabolic by its retained Levi simple root. Both Levis are isomorphic to $\GL_2$. The graded representations and evaluation points are given by the following table \cite[\S\S3.1--3.2]{HosseiniJafariG2}.
\begin{equation}\label{eq:g2-tuple}
\begin{array}{c|c|c|c|c}
 P&j&r_j&d_j&jk_P,\ jk_P+1\\ \hline
 P_\beta&1&\Ad^3&4&-5/2,\ -3/2\\
        &2&\det&1&-5,\ -4\\ \hline
 P_\alpha&1&\Std&2&-3/2,\ -1/2\\
         &2&\det&1&-3,\ -2\\
         &3&\Std\otimes\det&2&-9/2,\ -7/2
\end{array}
\end{equation}
Here $\Ad^3=\Sym^3\otimes\det^{-1}$. At coefficient shift $b$, both arguments in grade $j$ decrease by $jb$. 

For these parabolics, Theorem~\ref{thm:complex} recovers the equivalence between adjacent criticality and complementary Kostant degrees established in \cite[Lemmas~6.1--6.2]{HosseiniJafariG2}. The conjugate Kostant lengths sum to five, with
\[
 \ell_j(u^\eta)+\ell_j(u^{\bar\eta})=d_j.
\]
Thus the graded sums are $(4,1)$ for $P_\beta$ and $(2,1,2)$ for $P_\alpha$. Proposition~\ref{prop:interval} determines the simultaneous critical twists by intersecting the restrictions from all grades. 

\begin{example}[Grade three determines the interval]\label{ex:g2-trivial}
For $P_\alpha$, take 
\[
 \mu^\eta=(4,-3),\qquad
 \mu^{\bar\eta}=(4,-4)
\]
in Dynkin coordinates ordered by $(\alpha,\beta)$. The interval is centered at zero, and the grade widths are
\[
 R_{v,1}=4,\qquad R_{v,2}=1,\qquad R_{v,3}=2/3.
\]
Hence
\[
 \mathcal I_v=(-1/3,1/3),
 \qquad
 \mathcal I_v\cap\Lambda_{\rm alg}=\{0\}.
\]
The coefficient-shift formulation also incorporates the Tate-twist variation considered in \cite[\S6.3]{HosseiniJafariG2}.
\end{example}

A direct consequence of the complex local test is a normalization check. For a critical lattice point put 
\[
 h_{v,\beta}(b)=-\min\bigl(q^\eta_\beta(b),
                         q^{\bar\eta}_{w_M\beta}(b)\bigr)
                \in\Z_{>0}.
\]
The Gamma recurrence gives
\begin{equation}\label{eq:Gamma-quotient}
 \prod_j\frac{L_v(js,\sigma_v,r_j)}{L_v(js+1,\sigma_v,r_j)}
     =\frac{(2\pi)^{d_P}}
             {\displaystyle\prod_{\beta\in\Psi_P}h_{v,\beta}(b)}.
\end{equation}
Indeed $\Gamma_\C(t+1)=t\Gamma_\C(t)/(2\pi)$, and the dictionary has one factor for each normal root. This computes the archimedean Euler quotient entering the period convention of Remark~\ref{rem:rational-graph}.

\section{Global balance and common boundary degrees}\label{sec:cohomological-degrees}
A single degree equality can detect criticality at every archimedean place. The reason is that a common relative center forces all nonzero local length defects to have the same sign. The same count measures how much of the boundary degree range remains when the two Kostant lengths differ.

\subsection{The defect formula and the global combinatorial lemma}

After restriction of scalars from $F$ to $\Q$, a Kostant family $w=(w^\eta)_\eta$ is an element of the product Weyl group over the embeddings of $F$. Its length is the sum of the lengths of its components, so  
\[
 \ell(w)=\sum_{\eta:F\hookrightarrow\C}\ell(w^\eta),\qquad
 d_{P,F}=[F:\Q]d_P.
\]
A complex place counts twice. The associate family of Proposition~\ref{prop:associate-complement} satisfies $\ell(w')=d_{P,F}-\ell(w)$.

Define the signed length defect
\begin{equation}\label{eq:defect}
 \Def(w)=\ell(w)-\ell(w')=2\ell(w)-d_{P,F}.
\end{equation}
Reversing source and target changes its sign. Its parity is $\Def(w)\equiv d_{P,F}\pmod2$.

At a parameter $b$ for which every $x^\eta(b)$ is $G$-regular, let $w^\eta=w^\eta(b)\in W^P$ be the unique representative such that $(w^\eta)^{-1}x^\eta(b)$ is strictly $G$-dominant. At a real place $v$, write $w^v=w^\eta$ for its unique embedding. Put $\Def(b)=\Def(w(b))$ and define
\begin{equation}\label{eq:local-defect}
 \Def_v=\begin{cases}
  2\ell(w^v)-d_P,&v\text{ real},\\
  2(\ell(w^\eta)+\ell(w^{\bar\eta}))-2d_P,&v\text{ complex}.
 \end{cases}
 \qquad\Def=\sum_{v\mid\infty}\Def_v .
\end{equation}
The local quantities are computable by counting negative coordinates. Call the family globally balanced when $\Def=0$, and locally balanced at $v$ when $\Def_v=0$; the latter agrees with Corollary~\ref{cor:harder-interval}.

\begin{proposition}[Global defect at a regular shift] \label{prop:global-defect}
Assume that the coefficient family is locally opposition-pure with $a_v=a$ at every archimedean place. Let $b\in\mathbb R$ be regular at every embedding and put 
\[
 T(b)=a+2f_P+2b=a-2s,\qquad s=k_P-b.
\]
For a real place $v$, let $m_v(b)$ count the two-element $w_M$-orbits in $\Psi_P$ whose coordinates have the same sign, and put $n_{\rm fix,v}=n_{\rm fix}(P)$. For a complex place $v=\{\eta,\bar\eta\}$, let $m_v(b)$ count the pairs $(q_\beta^\eta(b),q_{w_M\beta}^{\bar\eta}(b))$, $\beta\in\Psi_P$, whose coordinates have the same sign, and put $n_{\rm fix,v}=0$. Set $m(b)=\sum_{v\mid\infty}m_v(b)$ and let $r_1(F)$ be the number of real places. Then
\begin{equation}\label{eq:regular-defect-count}
\begin{aligned}
  \Def_v(b) &= -\operatorname{sgn}(T(b))\bigl(n_{\rm fix,v}+2m_v(b)\bigr),\\
  \Def(b)   &= -\operatorname{sgn}(T(b))\bigl(r_1(F)n_{\rm fix}(P)+2m(b)\bigr).
\end{aligned}
\end{equation}
Here $\operatorname{sgn}(0)=0$. If $T(b)=0$, regularity forces $r_1(F)n_{\rm fix}(P)=m(b)=0$. Consequently, the local defects are either all nonnegative or all nonpositive. 
\end{proposition}

\begin{proof}
This is the quantitative form of Lemma~\ref{lem:common-sign}. By Proposition~\ref{prop:root-centering}, every paired sum is $j(\beta)T(b)$, with $j(\beta)>0$; a fixed real coordinate is $j(\beta)T(b)/2$. In the sign count \eqref{eq:local-defect}, an opposite-sign pair contributes zero. A same-sign pair contributes $-2\operatorname{sgn}(T(b))$, and a fixed coordinate contributes $-\operatorname{sgn}(T(b))$. This proves the formulas when $T(b)\ne0$. If $T(b)=0$, regularity excludes fixed coordinates and makes every remaining pair opposite in sign.
\end{proof}

\begin{theorem}[Global combinatorial lemma at regular shifts]\label{thm:global-converse}\label{thm:intro-global}
Let $P$ be maximal, let $b\in\Lambda_{\mathrm{alg}}$ be regular at every embedding, and suppose the coefficient family is locally opposition-pure with $a_v=a$ at every archimedean place. The following conditions are equivalent: 
\begin{enumerate}
\item $\Def(b)=0$, equivalently $\ell(w)=\ell(w')$;
\item $\Def_v(b)=0$ at every archimedean place;
\item every pair $(q_\beta^\eta(b),q_{w_M\beta}^{\bar\eta}(b))$ has opposite signs at a complex place; at a real place $\mathcal F_P=\varnothing$ and every two-element $w_M$-orbit has opposite signs.
\end{enumerate}
Suppose moreover that every inducing representation $\sigma_v$ is irreducible, essentially tempered, and cohomological for the stated weights, with the normalizations of the real and complex criteria (Theorems~\ref{thm:real-criterion} and~\ref{thm:complex}). Then these conditions are equivalent to simultaneous criticality at every archimedean place, for either twist family of Remark~\ref{rem:twist-families}: the absolute-value family, with $\sigma_v$ fixed, is tested at $s=k_P-b$, and the algebraic family, with $\sigma_v\otimes\theta_b^{-1}$, at $k_P$.

If both chosen Levi cohomology blocks have degree support $\{k_0\}$ (Definition~\ref{def:retained-Levi-block}), the conditions are also equivalent to overlap of their boundary terms in a common degree.
\end{theorem}

\begin{proof}
Proposition~\ref{prop:global-defect} shows that the local defects are either all nonnegative or all nonpositive, and that their sum vanishes exactly when every pair has opposite signs and no real fixed root occurs.

At a complex place the Gamma factors depend only on the root coordinates, which are the same for both families, and the criticality equivalence is Theorem~\ref{thm:complex}. At a real place, suppose first that $\mathcal F_P\ne\varnothing$. Then condition~(3) fails, and the fixed-line criterion in Theorem~\ref{thm:real-criterion} requires $q_\beta(b)=0$ for $\beta\in\mathcal F_P$ in either family, which a regular $b$ excludes; so neither family is simultaneously critical. If $\mathcal F_P=\varnothing$, the two families differ only in fixed-line signs, so their tests agree and reduce to part~(a) of Theorem~\ref{thm:real-criterion}.

The two boundary terms in the last assertion occur in degrees $k_0+\ell(w)$ and $k_0+\ell(w')$, respectively.
\end{proof}

\begin{corollary}[The degree defect counts archimedean poles]
\label{cor:graded-defect-poles}
Fix irreducible, essentially tempered cohomological inducing data at every archimedean place. Let \(b\in\Lambda_{\mathrm{alg}}\) be regular at every embedding and put \(s=k_P-b\). For each nonempty grade \(j\) define its contribution to the degree defect by 
\[
 \Def_j(b)=2\sum_{\eta:F\hookrightarrow\C}\ell_j(w^\eta)
                  -[F:\Q]d_j .
\]
In \eqref{eq:four-factors}, let \(p_j^{\mathrm{dir}}\) be the sum of the pole orders of the first two factors, evaluated at \(js,js+1\), over all archimedean places. Let \(p_j^{\mathrm{dual}}\) be the corresponding sum for the dual factors at \(1-js,-js\).
A finite factor contributes zero. Then
\[
 \Def_j(b)=p_j^{\mathrm{dual}}-p_j^{\mathrm{dir}},
 \qquad \Def(b)=\sum_j\Def_j(b).
\]
If the relative central coordinate \(a_v\) is independent of \(v\), then \(\lvert\Def_j(b)\rvert=p_j^{\mathrm{dual}}+p_j^{\mathrm{dir}}\) in each grade, and \(\lvert\Def(b)\rvert=\sum_j\lvert\Def_j(b)\rvert\). When the relative centers agree, \(js\) and \(js+1\) are both critical
in grade \(j\) at every archimedean place if and only if \(\Def_j(b)=0\).
\end{corollary}

\begin{proof}
By Proposition~\ref{prop:associate-complement}, passage to the associate parabolic complements the inversions in each grade. Thus \(\Def_j(b)\) is the difference of the two grade-\(j\) inversion counts. In the dictionaries of Propositions~\ref{prop:complex-dictionary} and~\ref{prop:real-dictionary}, an opposite-sign pair contributes zero to this difference and has four finite Gamma factors. Two positive coordinates contribute \(-2\) and give two simple poles on the direct side. Two negative coordinates contribute \(+2\) and give two simple poles on the dual side.

For a fixed real root \(\beta\), regularity gives \(q_\beta(b)\ne0\). If \(q_\beta(b)>0\), the root contributes \(-1\) to the defect. Exactly one of the direct arguments \(-q_\beta(b)+\epsilon\) and \(1-q_\beta(b)+\epsilon\) is a nonpositive even integer. If \(q_\beta(b)<0\), the root contributes \(+1\), and exactly one of the dual arguments \(1+q_\beta(b)+\epsilon\) and \(q_\beta(b)+\epsilon\) gives a pole. In either case the two factors on the other side are finite (Lemma~\ref{lem:parity}). Since Gamma factors have no zeros, the pole orders add. This proves the signed identities.

Under the common-center hypothesis, every paired sum is \(jT(b)\) and every fixed coordinate is \(jT(b)/2\). Hence all pole contributions, if any, occur on the same side across all places and grades. The absolute-value identities and the criticality assertion follow.
\end{proof}

\subsection{The common boundary degree range}

\begin{definition}[Retained cohomology block] \label{def:retained-Levi-block}
Fix an integral ambient coefficient $\lambda$ and the Kostant coefficient $\mu_w=w\cdot\lambda$, embedding by embedding. Choose a compatible boundary level $K_{M,f}$ and put 
\[
 S^M_{K_{M,f}}=M(F)\backslash M(\mathbb A_F)/
 (A_MK_{M,\infty}K_{M,f}),\qquad
 H_!^\bullet=\operatorname{im}(H_c^\bullet\longrightarrow H^\bullet),
\]
where $A_M$ is the positive rational split center of $\operatorname{Res}_{F/\Q}M$. A retained Levi block is a specified finite-dimensional graded subspace 
\[
 H_{M,w}^\bullet\subseteq
 H_!^\bullet(S^M_{K_{M,f}},\tcM_{\mu_w}^M),\qquad
 \mathcal D_w=\{k\in\Z:H_{M,w}^k\ne0\}.
\]
All Hecke selections, finite-level vectors, multiplicities and component projections are included in this choice. If a component character is retained on the connected-compact cover, use the corresponding cover and compatible boundary projection throughout. Make the corresponding choices for $w'$ on the associate Levi.
\end{definition}

For example, a full cuspidal $\sigma$-block is obtained from
\[
 \operatorname{Mult}_{\rm cusp}(\sigma)\otimes
 H^\bullet(\mathfrak m_\infty,K_{M,\infty}A_M;
       \sigma_\infty\otimes\cM_{\mu_w}^M)
 \otimes\sigma_f^{K_{M,f}},
\]
with $K_{M,\infty}^{\circ}$ in place of $K_{M,\infty}$ when working on that cover. The combined action of $A_M$ is trivial. The degree arguments concern these actual cohomology spaces; they do not supply an automorphic realization of a weight family.

Kostant's theorem \eqref{eq:Kostant} gives the selected boundary term 
\[
 \mathcal B_w^q=H_{M,w}^{q-\ell(w)}.
\]

For $H=G,M$, write $H_\infty=\prod_{v\mid\infty}H(F_v)$ and let $A_H$ be the positive rational split center of $\operatorname{Res}_{F/\Q}H$. Put  
\begin{equation}\label{eq:reduced-degree-invariants}
 D_H=\dim_{\R}(H_\infty/(K_{H,\infty}A_H)),\qquad
 \delta_H=\ell_0(H_\infty)-\dim A_H.
\end{equation}
Use $q_0(H)=q_0(H_\infty)$ in the convention of Theorem~\ref{thm:fixed-roots-fundamental-rank}, extended additively over archimedean places. Removing $A_H$ subtracts its dimension from both the symmetric-space dimension and $\ell_0$, so 
\begin{equation}\label{eq:reduced-degree-dimension}
 D_H=2q_0(H)+\delta_H.
\end{equation}
In particular, the value of $q_0(H)$ is unchanged.

\begin{corollary}[Common boundary degrees at a regular twist] \label{cor:regular-boundary-window}
Under the coefficient hypotheses of Theorem~\ref{thm:global-converse}, let $m(b)$ be the number of same-sign pairs in \eqref{eq:regular-defect-count}. If both retained Levi blocks have degree support 
\[
 [q_0(M),q_0(M)+\delta_M]\cap\Z,
\]
their common boundary degrees are exactly 
\begin{equation}\label{eq:regular-boundary-degrees}
 q_0(G)+m(b)\ \le q\le\ q_0(G)+\delta_G-1-m(b),
 \qquad q\in\Z.
\end{equation}
This interval is nonempty precisely when $2m(b)<\delta_G$. If the supports are merely contained in the stated Levi interval, \eqref{eq:regular-boundary-degrees} is a necessary condition. 
\end{corollary}

\begin{proof}
At real places, \eqref{eq:fixed-root-rank-degree} gives the fundamental-rank difference; at complex places the difference is zero. Since $\dim A_M-\dim A_G=1$, we obtain 
\[
 \delta_M-\delta_G=r_1(F)n_{\rm fix}(P)-1,
 \qquad D_G-D_M=d_{P,F}+1.
\]
Equation~\eqref{eq:reduced-degree-dimension} now gives $2(q_0(G)-q_0(M))=d_{P,F}+r_1(F)n_{\rm fix}(P)$. The two lengths sum to $d_{P,F}$ and differ in absolute value by $r_1(F)n_{\rm fix}(P)+2m(b)$. Intersecting the two translates of the Levi degree interval proves the claim. 
\end{proof}

The same count has an analytic interpretation. At a regular twist a same-sign exchanged pair excludes both endpoints in its grade. If all exchanged pairs of the grade have opposite signs, each fixed line selects one endpoint by Lemma~\ref{lem:parity}. Under this hypothesis an endpoint is critical in the grade exactly when all its fixed lines allow it. If the grade has no fixed lines and every exchanged pair has opposite signs, both endpoints are critical.

\begin{example}[Boundary degree four for the $(3,1)$ parabolic] \label{ex:gl4-unbalanced}
For $\GL_4/\Q$ take $\Pi\boxtimes1$, where $\Pi=\Sym^2\pi_\Delta$ is the unitary Gelbart--Jacquet lift of the level-one cusp form $\Delta$ of weight $12$ \cite{GelbartJacquet1978}. The Levi coefficient of this datum has highest weight $(10,0,-10;0)$. At $b=0$, \eqref{eq:gl-coordinates} gives $(13,2,-9)$: the outer coordinates form an opposite-sign pair and the middle root is fixed by $w_M$. The two associate lengths are $1$ and $2$. With the connected compact subgroup used and the positive split center removed, the Levi cohomology occurs in degrees $2,3$ \cite[\S3.1.5, Prop.~3.4]{HarderRaghuram2020}. Hence the two shifted ranges meet only in degree $4=1+3=2+2$.

There is one grade, with adjacent arguments $-2,-1$. The archimedean factor in the unitary normalization has the form \cite[\S2.4, Prop.~2.4 and Rem.~2.5]{LoefflerWilliams}
\[
 L_\infty(s,\Pi)=\Gamma_\C(s+11)\Gamma_\R(s+1).
\]
The lower argument is critical, whereas at the upper argument the real factor has a simple pole. The functional equation reflects $-2,-1$ to $3,2$, where the Euler product converges absolutely and is nonzero. Consequently, 
\[
 L_{\rm fin}(-2,\Pi)\ne0,\qquad
 \operatorname{ord}_{s=-1}L_{\rm fin}(s,\Pi)=1.
\]
The degree overlap therefore accompanies a critical value and a simple trivial zero. An interpretation of the derivative in cohomology requires a further construction. 
\end{example}

\subsection{The common center for automorphic data}

The defect formula requires only that the relative coordinates $a_v$ agree. The central character of a global cohomological datum supplies this equality under the following compatibility hypothesis.

\begin{lemma}[The center of a global cohomological datum] \label{lem:global-center}
Let $M$ be split over $F$, and let $\sigma$ be an irreducible automorphic representation of $M(\mathbb A_F)$ with central character $\omega_\sigma$. Suppose that every $\sigma_v$, $v\mid\infty$, is essentially unitary and cohomological for the coefficient weights $(\mu^\xi)_{\xi\mid v}$. Put $Z=Z_M^\circ$, and write $\omega_{\mu,\infty}$ for the character by which $Z(F_\infty)$ acts on the coefficient representation. Assume the connected central compatibility
\begin{equation}\label{eq:connected-central-compatibility}
 \left.\omega_{\sigma,\infty}\right|_{Z(F_\infty)^\circ}
 =\left.\omega_{\mu,\infty}^{-1}\right|_{Z(F_\infty)^\circ}.
\end{equation}
Then $\mu$ is locally opposition-pure and there is one $c\in\mathfrak a_{M,\mathbb Q}^*$ such that 
\[
 \mu^{\bar\eta}+w_M\mu^\eta=c
 \qquad(\eta\in\Sigma_F).
\]
In particular, $a_v=\langle c,\alpha_P^\vee\rangle$ is independent of $v$. 
\end{lemma}

\begin{proof}
Local opposition-purity follows from Lemma~\ref{lem:essential-unitarity-purity}. For an $F$-rational cocharacter $y:\mathbb G_m\to Z$, the pullback $\chi_y=\omega_\sigma\circ y$ is a Hecke character with connected infinity-type exponents 
\[
 n_{y,\eta}=-\langle\mu^\eta,y\rangle\in\mathbb Z.
\]
At real places its sign is unrestricted by \eqref{eq:connected-central-compatibility}, but $\chi_y^2$ is an algebraic Hecke character with infinity type $(2n_{y,\eta})_\eta$. Purity of algebraic Hecke characters \cite[\S4.1, Lemma~17]{RaghuramHecke2022} therefore gives an integer $\mathrm{w}_y$ such that 
\begin{equation}\label{eq:central-exponent-purity}
 n_{y,\tau^{-1}\circ\eta}+n_{y,\tau^{-1}\circ\bar\eta}=\mathrm{w}_y
 \qquad(\tau\in\operatorname{Aut}(\mathbb C),\ \eta\in\Sigma_F).
\end{equation}
If $F$ has a real place, the same lemma gives the stronger conclusion that all $n_{y,\eta}$ are equal, and $\mathrm{w}_y=2n_{y,\eta}$.

Taking $\tau$ to be the identity and using that $w_M$ fixes $Z$ yields $\langle c_v,y\rangle=-\mathrm{w}_y$, independently of $v$. Since $M$ is split, the $F$-rational cocharacters span $X_*(Z)\otimes\mathbb Q$. Moreover, restriction induces an isomorphism 
\[
 \mathfrak a_{M,\mathbb Q}^*
 \xrightarrow{\ \sim\ }X^*(Z)\otimes\mathbb Q,
\]
because $Z\to M/M_{\rm der}$ is an isogeny. Thus these pairings determine one common $c$. 
\end{proof}

The lemma applies to cuspidal and discrete automorphic data in the split-central normalization of Section~\ref{sec:purity}, provided \eqref{eq:connected-central-compatibility} holds. Its proof needs no cuspidality. The connected compatibility is imposed before quotienting by the positive split center; real sign characters remain part of the inducing data.

\begin{corollary}[The center under arithmetic conjugation] \label{cor:arithmetic-common-center}
Under the hypotheses of Lemma~\ref{lem:global-center}, suppose in addition that $\mu$ is strongly opposition-pure. Then every arithmetic reindexing has the same center $c$:
\[
 ({}^{\tau}\mu)^{\bar\eta}+w_M({}^{\tau}\mu)^\eta=c
 \qquad(\tau\in\operatorname{Aut}(\mathbb C),\ \eta\in\Sigma_F).
\]
In particular, this holds when the conjugate coefficient families satisfy the realization hypothesis of Corollary~\ref{cor:realized-strong-purity}. 
\end{corollary}

\begin{proof}
Strong opposition-purity makes each displayed sum Levi-central. Pairing it with $y$ gives $-\mathrm{w}_y$ by \eqref{eq:central-exponent-purity}. The restriction isomorphism in the proof of Lemma~\ref{lem:global-center} identifies it with $c$. 
\end{proof}

Under the hypotheses of Corollary~\ref{cor:arithmetic-common-center}, for the common central coefficient shift $\mu_b^\eta=\mu^\eta+\theta_b$ of Definition~\ref{def:admissible-shift}, the center becomes $c_b=c+2\theta_b$ at every place and on every arithmetic reindexing. Its relative coordinate is $a+2b$; this assertion concerns the coefficient family and does not assert an automorphic realization of every shift.

\begin{example}[Local defects can cancel without a common center] \label{ex:defect-cancellation}
Let $F$ be a quartic CM field with complex places $v_1,v_2$, and let $\eta_i,\bar\eta_i$ be the embeddings above $v_i$. Take $G=\GL_4$, $M=\GL_3\times\GL_1$, $b=0$, and 
\[
 \mu^{\eta_1}=\mu^{\bar\eta_1}=(2,1,0;0),\qquad
 \mu^{\eta_2}=\mu^{\bar\eta_2}=(2,1,0;6).
\]
Here $f_P=2$, $d_P=3$, and all normal roots $\beta_i=e_i-e_4$ have grade one. Since $w_M$ reverses the first three coordinates, 
\[
 c_{v_1}=(2,2,2;0),\qquad c_{v_2}=(2,2,2;12),
 \qquad a_{v_1}=2,\quad a_{v_2}=-10.
\]
The family is opposition-pure. This condition is not vacuous on the $\GL_3$ factor: for equal conjugate weights $(u,v,w)$ it requires $u+w=2v$.

The normal coordinates are $(5,3,1)$ at both embeddings above $v_1$ and $(-1,-3,-5)$ above $v_2$. All four shifted weights are regular, with Kostant representatives the identity and $w_P$, respectively, and ambient highest weights $(2,1,0,0)$ and $(3,3,2,1)$. The opposed pairs sum to $T_{v_1}(0)=6$ and $T_{v_2}(0)=-6$, as in Proposition~\ref{prop:root-centering}. Thus 
\[
 \Def_{v_1}=-6,\qquad \Def_{v_2}=6,\qquad \Def=0,
\]
although neither place is balanced.

These weights have essentially tempered cohomological local realizations. Let $\mathbb J$ denote the essentially tempered representation of $\GL_3(\C)$ specified in \cite[(2.8)--(2.10), Prop.~2.13]{Raghuram2025} for the coefficient pair $((2,1,0),(2,1,0))$. Take 
\[
 \sigma_{v_1}=\mathbb J\boxtimes\mathbf1,\qquad
 \sigma_{v_2}=\mathbb J\boxtimes|\cdot|_{\C}^{-6}.
\]
No automorphic representation of $M(\mathbb A_F)$ has these local components: its $\GL_1$ factor would have conjugate exponent sums $0$ and $-12$, contrary to the purity of algebraic Hecke characters \cite[Lemma~17(ii)]{RaghuramHecke2022}. Thus opposition-purity centers the root pairs within each place, while a common relative center is needed to prevent cancellation of local defects across places.
\end{example}

\subsection{Selected degree supports and the induced cohomology}

For a matched inducing module realizing this Levi block, relative Shapiro induction may retain a split Lie algebra $\mathfrak a_{\rm rel}$ after quotienting by $A_G$. If its dimension is $d$ and its action on the combined inducing module and coefficient is trivial, the term is
\begin{equation}\label{eq:induced-piece}
 \mathcal J_w^q=\bigoplus_{e=0}^d
 H_{M,w}^{q-\ell(w)-e}\otimes
 \Lambda^e\mathfrak a_{{\rm rel},\C}^*.
\end{equation}
This is the exterior factor in relative Shapiro induction \cite[Ch.~III, Prop.~2.5(ii) and Thm.~3.3(ii)]{BorelWallach}. The formula is used after component projection only when the projection preserves this tensor factor; otherwise use the support of the projected complex.

If the residual action is a nonzero scalar functional $\zeta$, the differential is $\zeta\wedge{}$, and contraction with an element $a$ satisfying $\zeta(a)=1$ contracts the complex. For a maximal rational parabolic, $d=1$ when precisely $\operatorname{Lie}(A_M/A_G)$ remains.

For finite sets $X,Y\subset\Z$, write $X-Y=\{x-y:x\in X,y\in Y\}$ and $X+Y=\{x+y:x\in X,y\in Y\}$.

\begin{lemma}[Degree overlap]\label{lem:degree-overlap}
For the selected graded terms above:
\begin{enumerate}
\item The boundary terms $\mathcal B_w^\bullet$ and $\mathcal B_{w'}^\bullet$ have a common nonzero degree if and only if
\begin{equation}\label{eq:boundary-window}
 \Def(w)\in\mathcal D_{w'}-\mathcal D_w.
\end{equation}
\item If both induced terms have the residual exterior factor of dimension $d$ in \eqref{eq:induced-piece}, they have a common nonzero degree if and only if
\begin{equation}\label{eq:induced-window}
 \Def(w)\in
 (\mathcal D_{w'}-\mathcal D_w)+\{-d,-d+1,\ldots,d\}.
\end{equation}
\end{enumerate}
Consequently these are necessary degree conditions for a nonzero degree-preserving comparison between the corresponding realized terms.
\end{lemma}

\begin{proof}
The boundary degrees coincide exactly when $\ell(w)+k=\ell(w')+k'$ with $k\in\mathcal D_w$ and $k'\in\mathcal D_{w'}$, or $\Def(w)=k'-k$. In \eqref{eq:induced-piece} an additional difference $e'-e$ occurs. Each integer from $-d$ to $d$ is realized by two nonzero exterior powers, which gives \eqref{eq:induced-window}.
\end{proof}

For the same singleton support $\{k_0\}$ on both sides, the boundary test becomes $\Def=0$, while induced overlap for $d=1$ allows $\Def=0,\pm1$ with its prescribed parity. If both Levi supports lie in \([k_{\min},k_{\max}]\), boundary overlap requires
\(\lvert\Def\rvert\le k_{\max}-k_{\min}\). For induced overlap the
upper bound increases by \(d\). For example, the support \(\{0,4\}\)
has difference set \(\{-4,0,4\}\).

In a common degree one can ask whether the intertwiner acts nontrivially and which classes lift to ordinary cohomology. The possible differences between a target degree and a source degree form the set
\((\mathcal D_{w'}-\mathcal D_w)-\Def(w)\).
Constructing a map in one of these degrees is a further step.

\subsection{Changes in boundary degree at a wall}
For a fixed embedding the Kostant length gives the degree shift in the boundary term. We describe its jump at a wall using the notation of Section~\ref{sec:chamber-intervals}.

At a crossing of a single wall $H_\beta$, let $w_-$ and $w_+$ denote the representatives before and after the crossing as $b$ increases. The chambers share a facet in $H_\beta$, and only $q_\beta$ changes sign. Thus $w_+=s_\beta w_-$ and $\ell(w_+)=\ell(w_-)-1$. Distinct crossings therefore give a saturated descending Bruhat chain. Within one grade, $b_{\beta'}<b_\beta$ whenever $\beta'^\vee-\beta^\vee$ is a nonzero nonnegative sum of Levi simple coroots, by strict $M$-dominance. For the relation to reflection orders see \cite[\S2]{Dyer1993}.

For a multiple crossing, let $x\in X\otimes_{\mathbb Z}\mathbb R$ be strictly $M$-dominant, and put 
\[
 N_-(x)=\{\beta\in\Psi_P:\langle x,\beta^\vee\rangle<0\},\qquad
 Z(x)=\{\beta\in\Psi_P:\langle x,\beta^\vee\rangle=0\}.
\]
The vanishing roots form a root subsystem $\Phi_x$ with $\Phi_x^+:=\Phi_x\cap\Phi_G^+=Z(x)$, since strict $M$-dominance excludes Levi roots. Let $W_x$ be its Weyl group, $w_{0,x}$ its longest element relative to $\Phi_x^+$, and $\mathcal W_P(x):=\{w\in W:x\in\overline{wC}\}$ the set of incident representatives, where $\overline{wC}$ denotes the closure of the chamber $wC$.

\begin{proposition}[Chambers at a singular parameter]\label{prop:face}
Every $w\in\mathcal W_P(x)$ belongs to $W^P$ and satisfies 
\[
 N_-(x)\subseteq\Inv_P(w)\subseteq N_-(x)\cup Z(x).
\]
The group $W_x$ acts simply transitively on the incident chambers. There are unique representatives $w_{\min}$ and $w_{\max}$ of minimum and maximum length, with 
\[
 \Inv_P(w_{\min})=N_-(x),\qquad
 \Inv_P(w_{\max})=N_-(x)\cup Z(x),\qquad
 w_{\max}=w_{0,x}w_{\min}.
\]
If $x=x(b_*)$ is singular, then for sufficiently small $\varepsilon>0$, 
\[
 w(b_*-\varepsilon)=w_{\max},\qquad
 w(b_*+\varepsilon)=w_{\min}.
\]
The degree drop at $b_*$ is $|Z(x)|$, with graded contributions $\ell_j(w_{\max})-\ell_j(w_{\min})=|Z(x)\cap\Psi_{P,j}|$. 
\end{proposition}

\begin{proof}
Small perturbations of $x$ preserve all nonzero root pairings and strict $M$-dominance, which gives membership in $W^P$ and the inclusions. Locally at $x$ only the hyperplanes of $\Phi_x$ occur, so facet crossings connect the incident chambers through reflections in $W_x$; the action is free because $W$ acts simply transitively on chambers \cite[\S\S1.8, 1.12]{Humphreys1990}.

For $\beta\in Z(x)$ the pairing $\langle\gamma_P,\beta^\vee\rangle=j(\beta)$ is positive, so $x\pm\varepsilon\gamma_P$ makes every vanishing pairing positive, respectively negative, and changes no other sign. These perturbations give the two extremal chambers; every other incident chamber has strictly intermediate length. The element $w_{0,x}$ fixes $x$ and exchanges the positive and negative chambers of $\Phi_x$, whence $w_{\max}=w_{0,x}w_{\min}$. For $x=x(b_*)$ the
perturbations are $x(b_*\pm\varepsilon)$.
\end{proof}

The drop counts positive vanishing roots, not the rank of $\Phi_x$: an $A_2$ subsystem contributes three. If $b_*\in\Lambda_{\rm alg}$ is a crossing parameter, then $\mu_{b_*}+\rho=x(b_*)+z_{b_*}$ is singular, and Proposition~\ref{prop:inverse-kostant} excludes a
realization $\mu_{b_*}=w\cdot\lambda$ with $\lambda\in X_G^+$.

\subsection{Boundary degrees in classical and exceptional examples}
For $\GL_{n+n'}/\Q$ put $\mathsf b_n=\lfloor n^2/4\rfloor$. When $nn'$ is even, $\mathsf b_{n+n'}=\mathsf b_n+\mathsf b_{n'}+nn'/2$. Thus a balanced representative carries the Levi bottom degree to the ambient bottom degree \cite[Ch.~3]{HarderRaghuram2020}.

\begin{example}[The complex $G_2$ bottom degree]
Let $F$ be totally imaginary with $r$ complex places. For either maximal $G_2$ parabolic, $M\simeq\GL_2$ and $d_{P,F}=10r$. For a full tempered cuspidal Levi block, each semisimple complex factor contributes degrees one and two, and the center modulo the positive rational split center contributes an exterior range of width $r-1$. The full Levi range is therefore $[r,3r-1]\cap\Z$ (see \cite[Prop.~2.15]{Raghuram2025}).

On the critical branch $\ell(w)=\ell(w')=5r$, so the Levi bottom degree gives the degree-$6r$ comparison in \cite{HosseiniJafariG2}. For the full range, \eqref{eq:regular-boundary-degrees} gives 
\[
 6r+m(b)\le q\le8r-1-m(b).
\]
An overlap with unequal lengths requires $1\le m(b)\le r-1$ and occurs above the bottom degree.
\end{example}

\begin{remark}[Harder's genus-three calculation]
At level one for $\GSp_6$, the middle parabolic $P_2$ stabilizes an isotropic two-plane. The two elliptic cuspidal inputs have total Levi degree two. The Kostant lengths four and three therefore give boundary degrees six and five.

The two grades of the dual nilradical have dimensions six and one. The second grade contributes $\zeta(n_i+1)/\zeta(n_i+2)$, with $n_i>0$ even in the case considered \cite[\S9.6.3, pp.~459--461, (9.235), (9.244)--(9.248), (9.252)]{HarderChapterNine}. The positive odd numerator argument is noncritical because the dual factor $\Gamma_\R(-n_i)$ has a pole, although the finite zeta value is nonzero. The positive even denominator argument is critical. For a noncritical endpoint with a trivial zero, see Section~\ref{sec:g2-real-endpoints}.
\end{remark}

\begin{example}[Klingen boundary degrees] \label{ex:klingen-mismatch}
Let $F$ be totally real of degree $r$, $n\ge2$, and $P$ the Klingen parabolic of $\Sp_{2n}$. There is one $w_M$-fixed normal root and $d_P=2n-1$. If the other orbit gaps are nonzero, the two chambers adjacent to the central wall have lengths $n$ and $n-1$ by Section~\ref{sec:real-places}. For global cohomological inducing data satisfying Lemma~\ref{lem:global-center}, this is a common wall $b^\circ$ at all real places. In a sufficiently small common neighborhood, its lower side has
\[
 \ell(w)=rn,\qquad \ell(w')=r(n-1),\qquad\Def(w)=r;
\]
the upper side has defect $-r$. An admissible lattice point on either side need not lie in this neighborhood. Proposition~\ref{prop:klingen-wall} realizes the common singular wall $b^\circ=0$ by globally generic cuspidal Levi data for suitable weights.

For these centered data, the fixed root has coordinate $q_{2e_1}(b)=b$ and odd sign at the wall. At a real place, let $b$ be an admissible shift in either adjacent chamber, so $0<|b|<\min_{i\ge2}t_i$. Each exchanged pair then permits both endpoints. The fixed line selects one by Lemma~\ref{lem:parity}. In the algebraic family its sign exponent is $\epsilon\equiv1+b\pmod2$, so $q_{2e_1}(b)-\epsilon$ is odd throughout. The upper endpoint $k_P+1$ is critical below the wall and the lower endpoint $k_P$ above it. For the fixed representation the sign exponent stays equal to one. Among its two arguments $s,s+1$, the critical endpoint occupies the same position as in the algebraic family for even $b$ and the opposite position for odd $b$ (Remark~\ref{rem:twist-families}).

At $-1$ in the $\GL_1$ factor, the shifted coefficient acquires the sign $(-1)^b$. The inverse algebraic twist cancels this sign for every $b$, and the absolute-value twist satisfies the same full-compact parity condition exactly for even $b$ (Remark~\ref{rem:compact-conventions}). For $n=2$, the inducing exponent in the companion paper is $-b=s+2$ in the length-two chamber \cite{HosseiniJafariKlingen}.

Farther from the common wall, the common-center sign argument gives $|\Def|\ge r$. For example, for Klingen in $\Sp_4/\Q$, the line through $\mu=(0;0)$ at $b=-6$ has $x=(-4,1)$. All three normal coordinates are negative, giving $\Def=3$.

Suppose the retained Levi block has symplectic discrete-series cohomology in one degree and retains the center in degrees $0,\ldots,r-1$. Its width $r-1$ is smaller than $|\Def|$, so the two boundary ranges are disjoint. In the neighboring chamber with $\Def=r$, retaining one residual split direction in induced cohomology permits  
\[
 rn+t_0=r(n-1)+(t_0+r-1)+1,
\]
where $t_0$ is the common Levi bottom degree and the final $1$ is the residual exterior degree in \eqref{eq:induced-piece}. For $\Sp_4/\Q$, the reduced defect $\delta_G$ is zero, so \eqref{eq:regular-boundary-degrees} is empty for both maximal parabolics. For Klingen, a cuspidal Levi class in degree one gives the boundary degrees $3$ and $2$. 
\end{example}

\begin{remark}[Denominators require an integral comparison]\label{rem:denominators}
For even $n\ge2$, the denominator of Harder's level-one Eisenstein class for $\mathrm{SL}_2(\Z)$ is the positive numerator of $\zeta(-n-1)$ \cite[Thm.~2.13]{BekkiSakamoto2025}. Here $w_M$ fixes the unique positive root, the two Kostant lengths are $0$ and $1$, and of the adjacent zeta arguments $-n-1$ and $-n$, only the first is critical.

Harder's Siegel calculation relates a secondary cohomological construction to a modular $L$-value ratio and $\zeta'(-n_1)/\pi$, and asks how the denominator of the period-normalized value is reflected in the Eisenstein class \cite[\S\S2.2.1--2.2.3]{HarderSiegel2013}, \cite[\S9.6.2, pp.~453--455]{HarderChapterNine}. For a separate conjectural denominator ideal and a proposed height pairing involving a central derivative, see \cite[(93), (112)]{HarderDenominators}.

The Klingen parabolic of $\GSp_4$ restricts to the Klingen parabolic of $\Sp_4$. Example~\ref{ex:klingen-mismatch} identifies the upper endpoint as critical in the algebraic family in the length-two chamber. 

In a companion manuscript in preparation \cite{HosseiniJafariKlingen}, the author conjectures that, for the primitive degree-three boundary class in this chamber, the period-normalized twisted adjoint $L$-value at that endpoint generates the denominator ideal of the rational Eisenstein lift away from the bad primes. This integral comparison requires a choice of lattices. All results of the present paper are independent of the conjecture.

\end{remark}

\subsection{Inner Levi classes and global Hecke cohomology} \label{sec:clozel-boundary}
Clozel extends Scholze's construction \cite[\S5.2]{Scholze2015} to inner Levi classes with algebraic coefficients \cite[Thms.~3.4, 3.7 and Prop.~3.5]{Clozel2022}. Each resulting boundary class either lifts to ordinary cohomology of $G$ or has nonzero connecting image in compactly supported cohomology. The formulation below records the generalized eigenspaces and their multiplicities in the Hecke normalization of Section~\ref{sec:root-setup}.

For this subsection let $G/\Q$ be connected reductive and let $P=MU_P$ be a maximal proper rational parabolic. Groups over number fields are treated by restriction of scalars. Choose compatible maximal compact subgroups $K_{G,\infty}$ and $K_{M,\infty}$, and a product level $K_f=\prod_pK_p$ that is neat in the strong sense used in \cite[\S2.1]{Clozel2022}. Put 
\[
 K_{P,f}=K_f\cap P(\A_f),\qquad
 K_{M,f}=\operatorname{pr}_M(K_{P,f}),
\]
and set 
\[
 S^G_{K_f}=G(\Q)\backslash G(\A)/(A_GK_{G,\infty}K_f),\qquad
 S^M_{K_{M,f}}=M(\Q)\backslash M(\A)/(A_MK_{M,\infty}K_{M,f}).
\]
Here $A_H$ is the identity component of the real points of the maximal $\Q$-split central torus of $H$. The Borel--Serre compactification $\overline S^G_{K_f}$ \cite{BorelSerre1973} has an open boundary piece
\[
 Y_P=P(\Q)\backslash P(\A)/(A_MK_{M,\infty}K_{P,f}),\qquad
 j:Y_P\hookrightarrow\partial\overline S^G_{K_f}.
\]
The full stratum of type $P$ is a finite disjoint union of such pieces; the projection $Y_P\to S^M_{K_{M,f}}$ has compact nilmanifold fibers. Let $\tcM$ be the local system of an algebraic $G$-module $\cM$ over $\overline\Q$. All cohomology below is over $\overline\Q$, and $H_!$ denotes the image of compactly supported in ordinary cohomology.

Choose a finite set $S$ so that outside $S$ the levels are the chosen compatible hyperspecial subgroups and $G(\Q_p)=P(\Q_p)K_p$. In the split setting use the pinned model of Section~\ref{sec:factor-conventions}. Give each integral subgroup volume one. The spherical Hecke transfer is the unnormalized constant term 
\[
 C_P^{\mathrm a}:\mathcal H_G^S\longrightarrow\mathcal H_M^S,
 \qquad C_P^{\mathrm a}(f)(m)=\int_{U_P(\A_f^S)}f(mu)\,du,
\]
the map $\lambda\circ\rho$ of \cite[\S2.2]{Clozel2022}. For a character $\chi$ of $\mathcal H_M^S$, put $\chi'=\chi\circ C_P^{\mathrm a}$. Subscripts $\chi$ and $\chi'$ denote generalized eigenspaces. The theorem below uses the full compact quotients; a component projection can be retained when it is compatible with its maps.

\begin{theorem}[Clozel] \label{thm:clozel-boundary}
Let $\cM_\mu^M$ be a fixed Kostant summand of $H^{q_U}(\mathfrak u_P,\cM)$, put $q=k+q_U$, and suppose
\[
 A_\chi:=H_!^k(S^M_{K_{M,f}},\tcM_\mu^M)_\chi\ne0.
\]
Through the nilmanifold-fibration decomposition, regard $A_\chi$ as a subspace of $H_!^q(Y_P,\tcM)$. Every nonzero $\alpha\in A_\chi$ has an extension $\gamma\in H^q(\partial\overline S^G_{K_f},\tcM)_{\chi'}$ with $j^*\gamma=\alpha$. For each such extension, exactly one of the following holds:
\begin{enumerate}
\item $\gamma=\operatorname{res}\beta$ for some nonzero $\beta\in H^q(S^G_{K_f},\tcM)_{\chi'}$;
\item the connecting image $\partial\gamma\in H_c^{q+1}(S^G_{K_f},\tcM)_{\chi'}$ is nonzero.
\end{enumerate}
Moreover,
\begin{equation}\label{eq:clozel-multiplicity}
 \dim H^q(S^G_{K_f},\tcM)_{\chi'}+
 \dim H_c^{q+1}(S^G_{K_f},\tcM)_{\chi'}
 \geq\dim A_\chi.
\end{equation}
These statements follow from \cite[Thms.~3.4, 3.7 and Prop.~3.5]{Clozel2022}; no essential-temperedness or regularity of the highest weight is required. 
\end{theorem}

\begin{proof}
Clozel's Hecke-equivariant form of the Harder--Schwermer fibration decomposition applies to both ordinary and compactly supported cohomology \cite[Lemmas~3.1--3.3, \S3.4]{Clozel2022}. Passing to the image gives 
\[
 H_!^q(Y_P,\tcM)
 \simeq\bigoplus_{q_M+q_U=q}
 H_!^{q_M}(S^M_{K_{M,f}},\widetilde{H^{q_U}(\mathfrak u_P,\cM)}).
\]
Extension by zero defines $j_*:H_c^q(Y_P,\tcM)\to H^q(\partial\overline S^G_{K_f},\tcM)$. Its image $\mathcal E$ maps surjectively to $H_!^q(Y_P,\tcM)$ under $j^*$, with Hecke action on $H_!^q(Y_P,\tcM)$ through $C_P^{\mathrm a}$. Taking generalized eigenspaces is exact for finite-dimensional modules over a commutative algebra. Thus every specified $\alpha$ has a lift in $\mathcal E_{\chi'}$.

The Hecke-equivariant exact sequence
\begin{equation}\label{eq:clozel-exact-sequence}
 H_c^q(S^G_{K_f},\tcM)\longrightarrow H^q(S^G_{K_f},\tcM)
 \xrightarrow{\operatorname{res}}
 H^q(\partial\overline S^G_{K_f},\tcM)
 \xrightarrow{\partial}H_c^{q+1}(S^G_{K_f},\tcM)
\end{equation}
remains exact on the $\chi'$-generalized eigenspaces. For a fixed $\gamma$, the two cases are precisely $\partial\gamma=0$ and $\partial\gamma\ne0$; any ordinary lift is nonzero since $j^*\gamma=\alpha\ne0$.

Projection of $j^*$ to the chosen Kostant and $\chi$-summand gives a surjection $\mathcal E_{\chi'}\twoheadrightarrow A_\chi$. The dimension of the boundary generalized eigenspace is at most the sum on the left of \eqref{eq:clozel-multiplicity}, by exactness, proving the bound.
\end{proof}

The alternatives concern each chosen extension $\gamma$; different extensions of the same $\alpha$ can have different connecting images. They do not assert that exactly one of the two global cohomology spaces is nonzero.

For an automorphic constituent use the convention
\[
 H^\bullet(\mathfrak m,K_{M,\infty}A_M;
            \sigma_\infty\otimes\cM_\mu^M)\ne0,
\]
with the split-central character matched as above. The transfer $C_P^{\mathrm a}$ is the Hecke transfer for unnormalized induction \cite[\S2]{Clozel2022}. In the split maximal-parabolic setting of Section~\ref{sec:root-setup}, this is the finite-part Hecke normalization of $I_P(k_P,\sigma)$ by \eqref{eq:induction-conversion}. A coefficient twist is accompanied by the inverse representation twist; no dual coefficient is inserted into the Kostant decomposition.

For a nonmaximal rational face, Poincar\'e duality defines a Gysin map $j_*$ that raises degree by $\dim(A_M/A_G)-1$. With constant complex coefficients, Clozel deduces $j^*j_*=0$ from Saper's deformation theorem \cite[\S5.1, Prop.~5.1, Thm.~5.2]{Clozel2022}. Thus every class in the image of this Gysin map restricts to zero on the face.

For the separate question of a common sign in higher relative rank, see the $\GL_3$ Borel example in Section~\ref{sec:complex-criterion}.

\begin{lemma}[Strictly dominant coefficients and inner degrees] \label{lem:regular-inner-degrees} 
Let $\lambda$ be a strictly dominant integral ambient highest weight and let $\mu=w\cdot\lambda$ be a Kostant coefficient. With the invariants of \eqref{eq:reduced-degree-invariants},
\[
 H_!^k(S^M_{K_{M,f}},\tcM_\mu^M)=0
 \quad\text{unless}\quad
 q_0(M)\le k\le q_0(M)+\delta_M.
\]
For restriction of scalars, strict dominance of the ambient highest weight is required at every embedding.
\end{lemma}

\begin{proof}
For a simple Levi coroot $\alpha^\vee$, the coroot $w^{-1}\alpha^\vee$ is positive because $w\in W^P$. Every simple-coroot pairing of $\lambda+\rho$ is at least two; therefore
\[
 \langle\mu,\alpha^\vee\rangle
 =\langle\lambda+\rho,w^{-1}\alpha^\vee\rangle-1\ge1.
\]
Thus both $\mu$ and its contragredient highest weight are strictly Levi dominant.

Li and Schwermer prove that for regular coefficients ordinary cohomology vanishes in degrees below $q_0(M)$ \cite{LiSchwermer2004}. Saper proves the equal-rank case and explains how his method also gives the general bound \cite[Thm.~5 and n.~4, pp.~332--333]{Saper2005}.

Poincar\'e duality for the dual coefficient gives vanishing of compactly supported cohomology above $D_M-q_0(M)=q_0(M)+\delta_M$. One may pass to an orientable neat finite cover and descend in characteristic zero. Taking the image of compactly supported in ordinary cohomology gives the asserted interval. 
\end{proof}

For a strictly dominant ambient highest weight the lemma bounds both Levi supports in Corollary~\ref{cor:regular-boundary-window}. The reflection $q\mapsto D_G-1-q$ exchanges the two endpoints of \eqref{eq:regular-boundary-degrees}, since $D_G=2q_0(G)+\delta_G$. For any degree in this interval the dual obstruction degree therefore also lies at or above $q_0(G)$. The ordinary lower-degree vanishing bound alone therefore does not decide which branch of Theorem~\ref{thm:clozel-boundary} occurs.

\section{Eisenstein classes for the half-spin parabolics of split \texorpdfstring{$E_6$}{E6}}\label{sec:realclass-E6}

The two end parabolics of split $E_6/\Q$ are the exceptional entries in Theorem~\ref{thm:realclass-classification}. For these parabolics the regular branch of the real combinatorial lemma leads to actual Eisenstein classes. Use Bourbaki numbering, take $P=P_1$, and put
\[
 u=s_1s_3s_4s_5s_2s_4s_3s_1,\qquad
 \lambda_t=2t\rho,\qquad \mu_t=u\cdot\lambda_t
 \quad(t\in\Z_{\ge0}).
\]
The representative $u$ has length eight. At the centered shift $b=0$, the cohomological evaluation point is $s=-6$. For a unitary inducing datum $\pi$, write $\sigma_t=\pi\otimes|\xi|_{\mathbb A}$, where $\xi$ is the Levi character of weight $6\gamma_1$. Then
\[
 I_P(-6,\sigma_t)=I_P(0,\pi).
\]
The centered Eisenstein series is therefore evaluated on the unitary axis.

\begin{theorem}[Eisenstein classes at the centered twist] \label{thm:E6-ordinary}
Let $G/\Q$ be split semisimple of type $E_6$ (either isogeny form), $P=P_1$, and $t\ge1$. Use the level and rational structures of Theorem~\ref{thm:clozel-boundary}. Write $A_\chi= H_!^k(S^M_{K_{M,f}},\tcM_{\mu_t}^M)_\chi$, $q=k+8$ and $\chi'=\chi\circ C_P^{\mathrm a}$. Through the $u$-Kostant summand, regard $A_\chi$ as a subspace of $H^q(Y_P,\tcM_{\lambda_t})$, as in Theorem~\ref{thm:clozel-boundary}. Then
\begin{equation}\label{eq:E6-face-image}
 A_\chi\ \subseteq\
 \operatorname{im}\!\left[H^q(S^G_{K_f},\tcM_{\lambda_t})_{\chi'}
         \xrightarrow{\ j^*\operatorname{res}\ }H^q(Y_P,\tcM_{\lambda_t})\right]
 \quad\text{over }\overline{\Q}.
\end{equation}
In particular $\dim H^q(S^G_{K_f},\tcM_{\lambda_t})_{\chi'}\ge\dim A_\chi$. For each $t$, at a common sufficiently small finite level, both $k=12$ and $k=13$ have nonzero inner Levi cohomology. In each of degrees $20$ and $21$, the value of a holomorphic Eisenstein family at the centered parameter gives a cohomology class with nonzero restriction to its inducing face. The Hecke characters may differ. The diagram involution gives the corresponding statement for $P_6$, transporting the level, weights and Hecke character.
\end{theorem}

The length-eight shift carries Levi classes of degrees $12$ and $13$ to ambient degrees $20$ and $21$. In degree $20$, one cuspidal inducing datum gives an Eisenstein class with nonzero restrictions to both associate face types (Corollary~\ref{cor:E6-two-faces}). The root calculation gives the full critical interval. At its center, a twisted Lefschetz number supplies the nonvanishing needed for cuspidal Levi cohomology, and holomorphic evaluation of the Eisenstein series lifts the Levi classes.

\subsection{Root coordinates and the critical interval}
Write weights by their simple-coroot pairings, in the numbering of Appendix~\ref{app:fixed-roots}. The algebraic twist lattice is $\Z$ for the simply connected group and $3\Z$ for the adjoint group. Indeed,
\[
 \gamma_1=\frac{4\alpha_1+3\alpha_2+5\alpha_3+6\alpha_4
                    +4\alpha_5+2\alpha_6}{3},
\]
so $\gamma_1$ has order three modulo the root lattice. The coefficient constructed below lies in the root lattice, so its centered member is valid for both isogeny forms. Set
\[
 \mu=(-8,0,0,0,0,4),\qquad
 x=\mu+\rho=(-7,1,1,1,1,5).
\]
Applying the defining reflections for $u$ to $\rho=(1,1,1,1,1,1)$ in the order $1,3,4,2,5,4,3,1$ gives $x=u\rho$. Thus $\mu=u\cdot0\in\Z\Phi_G$. The positive Levi pairings put $u$ in $W^P$, so $\mu$ is the highest weight of a Kostant summand for the trivial ambient coefficient.

The eight roots below have sixth coefficient zero. Their images under $w_M$, given by \eqref{eq:realclass-E6-tau}, are the other eight normal roots.
\[
\begin{array}{c|r@{\qquad}c|r}
 \beta&q_\beta&\beta&q_\beta\\\hline
 (1,0,0,0,0,0)&-7&(1,1,1,1,0,0)&-4\\
 (1,0,1,0,0,0)&-6&(1,1,1,1,1,0)&-3\\
 (1,0,1,1,0,0)&-5&(1,1,1,2,1,0)&-2\\
 (1,0,1,1,1,0)&-4&(1,1,2,2,1,0)&-1
\end{array}
\]
The reflection rule $s_i(x)_j=x_j-C_{ji}x_i$ gives the displayed coordinates. By \eqref{eq:realclass-E6-tau}, $w_M x=-x$, so the other eight coordinates have the opposite signs and the same absolute values. Thus $d_P=16$ and $\ell(u)=8$. All normal roots have coroot grade one, and none is fixed by $w_M$.

The sum of the normal roots is $12\gamma_1$, so $f_P=6$. Also $\mu+w_M\mu=-12\gamma_1$. Projection to the Levi center gives $a_0=-6$. In Harder's convention, $b_{\rm H}=b$, so the centered shift $b=0$ corresponds to $s_0=-6$.

On the simply connected cover of the derived dual Levi of type $D_5$, the action on the dual nilradical is one $16$-dimensional half-spin representation \cite[\S2.5.4]{Kim2005}. The action of the full dual Levi on the dual nilradical specifies the central character of the half-spin factor $r_1$. For simply connected $E_6$, the corrected description of the full dual Levi is given in \cite[\S5, p.~293]{KimKim2011}. The orientation of $r_1$ is fixed in Remark~\ref{rem:factor-conventions}. In increasing order with multiplicity, the positive root coordinates are $(h_1,\ldots,h_8)=(1,2,3,4,4,5,6,7)$.

\begin{proposition}[Critical twists of the half-spin family]\label{prop:E6-regular-family}
Write $m_t=2t+1$. The shifted weight and the Levi coefficient satisfy
\[
 x_t=\mu_t+\rho=m_t u\rho,\qquad
 \mu_t=x_t-\rho=(-14t-8,2t,2t,2t,2t,10t+4).
\]
Both $\lambda_t$ and $\mu_t$ belong to the root lattice, $\mu_t=u\cdot\lambda_t$, and $\cM_{\mu_t}^M$ occurs in $H^8(\mathfrak u_P,\cM_{\lambda_t})$. For $t\ge1$ both highest weights are strictly dominant for their respective groups.

In the representation convention of Section~\ref{sec:root-setup}, fix an irreducible essentially tempered real representation $\sigma_{t,\infty}$ cohomological with coefficient $\cM_{\mu_t}^M$. For $t\ge1$ essential temperedness follows already from essential unitarity and cohomology, with the central character normalized as in Lemma~\ref{lem:alignment}. For the fixed inducing representation, among algebraically admissible coefficient shifts $b\in\Lambda_{\rm alg}$, those for which the half-spin Artin factor is critical at both $s_b=-6-b$ and $s_b+1$ are exactly
\[
 \mathcal I_t=[-2t,2t]\cap\Lambda_{\rm alg}.
\]
Since there are no fixed normal roots, the algebraic twist family has the same critical set, by Theorem~\ref{thm:real-criterion}. The number of critical shifts is $4t+1$ on the simply connected form and $2\lfloor 2t/3\rfloor+1$ on the adjoint form. For each such $b$ the shifted coefficient is obtained from the same $u$ with ambient weight
\[
 \lambda_{t,b}=(2t-b,2t,2t,2t,2t,2t+b).
\]
For $t\ge1$ precisely the twists $|b|<2t$ in $\mathcal I_t$ give strictly dominant ambient highest weights, hence $4t-1$ on the simply connected form.

With $h_1,\ldots,h_8$ as above, the archimedean quotient is
\begin{equation}\label{eq:E6-regular-ratio}
 \frac{L_\infty(s_b,\sigma_{t,\infty},r_1)}{L_\infty(s_b+1,\sigma_{t,\infty},r_1)}
 =\frac{(2\pi)^8}{\prod_{i=1}^8(m_th_i-b)}.
\end{equation}
In particular its centered value is $4\pi^8/(315m_t^8)$.
\end{proposition}

\begin{proof}
Both $2\rho$ and $u\rho-\rho$ are in the root lattice. The displayed coordinates and $\ell(u)=8$ give the assertion by Kostant's theorem. Since $w_M u\rho=-u\rho$, one has $\mu_t+w_M\mu_t=-\rho-w_M\rho=-12\gamma_1$. Thus $a_0=-6$ and $b_{\rm H}=b$ throughout the family. All grades equal one, so after shifting the coefficient, the coordinates in an orbit are $b+m_t h_i,b-m_t h_i$ for some $1\le i\le 8$. The real pole test is $(b+m_t h_i)(b-m_t h_i)<0$ for all $i$, or $|b|<m_t$. In particular the walls $|b|=m_t$ are excluded. Intersecting with $\Z$ or $3\Z$ gives the two counts.

The chamber is constant for $|b|<m_t$. The reflection word gives $u^{-1}\gamma_1=-\gamma_1+\gamma_6$. Thus $u\cdot\lambda_{t,b}=\mu_t+b\gamma_1$, as asserted. Every $b\in\mathcal I_t$ is a regular shift. For $t\ge1$, its ambient highest weight is strictly dominant exactly when $|b|<2t$. At an admissible endpoint $|b|=2t$ the shift is still regular, while one simple-coroot pairing of $\lambda_{t,b}$ is zero.

Strict dominance of the Levi highest weight for $t\ge1$ makes every unitary cohomological realization tempered modulo center by the Vogan--Zuckerman classification \cite{VoganZuckerman1984}. For the regular-weight argument, see \cite[\S7.3.1]{Grobner2013}. Lemma~\ref{lem:alignment} therefore gives the half-spin factor $\prod_{i=1}^8\Gamma_\C(s+6+m_t h_i)$. Its quotient at adjacent points gives \eqref{eq:E6-regular-ratio}.
\end{proof}

\subsection{Cuspidal Levi cohomology and proof of the restriction theorem}
\begin{lemma}[Nonzero centered Levi cohomology] \label{lem:E6-centered-existence}
Let $t\geq1$ and put $\nu_t=\mu_t+6\gamma_1$. For either isogeny form of $G$, there is a sufficiently small strongly neat ambient level, with induced Levi level $K_{M,f}$, such that
\[
 H_!^{12}(S^M_{K_{M,f}},\tcM_{\mu_t}^M)\ne0,
 \qquad
 H_!^{13}(S^M_{K_{M,f}},\tcM_{\mu_t}^M)\ne0.
\]
\end{lemma}

\begin{proof}
Put $H=M_{\mathrm{der}}$. This is the split simply connected group of type $D_5$: passing from simply connected to adjoint $E_6$ does not change it, since the kernel $\mu_3$ has trivial intersection with its center $\mu_4$. The weight $\nu_t$ has zero central component and is self-dual. With the Levi simple roots ordered as $(6,5,4,3,2)$, its Dynkin coordinates are $(10t+4,2t,2t,2t,2t)$; in the standard orthogonal coordinates for $D_5$ it is
\[
 \nu_t=(16t+4,6t,4t,2t,0).
\]
These coordinates are integral with even sum, so $\nu_t$ lies in the $D_5$ root lattice. In particular $F_t:=\cM_{\nu_t}^M|_H$ is a self-dual algebraic coefficient.

To obtain cuspidal cohomology for $H$, first consider the local Lefschetz number. Let $\vartheta$ be the rational Chevalley Cartan involution of $H$. Choose a fundamental Cartan and a $\vartheta$-stable Borel in $\mathfrak h_{\mathbb C}$. In the corresponding $D_5$ coordinates, $\vartheta$ negates the
fifth coordinate and fixes the other four. Thus $\nu_t$ is $\vartheta$-invariant and vanishes on the split direction of this Cartan. It defines a unitary character of the Cartan, admissible for this Borel by dominance.

The Vogan--Zuckerman construction and cohomology formula give an irreducible unitary representation $A_{\mathfrak b}(\nu_t)$ with nonzero $F_t$-cohomology \cite[Thms.~5.3 and~5.5]{VoganZuckerman1984}. For the Cartan Levi, the hypothesis on noncompact simple factors in \cite[Prop.~6.5]{VoganZuckerman1984} is vacuous. The even Dynkin coordinates put the highest weight in the lattice specified in \cite[Def.~3.1.1]{RohlfsSpeh1989}. The Rohlfs--Speh Lefschetz formula \cite[Prop.~4.3]{RohlfsSpeh1989} gives
\[
 \bigl|\operatorname{Lef}
       (\vartheta,A_{\mathfrak b}(\nu_t)\otimes F_t)\bigr|
 =\sum_q\dim H^q(\mathfrak h,K_H;
                  A_{\mathfrak b}(\nu_t)\otimes F_t)>0,
\]
where $K_H=H(\mathbb R)\cap K_{M,\infty}$ is connected.

The group $H$ is absolutely almost simple and has positive $\mathbb Q$-rank. The Borel--Labesse--Schwermer nonvanishing theorem, applied with $S=\{\infty\}$ and coefficient $F_t$, therefore gives nonzero cuspidal cohomology for some congruence subgroup $\Gamma_H\subset H(\mathbb Q)$ \cite[Thm.~10.4]{BorelLabesseSchwermer1996}. This theorem applies to nontrivial coefficients through the displayed local Lefschetz condition. The passage from the enlarged set $S$ containing a finite place back to $S=\{\infty\}$ is explained in \cite[Thm.~6.5 and proof of Thm.~10.4]{BorelLabesseSchwermer1996}.

To pass to the full Levi quotient, let $\zeta_M:M\longrightarrow M/H\simeq\mathbb G_m$ be a primitive quotient character. Choose the finite level so that
\[
 \Gamma':=H(\mathbb Q)\cap K_{M,f}\subseteq\Gamma_H,
 \qquad
 \zeta_M(K_{M,p})\subseteq1+p\mathbb Z_p
\]
at one odd prime $p$. These conditions persist on shrinking the ambient level. For every arithmetic stabilizer $\Gamma_g=M(\mathbb Q)\cap gK_{M,f}g^{-1}$, compactness of the local character images gives $\zeta_M(\Gamma_g)\subseteq\{\pm1\}$; the condition at $p$ therefore gives $\Gamma_g\subset H(\mathbb Q)$.

Let $A_M$ be the connected split center at infinity. The character $\zeta_M$ is primitive on the split maximal torus, so it maps $M(\mathbb R)$ onto $\mathbb R^\times$. It maps $A_M$ onto $\mathbb R_{>0}$, and Cartan decomposition gives $\zeta_M(K_{M,\infty})=\{\pm1\}$. The restriction of $\zeta_M$ to $K_{M,\infty}$ has connected kernel
$K_H=H(\R)\cap K_{M,\infty}$, so $K_{M,\infty}$ has exactly two
components. Here $K_{M,\infty}=K_{G,\infty}\cap M(\R)$ is the full
compact subgroup fixed in Section~\ref{sec:factor-conventions}.
Both $S^M_{K_{M,f}}$ and the induced cohomology below use this full
compact subgroup, as in Remark~\ref{rem:compact-conventions}. It follows that
\[
 M(\mathbb R)/(K_{M,\infty}A_M)
       \simeq H(\mathbb R)/K_H.
\]
Indeed $M(\mathbb R)=H(\mathbb R)K_{M,\infty}A_M$ and $H(\mathbb R)\cap K_{M,\infty}A_M=K_H$. Consequently the component indexed by the identity finite double coset is $\Gamma'\backslash H(\mathbb R)/K_H$.

Cuspidal cohomology injects into inner cohomology. Pullback from $\Gamma_H$ to its finite-index subgroup $\Gamma'$ is injective, by the finite-cover trace, and preserves inner cohomology. Extension by zero to the other open-and-closed components gives a nonzero inner class on $S^M_{K_{M,f}}$. On each component $\tcM_{\mu_t}^M$ agrees with the self-dual coefficient $F_t$, because $\mu_t-\nu_t=-6\gamma_1$ vanishes on $H$.

The group $H(\mathbb R)$ is connected \cite[Cor.~4.7]{BorelTits1972}, so these neat components are orientable manifolds of dimension $25$. The reduced defect is $\delta_M=1$, since the complexified derived group has rank five and its maximal compact subgroup has rank four. Lemma~\ref{lem:regular-inner-degrees} therefore restricts inner cohomology to degrees $12$ and $13$. Poincar\'e duality for the self-dual coefficient pairs these two inner spaces perfectly: the maps from compactly supported to ordinary cohomology are transposes under duality. Nonvanishing in one degree therefore gives nonvanishing in both.
\end{proof}

\begin{proof}[Proof of Theorem~\ref{thm:E6-ordinary}]
Fix the level in the theorem and take an arbitrary class $\alpha\in A_\chi\otimes_{\overline\Q}\C$ in the $u$-Kostant summand. We construct a lift with $P_1$-restriction $\alpha$.

\emph{Cuspidal decomposition on the Levi.}
On either isogeny form, $\xi$ is the square of the character of weight $3\gamma_1$. Thus $\xi=|\xi|$ on $M(\R)$ and $\xi$ is trivial on each arithmetic stabilizer, since its value there is a positive rational unit. Center the coefficient by $\nu_t=\mu_t+6\gamma_1$. It has trivial split-central weight by Proposition~\ref{prop:E6-regular-family}. For unitary data $\pi$ with this coefficient, the twist defining $\sigma_t$ gives an isomorphism respecting the full compact subgroup, 
\[
 \sigma_{t,\infty}\otimes\cM_{\mu_t}^M
   \simeq \pi_\infty\otimes\cM_{\nu_t}^M.
\]

Over $\C$, cuspidal cohomology lies in inner cohomology, which lies in the image of discrete $L^2$ cohomology \cite[\S1.3.3]{Urban2011}. The Levi highest weight $\nu_t$ is strictly dominant because $t\ge1$. The regular-weight argument used above shows that every cohomological discrete constituent for this coefficient is tempered at infinity modulo the center \cite[\S7.3.1]{Grobner2013}. Wallach's cuspidality criterion on the semisimple arithmetic components then puts its full discrete multiplicity in the cuspidal spectrum \cite[Thm.~4.3]{Wallach1984}. At the given finite level we may therefore write $\alpha=\sum_l\alpha_l$ as a finite sum of cuspidal classes, retaining their multiplicity spaces and $K_{M,f}$-fixed vectors.

\emph{Lifting the prescribed class.}
For each $\alpha_l$ choose a cuspidal cocycle of type $(\pi_l,u)$. In the finite induced representation take the section having the prescribed value on the double coset defining $Y_P$ and supported on that coset. Such a section exists by $K_{M,f}$-invariance. Since $w_Mx_t=-x_t$, the projection of $x_t=u(\lambda_t+\rho)$ to the $W_M$-fixed central direction is zero. Hence the cohomological evaluation point in $I_P(\cdot,\pi_l)$ is zero (Proposition~\ref{prop:coefficient-spectral-parameters}). On this Kostant summand the relative split line acts trivially. Lift $\alpha_l$ in exterior degree zero of \eqref{eq:induced-piece}, as $\alpha_l\otimes 1$, and denote the inducing cocycle by $\phi_l$. Tangential restriction to the inducing face takes the exterior-degree-zero component, so it recovers $\alpha_l$ in both Levi degrees. In degree $13$ an additional exterior-degree-one component may be present.

The Eisenstein series is holomorphic at zero by Langlands' unitary-axis theorem, in the formulation of \cite[Thm.~7.2(a)]{Arthur2005}. Schwermer's construction gives a closed Eisenstein form $E(\phi_l,0)$ whose tangential restriction is computed by the cuspidal constant term \cite[Thm.~4.11]{Schwermer1983}. For the restriction formula,
see \cite[Thm.~3.1 and \S3.2]{Schwermer2021}. Here $P_1$ has distinct standard associate $P_6$, so $W(A_P,A_P)=\{1\}$. The only term on the inducing face is the identity, and $j^*\operatorname{res}[E(\phi_l,0)]=\alpha_l$. Summing over $l$ lifts $\alpha$.

Restriction is defined over $\overline\Q$, and the unramified Hecke action on the face is through $C_P^{\mathrm a}$ as in Theorem~\ref{thm:clozel-boundary}. Projecting the lift to the $\chi'$-primary space still gives $\alpha$ on the face. Images and primary spaces commute with extension of scalars, which proves \eqref{eq:E6-face-image} over $\overline\Q$.

\emph{Occurrence in both degrees.}
Lemma~\ref{lem:E6-centered-existence} gives a common level with nonzero inner classes in degrees $12$ and $13$. Applying the inclusion to a nonzero Hecke summand in each degree gives the assertions in degrees $20$ and $21$ since $\ell(u)=8$. The pinned diagram involution transports the argument to $P_6$.
\end{proof}

\subsection{Associate faces and archimedean normalization}
\begin{corollary}[The two restrictions at the centered twist] \label{cor:E6-two-faces}
Put $P=P_1$, $Q=P_6$ and $v_P=w_Gw_M$. Let $\Theta$ be the pinned diagram automorphism of $G$ inducing $-w_G$ on weights; use the same symbol for that action. First work over $\C$. For $R=P,Q$, write $Y_{[R]}$ for the union of all open faces of type $R$ at level $K_f$, and $\operatorname{res}_{[R]}$ for tangential restriction to that union. Let $\phi$ be a cuspidal inducing cocycle as in Theorem~\ref{thm:E6-ordinary}, with inducing face class $\alpha$. The two restrictions of $e_\phi=[E(\phi,0)]$ are
\[
  \bigl(\operatorname{res}_{[P]}e_\phi,\operatorname{res}_{[Q]}e_\phi\bigr)
  =\bigl(\alpha,\,[M(v_P,0)\phi]|_{Y_{[Q]}}\bigr).
\]
Here $\alpha$ is extended by zero to the other open pieces of type $P$. The operator $M(v_P,0)$ is the standard global intertwining operator in normalized induction, defined using the fixed rational representative $\dot v_P$ of Section~\ref{sec:factor-conventions}. The target Kostant weight is $\Theta\mu_t$, and its representative has length eight. For a fixed contributing cusp $\pi$, take $k$ least with nonzero reduced Levi cohomology, including the compact component action. In degree $k+8$ the displayed comparison from the full $\pi$-summand to its transported cuspidal summand on $Q$ is an isomorphism. In particular a nonzero $\alpha$ in this least degree gives two nonzero restrictions at any level where it occurs. For every $t\ge1$ this gives an Eisenstein class detected on both face types in degree $20$ at a suitable level.
\end{corollary}

\begin{proof}
The set $W(A_P,A_Q)$ has one element, represented by $v_P$, since $W(A_P,A_P)=1$. The constant-term formula for boundary restriction gives the displayed pair \cite[\S3.2, p.~137]{Schwermer2021}. Since $w_Mx_t=-x_t$, we have $v_Px_t=\Theta x_t$. As $\Theta$ fixes $\lambda_t+\rho$ and $\rho$, the target representative is $\Theta u\Theta^{-1}$, with weight $\Theta\mu_t$ and length eight. The centered Levi coefficient also transports: $v_P$ carries the positive Levi roots to the positive Levi roots of $Q$, so $v_P\nu_t=v_P(x_t-\rho_M)=\Theta x_t-\rho_{M_Q}$.

To see why both face restrictions detect the comparison in degree $k+8$, consider normalized induction $I_P(0,\pi)$. A contributing Kostant row must have the derived-Levi infinitesimal character fixed by $\pi_\infty$. Its shifted Levi-dominant weight thus lies on $x_t+\R\gamma_1$, and has norm $\|x_t\|$ because it is Weyl-conjugate to $\lambda_t+\rho$. In a Weyl-invariant metric $x_t\perp\gamma_1$: the involution $w_M$ negates $x_t$ and fixes $\gamma_1$. Hence only $x_t$ occurs. Transport proves uniqueness on $Q$ and identifies the reduced Levi degree ranges.

In relative Shapiro induction \eqref{eq:induced-piece}, evaluation on the faces omits the relative split tangent direction. Only exterior degree zero contributes in this row in degree $k+8$, since exterior degree one would require reduced Levi degree $k-1$. Face evaluation is therefore an isomorphism on each side. This calculation retains the cuspidal multiplicity spaces and all finite sections; the Levi components act trivially on the relative split line and taking compact component invariants is exact.

By Langlands' functional equation, $M(v_P^{-1},0)$, formed using $\dot v_P^{-1}$, is inverse to $M(v_P,0)$ \cite[Thm.~7.2(a)]{Arthur2005}. These operators intertwine the ambient group actions and preserve $K_f$-invariants, hence induce inverse maps on relative cohomology. The two face evaluations identify the latter maps with the asserted comparison and its inverse.

At the level of Theorem~\ref{thm:E6-ordinary}, choose a cusp contributing in degree $12$ on the Levi. A lower contributing reduced degree would give inner cohomology there with the same finite vector, contrary to Lemma~\ref{lem:regular-inner-degrees}. Thus $k=12$ gives the asserted two restrictions in degree $20$.
\end{proof}

\begin{remark}[Rational Betti classes on the two face types]
At a common level, each of degrees $20$ and $21$ also contains a rational Betti class restricting nontrivially to both face types. Retain the pinned archimedean maximal compact of Section~\ref{sec:factor-conventions} and replace $K_f$ by $K_f\cap\Theta(K_f)$, using the pinned diagram automorphism. The same finite vectors still contribute at this smaller level. The weight $2t\rho$ is $\Theta$-invariant and has a $\Q$-model by split highest-weight theory. In each degree Theorem~\ref{thm:E6-ordinary} gives a nonzero restriction to $Y_{[P]}$; applying $\Theta$ gives a nonzero restriction to $Y_{[Q]}$. With the rational local system, both maps commute with extension of scalars by the finite CW models of the Borel--Serre compactification and its open strata \cite{BorelSerre1973}. Thus each rational restriction has proper kernel. Two proper subspaces cannot cover a vector space over $\Q$, so in each degree a rational class has both restrictions nonzero.

In either degree the rational class may involve several Hecke eigensystems, since \(\Theta\) can carry \(\chi'\) to a different Hecke character. In degree \(20\), Corollary~\ref{cor:E6-two-faces} also gives a complex Eisenstein class from a single cuspidal inducing datum whose restrictions to both face types are nonzero.
\end{remark}

For the same cuspidal datum, consider the local standard intertwining operator at infinity. Use Arthur's local normalization and Haar measures for the tempered representation $\pi_\infty$ at unitary parameter zero \cite[Thm.~2.1 and \S3, (3.2)]{Arthur1989Intertwining}. In the opposite-parabolic direction, the contragredient adjoint representation in the normalizing factor is $r_1$.

For the standard additive character of $\R$, the archimedean $\varepsilon$-factor is $1$. Indeed the eight two-dimensional summands contribute $i^{2m_t h_i+1}$, and $\sum_i h_i=32$. With these measures the normalizing scalar is
\[
 \mathcal R_\infty(t)
 =\frac{L_\infty(0,\pi_\infty,r_1)}
        {L_\infty(1,\pi_\infty,r_1)}
 =\frac{4\pi^8}{315m_t^8}\ne0.
\]
The normalized operator is holomorphic and unitary at zero. Multiplying it by $\mathcal R_\infty(t)$ gives a holomorphic isomorphism $A_\infty(-6,\sigma_{t,\infty})$. The change to the Haar measures fixed in Section~\ref{sec:factor-conventions} multiplies this operator by a nonzero constant.

\begin{remark}[Rational graph and period convention]\label{rem:rational-graph}
For the arithmetic comparison in the sense of Harder and Raghuram \cite[Chs.~6--7]{HarderRaghuram2020}, fix an embedding $\iota:E\hookrightarrow\C$ and two selected boundary Hecke summands $U_E,V_E$ in the same degree. Suppose projection maps the restriction image in $U_E\oplus V_E$ isomorphically onto $U_E$ and that the Hecke-equivariant maps from $U_E$ to $V_E$ form an $E$-line generated by an isomorphism $A_E$. The image is then the graph of $t_E A_E$ for some $t_E\in E$.

Suppose $A_E=B_E\otimes(e_U\mapsto e_V)$ in compatible rational models, where $B_E$ identifies the finite parts and $e_U,e_V$ generate the one-dimensional comparison lines. Write the analytic generators as $p_Ue_U,p_Ve_V$, and let $c_\infty\ne0$ be the scalar of the standard archimedean operator between them. Its Euler-normalized scalar is $c_\infty/\mathcal R_\infty$ when the archimedean Euler quotient $\mathcal R_\infty$ is finite and nonzero.

If the constant term identifies the complexified graph with the standard global operator on these summands, the unramified quotient $\mathcal R^S$ is finite and nonzero, and the remaining finite operators compare with $B_E$ by $\iota(\kappa_E)$, $\kappa_E\in E^\times$, comparison on $e_U,e_V$ gives 
\[
 D_F^{-\dim U_P/2}\,\frac{\mathcal R^S}{\mathcal P_\iota}
 =\iota(t_E/\kappa_E),\qquad
 \mathcal P_\iota=\frac{p_U}{c_\infty p_V}.
\]
Here $D_F=|\operatorname{disc}(F/\Q)|$ and $U_P$ is the unipotent radical. The measure factor uses covolume one globally and root coordinate measures of volume one on $\mathcal O_v$ at finite places, $dx$ at real places, and $2\,dx\,dy$ at complex places \cite[\S5.2.3, (5.12)--(5.14)]{Raghuram2025}. 

For the half-spin family, the two-face result gives the analytic comparison in degree $20$. Applying the displayed identity requires the rational graph and
the finite and period comparisons for that inducing datum.
\end{remark}

\appendix
\section{Fixed normal roots: verification by type}\label{app:fixed-roots}
\begin{lemma}[Fixed-root lemma for maximal parabolics]
\label{lem:fixed-roots-strongly-orthogonal}
Let $\Phi$ be a reduced crystallographic root system, with chosen simple
system $\Delta$, and let $I=\Delta\setminus\{\alpha_k\}$. The positive roots
fixed by the longest element $w_I$ of $W_I$ are pairwise strongly
orthogonal. On the irreducible component containing $\alpha_k$ they are
exactly the roots listed in the proof below.
\end{lemma}

\begin{proof}
Only the irreducible component containing $\alpha_k$ contributes, because
the longest element of every other component has no fixed roots. Thus
assume $\Phi$ irreducible. The proof uses the classification and standard
coordinate models of root systems \cite{Humphreys1990}, with finite
lattice certificates for the exceptional systems.
Every fixed positive root lies outside $\Phi_I$, since
$w_I\Phi_I^+=-\Phi_I^+$.

For $A_{N-1}$, use roots $e_i-e_j$ and the block decomposition
$(k,N-k)$. The element $w_I$ reverses each block. There is one fixed
positive root if both block sizes are odd, joining their middle
coordinates, and none otherwise.

For $B_n$ and $C_n$, $w_I$ reverses $e_1,\ldots,e_k$ and negates all
remaining coordinates. In either case the fixed roots include
\[
 e_i+e_{k+1-i}\quad(1\leq i\leq\lfloor k/2\rfloor).
\]
If $k$ is odd there is one further root, $e_{(k+1)/2}$ for $B_n$ and
$2e_{(k+1)/2}$ for $C_n$. These are all the fixed roots, directly from the
coordinate root lists. Their supports are disjoint, and neither the sum
nor difference of two listed roots is a root.

For $D_n$ with $n\geq4$ and $k\leq n-2$, the first $k$ coordinates are reversed.
On the remaining $n-k$ coordinates the longest element negates every
coordinate when $n-k$ is even, and fixes $e_n$ while negating the others
when $n-k$ is odd. Thus the fixed roots are the mirror sums displayed
above, together with
\[
 e_{(k+1)/2}+e_n,\quad e_{(k+1)/2}-e_n
 \quad\text{if both $k$ and $n-k$ are odd}.
\]
These roots are strongly orthogonal: the two additional roots have sum
and difference twice a coordinate vector, which is not a $D_n$ root,
and all other supports are disjoint. For the end parabolic $P_n$ the
fixed roots are $e_i+e_{n+1-i}$, $1\leq i\leq\lfloor n/2\rfloor$.
The other end parabolic follows by the diagram automorphism
$e_n\mapsto-e_n$.

For the exceptional systems use Bourbaki numbering: the $E_n$ edges are
$1$--$3$--$4$--$5$--$6$--$7$--$8$, truncated as appropriate, together
with $2$--$4$; in $F_4$, roots $\alpha_1,\alpha_2$ are long and
$\alpha_3,\alpha_4$ short; in $G_2$, $\alpha_1$ is short. A digit string
denotes simple-root coefficients.

For the simply laced types we use the following fact. A nonzero
vector $x$ of norm squared two in an $ADE$ root lattice, with
nonnegative simple-root coefficients $c_i$, is a root. Indeed
$\sum_ic_i(x,\alpha_i)=(x,x)=2$, so $(x,\alpha_i)>0$ for some $i$ in
the support of $x$. By Cauchy--Schwarz this integer is at most two,
with equality only for $x=\alpha_i$. Otherwise it equals one, and the
reflection $x-\alpha_i$ has smaller height, nonnegative coefficients
and the same norm; induct on the height.

In $E_6$, a positive root $c=(c_1,\ldots,c_6)$ satisfies
$0\leq c\leq(1,2,2,3,2,1)$ coordinatewise. The fixed-vector equations
for one representative of each diagram-symmetry orbit are
\[
\begin{array}{c|l}
k& w_Ic=c\\\hline
1&c_2+c_3=2c_1,\quad 2c_4=3c_1,\quad c_5=c_1,\quad2c_6=c_1,\\
2&c_1+c_6=c_2,\quad c_3+c_5=2c_2,\quad2c_4=3c_2,\\
3&2c_1=c_3,\quad c_2+c_6=c_3,\quad c_4+c_5=2c_3,\\
4&c_1+c_3=c_4,\quad2c_2=c_4,\quad c_5+c_6=c_4.
\end{array}
\]
For later reference, at $P_1$ the longest Levi element
$w_I=w_M$ acts on the simple roots by
\begin{equation}\label{eq:realclass-E6-tau}
 w_I\alpha_1=(1,2,2,3,2,1),\qquad
 w_I\alpha_2=-\alpha_3,\quad w_I\alpha_3=-\alpha_2,\quad
 w_I\alpha_i=-\alpha_i\quad(i=4,5,6).
\end{equation}
Together with $c_k>0$ and $(c,c)=2$, these give, respectively,
$\varnothing$, $\{122321\}$, $\{112221\}$, and $\{111211\}$.
For example, the third row forces $c_3=2,c_1=1$ and leaves four
possibilities of squared norms $2,4,6,4$; the fourth likewise leaves
four possibilities of squared norms $4,6,2,4$. Each surviving vector is a
root by the fact above. The diagram symmetry gives the cases $k=5,6$.
Every fixed set has at most one member.

For the remaining exceptional types, the last column of each table
lists row indices of maximal minors of the matrix of fixed-root
columns. A single unmarked minor has
determinant $\pm1$. Where two minors are given, their displayed
determinants are coprime. In either case the columns span a primitive
sublattice of the root lattice.
\[
\begin{array}{c|c|l|l}
\Phi&k&\text{fixed roots}&\text{rows}\\\hline
E_7&1&2234321&7\\
&2&0112111,1112110,1122100,1223321&1234\\
&3&1122221,1123211,1223210&124\\
&4&1112111,1112210,1234321&125\\
&5&0111110,1011110,1122211,1224321&1236\\
&6&1123321,1223221&12\\
&7&0112221,1112211,1122111&123\\\hline
E_8&1&22454321,23354321&28\\
&2&12233321,12243221,12343211,22343210&1234:2;\ 1345:3\\
&3&11222221,11232211,12232210,23465321&1245\\
&4&11121111,11122110,12343221,23465421&1256\\
&5&11222111,11222210,13354321,22354321&1367\\
&6&11233211,12232211,23465431&127\\
&7&12244321,12343321,22343221&138\\
&8&23465432&1:2;\ 2:3
\end{array}
\]
For the two exceptional systems with roots of unequal length the
corresponding lists are
\[
\begin{array}{c|c|l|l}
\Phi&k&\text{fixed roots}&\text{rows}\\\hline
F_4&1&2342&1:2;\ 2:3\\
&2&1222,1231&13\\
&3&1121,1342&14\\
&4&1232&1\\\hline
G_2&1&21&2\\
&2&32&1:3;\ 2:2
\end{array}
\]
Fixed spaces are taken in the real span of the root system.
The longest Levi element fixes the omitted fundamental weight and acts
on each component of the deleted diagram by its longest element, whose
fixed space has dimension $\lfloor m/2\rfloor$ for $A_m$, one for $D_m$
with $m$ odd, two for $E_6$, and zero otherwise. Adding these to the line
of the omitted fundamental weight gives
\[
 \dim\operatorname{Fix}(w_I)=
 \begin{cases}
 (1,4,3,3,4,2,3)&E_7,\\
 (2,4,4,4,4,3,3,1)&E_8,\\
 (1,2,2,1)&F_4,\\
 (1,1)&G_2,
 \end{cases}
 \qquad k=1,2,\ldots,
\]
which is the number of columns listed for each $k$. This count does not
use Theorem~\ref{thm:realclass-classification}; it agrees with the
table there. It is not used for $E_6$, where the fixed roots do not
span the fixed space. Simple-reflection
multiplication verifies that the listed columns are fixed.
For $E_7,E_8$ their Gram matrices are $2I$, so each column is a root
by the fact above. The columns thus give strongly orthogonal roots.
They span the fixed space by the dimension count. By primitivity they
form a basis of its integral lattice, whose Gram matrix $2I$ permits
only the signed basis vectors to have norm squared two.

For $F_4$ take $\alpha_1=e_2-e_3$, $\alpha_2=e_3-e_4$, $\alpha_3=e_4$
and $\alpha_4=(e_1-e_2-e_3-e_4)/2$. The listed columns are roots,
with Gram matrices $(2)$, $\operatorname{diag}(2,1)$,
$\operatorname{diag}(1,2)$ and $(1)$, respectively. Again they are
integral bases of the fixed lattices. A vector using both basis
directions has norm squared at least three, and a nontrivial multiple
of absolute coefficient at least two has norm squared at least four.
Since roots have norm squared one or two, only the signed basis vectors
are roots. For $G_2$ each fixed space is a line with the displayed
primitive root as generator. This proves completeness and strong
orthogonality in every type.
\end{proof}

\bibliographystyle{amsplain}
\bibliography{bib}

\end{document}